\documentclass[10pt]{article}
\usepackage[mathscr]{eucal}
\usepackage{xypic}
\usepackage{amsfonts}
\usepackage{amsmath}
\usepackage{amsthm}
\usepackage{amssymb}
\usepackage{latexsym}
\usepackage{eufrak}
\newcommand{\lc}{{}_{c}}
\newcommand{\rc}{_{\!c'}}

\newtheorem{theorem}[equation]{Theorem}
\newtheorem{proposition}[equation]{Proposition}
\newtheorem{corollary}[equation]{Corollary}
\newtheorem{lemma}[equation]{Lemma}
\theoremstyle{definition}                               
\newtheorem{definition}[equation]{Definition}
\newtheorem{remark}[equation]{Remark}
\newtheorem{conjecture}[equation]{Conjecture}

\newcommand{\dis}{{\displaystyle}}

\newcommand{\beq}{\begin{equation}\label}

\newcommand{\iso}{{\;\;\stackrel{_\sim}{\longrightarrow}\;\;}}
\newcommand{\cd}{\!\cdot\!}
\newcommand{\XX}{{\mathbf{x^2}}}
\newcommand{\YY}{{\mathbf{y^2}}}
\newcommand{\vi}{${\sf {(i)}}\;$}
\newcommand{\vii}{${\sf {(ii)}}\;$}
\newcommand{\viii}{${\sf {(iii)}}\;$}
\newcommand{\iv}{${\sf {(iv)}}\;$}
\newcommand{\mr}{|_{\hreg}}
\newcommand{\mrw}{|_{_{\hreg\!\!/W}}}
\newcommand{\Reg}{{\sf Reg}}
\newcommand{\sset}{\subset}
\newcommand{\ssminus}{\smallsetminus}
\newcommand{\G}{\Gamma}

\newcommand{\Mon}{{\mathsf{Mon}}}
\newcommand{\kk}{{\mathscr{K}}_c}
\def\hp{\hphantom{x}}
\def\ccr{\C[R]^W_{_{\sf reg}}}

\def\ccp{\Z[R]^W_+}
\newcommand{\eHe}{\e\hh\e}
\newcommand{\Ch}{{\mathtt{Ch}}}
\newcommand{\BK}{{\mathbb{K}}}
\newcommand{\SH}{{\mathsf{S}}}
\newcommand{\minus}{{_{\!\pmb{\boldsymbol{-}}}}}
\newcommand{\spher}{^{^{\!\mathtt{spher}}}\!^{\!}}
\newcommand{\sh}{{\mathbf{h}}}
\newcommand{\irrep}{{\sf{Irrep}}}
\newcommand{\Rep}{{\sf{Rep}}}
\newcommand{\id}{{{\mathtt {Id}}}}
\newcommand{\Ind}{{{\mathtt {Ind}}}}

\newcommand{\into}{\,\,\hookrightarrow\,\,}
\newcommand{\too}{\,\,\longrightarrow\,\,}
\newcommand{\onto}{\,\,\twoheadrightarrow\,\,}

\newcommand{\Hol}{{\mathcal{H}}ol}
\newcommand{\ad}{{\mathtt{{ad}}^{\,}}}

\newcommand{\Spec}{{\mathtt{Spec}}}
\newcommand{\End}{{\mathtt{End}}}
\newcommand{\rk}{{\mathtt{rk}}}
\newcommand{\supp}{{\mathtt{supp}}}

\newcommand{\Hom}{{\mathtt{Hom}}}

\newcommand{\grd}{{\mathtt{gr}}}

\newcommand{\Tr}{{{\mathsf {Tr}}}}

\newcommand{\hr}{{\mathfrak{h}^{^{\mathsf{reg}}}}}

\newcommand{\dd}{{\mathcal{D}}}
\newcommand{\GL}{\operatorname{GL}}

\newcommand{\cont}{{\mathsf {cont}}}
\newcommand{\pp}{{\mathsf{P}}}

\newcommand{\hh}{{\mathsf{H}}}

\newcommand{\ehe}{{\mathbf{e}\mathsf{H}_c\mathbf{e}}}

\newcommand{\e}{{\mathbf{e}}}
\newcommand{\hcc}{\hh_{c'}}
\newcommand{\hhcc}{(\hh_c\mbox{-}\hh_{c'})}
\newcommand{\ehhcce}{(\e\hh_c\e\mbox{-}\e\hh_{c'}\e)}
\newcommand{\og}{\oo_{{\mathfrak{g}}}}

\newcommand{\ff}{^\flat}
\newcommand{\fx}{^{(x)}\!}

\newcommand{\mm}{{\mathcal{M}}}
\newcommand{\cyc}{{\mathtt{cycl}}}

\newcommand{\Ker}{{\mathtt {Ker}}}

\newcommand{\h}{{\mathfrak{h}}}

\newcommand{\hreg}{{\mathfrak{h}^{^{_{\mathsf{reg}}}}}}

\newcommand{\triv}{{\mathsf{triv}}}
\newcommand{\HC}{{\mathscr{HC}}}

\newcommand{\PPP}{{\mathbb{P}}}

\newcommand{\ddx}{\dd(X_c)}
\def\C{{\mathbb{C}}}

\def\qq{{\mathcal{Q}}}
\def\ms#1{\mathcal{#1}}

\def\Z{{\mathbb{Z}}}

\def\oo{{\mathscr O}}
\def\co{{\mathcal O}}

\def\ll{{\mathcal L}}
\def\LL{{\mathsf{L}}}

\def\eu{{\mathsf{eu}}}

\def\fin{_{_{\,\mathsf{fin}}}}
\def\bb{{\mathcal{B}}}

\def\sll2{{\mathfrak{s}\mathfrak{l}}_2}
\def\cc{{\mathcal C}}
\def\CC{\C[R]^W}

\def\ccirc{{{}_{\,^{^\circ}}}}

\makeatletter
\def\downroundfill{$\m@th \setbox\z@\hbox{$\braceld$}%
  \braceld\leaders\vrule height\ht\z@ depth\z@\hfill\bracerd$}
\def\uproundfill{$\m@th \setbox\z@\hbox{$\braceld$}%
 \bracelu\leaders\vrule height\ht\z@ depth\z@\hfill\braceru$}
\def\overround#1{\mathop{\vbox{\m@th\ialign{##\crcr\noalign{\kern3\p@}
      \downroundfill\crcr\noalign{\kern3\p@\nointerlineskip}
      $\hfil\displaystyle{#1}\hfil$\crcr}}}\limits}
\def\underround#1{\mathop{\vtop{\m@th\ialign{##\crcr
      $\hfil\displaystyle{#1}\hfil$\crcr\noalign{\kern3\p@\nointerlineskip}
      \uproundfill\crcr\noalign{\kern3\p@}}}}\limits}
\makeatother
\begin{document}
\setlength{\parindent}{6mm}
\setlength{\parskip}{3pt plus 5pt minus 0pt}
\centerline{\Large {\bf
CHEREDNIK ALGEBRAS AND DIFFERENTIAL}}
\vskip 2pt
\centerline{\Large {\bf OPERATORS ON QUASI-INVARIANTS}}
\vskip 4mm
\centerline{\large {\sc {Yuri Berest, Pavel Etingof and Victor Ginzburg}}}
%
\vskip 2pt

\begin{abstract}{\footnotesize{\noindent
We develop representation theory of the rational Cherednik algebra $\hh_c$
associated to a finite Coxeter group $W$
in a vector space $ \h $, and a parameter  `$c$'.
We use it to 
show that, for integral values of `$c$', the
algebra $\hh_c$ is simple and Morita equivalent
to $\dd(\h)\#W$, the cross product of $W$ with the algebra of  polynomial 
differential operators on $\h$.

Chalykh, Feigin, and Veselov [CV], [FV], introduced
an  algebra,  $Q_c$,  of
{\it quasi-invariant} polynomials on $\h$,
such that $\C[\h]^W\subset Q_c\subset \C[\h]$.
We prove that 
the algebra $\dd(Q_c)$ of differential operators on
quasi-invariants is a
 simple algebra, Morita equivalent to $\dd(\h)$.
The subalgebra $\dd(Q_c)^W\subset\dd(Q_c) $ of
$W$-invariant operators turns out to be
isomorphic to the spherical subalgebra $\ehe\subset \hh_c$.
We  show that  $\dd(Q_c)$ is 
generated, as an algebra,  by  $Q_c$ and
 its
`Fourier dual' $Q_c\ff$,
and that $\dd(Q_c)$ is a rank one projective
$Q_c\otimes Q_c\ff$-module (via  multiplication-action on $\dd(Q_c)$ on
opposite sides).
}}
\end{abstract}

{\centerline{\bf Table of Contents}
\vskip -5mm
$\hspace{20mm}$ {\footnotesize \parbox[t]{115mm}{\,

\hp${}_{}$\hp1.{ $\;\,\,$} {\tt Introduction} \newline
\hp2.{ $\;\,\,$} {\tt Standard modules over the 
rational Cherednik algebra}\newline
\hp3.{ $\;\,\,$} {\tt Harish-Chandra $\hh_c$-bimodules}\newline
\hp4.{ $\;\,\,$} {\tt The spherical subalgebra $\e\hh_c\e$}\newline
\hp5.{ $\;\,\,$} {\tt A trace on the Cherednik algebra}\newline
\hp6.{ $\;\,\,$} {\tt The $\e\hh_c\e$-module 
structure on quasi-invariants}\newline
\hp7.{ $\;\,\,$} {\tt Differential operators on quasi-invariants}\newline
\hp8.{ $\;\,\,$} {\tt Translation functors and Morita equivalence}\newline
\hp9.{ $\;\,\,$} {\tt Applications of the shift operator}\newline
10.{ $\;\,\,$} {\tt Appendix: A filtration on differential
operators}\newline
}}}

\section{Introduction}
Let $ W $ be a finite Coxeter group in a
complex vector space $\h$,
and $R\subset \h^*$ the corresponding set of roots.
To each 
 $W$-invariant function $c: R \to \C\,,\, c\mapsto c_\alpha,$
 one can  attach an associative
algebra
$\hh_c$, called  the {\it rational Cherednik algebra}.
 This is a very interesting algebra tied to
exciting works in combinatorics, completely integrable systems, and
generalized McKay correspondence.
Historically, the
rational Cherednik  algebra appeared as a
`rational' degeneration of the {\it double-affine 
Hecke algebra} introduced by Cherednik [Ch2].
Thus, the latter may (and should) be  thought of as a deformation
of the former. From this point of view,
representation theory of the
rational Cherednik  algebra  is perhaps
 `more basic' than (or at least should be studied before)
 that  of the double-affine 
Hecke algebra in the same sense as 
the representation theory of semisimple Lie algebras
is  `more basic'  than that  of the corresponding
quantum groups. Thus, one of our goals is to begin
a systematic study of $\hh_c$ representation theory.

Given $\alpha\in\h^*,$ write $\alpha^{\vee} \in \h $
 for the  coroot,
and $ s_\alpha \in \mbox{\rm GL}(\h) 
$ for the  reflection corresponding to $ \alpha$.
Recall  from [EG], that the rational Cherednik algebra $\hh_c$
(which was denoted $\hh_{1,c}$ in  [EG]) 
is generated
by the vector spaces $\h$, $\h^*,$ and the set
$W,$ with defining relations  
(cf.  formula (1.15) of [EG] for $t=1$) given
by
\begin{equation}
\label{defrel}
\begin{array}{lll}\displaystyle
&{}_{_{\vphantom{x}}}w\cd x\cd w^{-1}= w(x)\;\;,\;\;
w\cd y\cd w^{-1}= w(y)\,,&
\forall y\in \h\,,\,x\in \h^*\,,\,w\in W\break\medskip\\
&{}^{^{\vphantom{x}}}{}_{_{\vphantom{x}}}x_1\cd x_2 = 
 x_2\cd x_1\enspace,\enspace
y_1\cd y_2=y_2\cd y_1\,, &
\forall y_1,y_2\in \h,\;x_1\,,\,x_2 \in \h^*\,\break\medskip\\
&{}^{^{\vphantom{x}}}y\cd x-x\cd y = \langle y,x\rangle
-\!\!\!\sum\limits_{^{_{\alpha\in R/{\{\pm 1}\}}}}
c_\alpha\cd\langle y,\alpha\rangle
\langle\alpha^\vee,x\rangle \cd s_\alpha\,,& \forall y\in
\h\,,\,x\in \h^*\,.
\end{array}
\end{equation}
 Thus, the elements  $x\in\h^*$ generate a 
subalgebra $\C[\h]\subset \hh_c$ of polynomial functions on $\h$,
 the elements  $y\in\h$ generate a 
subalgebra $\C[\h^*]\subset \hh_c$,
and the elements $w\in W$
span  a copy of the group
algebra $\C W$ sitting naturally inside $\hh_c$. Furthermore,
it has been shown by Cherednik, see also  [EG],
 that multiplication in $\hh_c$
induces a vector space isomorphism:
\begin{equation}
\label{pbw}
\C[\h]\,\otimes_{_\C}\,\C W\,
\otimes_{_\C}\,\C[\h^*]
\iso \hh_c\qquad\text{\sf{(Poincar\'e-Birkhoff-Witt isomorphism for $\hh_c$)}}.
\end{equation}
 The name for the  isomorphism above comes from its 
analogy with the well-known isomorphism:
$U({\mathfrak{n}}_-)\,\otimes_{_\C}\,
U(\h)\,\otimes_{_\C}\,U({\mathfrak{n}}_+)
\iso U({\mathfrak{g}})$
for the enveloping algebra 
of a complex semisimple Lie
algebra  ${\mathfrak{g}}$ with
triangular decomposition:
${\mathfrak{g}}={\mathfrak{n}}_+
+\h+{\mathfrak{n}}_-\,$. 

The Poincar\'e-Birkhoff-Witt isomorphism
for $\hh_c$ allows one to introduce
a category $\oo_{_{\hh_c}}$ of  modules
over the algebra  $\hh_c$ similar to
the category $\og$ of highest weight modules
over
$U({\mathfrak{g}}),$ considered
by Bernstein-Gelfand-Gelfand, see [BGG]. 
The category $\oo_{_{\hh_c}}$ splits up
into a direct sum of its subcategories
$\oo_{_{\hh_c}}(\bar{\lambda}),\,$
one for each
$\bar{\lambda}\in\h^*/W$. 
We will be mainly concerned below with the category  $\oo_{_{\hh_c}}(0)$,
 which is  most interesting among
all $\oo_{_{\hh_c}}(\bar{\lambda})$'s.
 The  category $\oo_{_{\hh_c}}(0)$ is the
Cherednik algebra counterpart
of the subcategory  $\og(\chi)\subset
\og$ corresponding,
in the Bernstein-Gelfand-Gelfand setting\footnote{
We emphasize that it is the  parameter `$\,c\,$',
and not $\bar{\lambda}\in\h^*/W$ that plays the role of central
character
in representation theory of $\hh_c$. The  parameter
$\lambda\in\h^*$ plays,
in our present situation,
 the role of `Whittaker character',
see remark after Corollary \ref{cor_eHe}.}, to a fixed
character $\chi$ of the center of $U({\mathfrak{g}})$.

The isomorphism (\ref{pbw})
 shows that the group algebra $\C W$ plays
a role of the subalgebra $U(\h)\subset
U({\mathfrak{g}})$.
Thus, following the classical construction due to Verma,
to each irreducible representation $\tau\in\irrep(W)$
one can associate
a `standard' module $M(\tau)\in \oo_{_{\hh_c}}(0),$
an analogue of Verma module. It is easy to show, see [DO],
that each standard module $M(\tau)$ has a unique
simple quotient, $L(\tau)$.
Furthermore, any object of
 the category $\oo_{_{\hh_c}}(0)$ has finite length,
and the collection $\{L(\tau)\}_{\tau\in\irrep(W)}\,$
is a complete collection of  isomorphism
classes of simple objects of  $\oo_{_{\hh_c}}(0)$.
 Thus, the simple objects in 
$\oo_{_{\hh_c}}(0)$ are parametrized by the
set $\irrep(W)$, while  the simple objects in 
$\og(\chi)$
are parametrized (for regular $\chi$)
 by  elements of the Weyl group $W$ itself.

The structure of the category $\oo_{_{\hh_c}}(0)$
depends crucially on the value of the  parameter `$\,c\,$'. 
The category is semisimple for almost all  `$\,c\,$',
see [OR], in which case $M(\tau)=L(\tau),$ for any $\tau$.
However, for a certain  set
of "singular"
 values of  `$\,c\,$', the multiplicities
$[M(\tau):L(\sigma)]$ are unknown.
These multiplicities are Cherednik algebra analogues of
Kazhdan-Lusztig type multiplicities for affine Hecke algebras,
see [CG],
 and they are expected to be provided
by some Intersection cohomology.

In this paper we are concerned with the case
of {\it integral} values of  `$\,c\,$', which
is, in a sense,  intermediate between the two extreme
cases above. 
We will see,
 although it is not
a priori obvious\footnote{See remark at the end of \S2,
and also [OR].},
 that, for any such  `$\,c\,$', the category  $\oo_{_{\hh_c}}(0)$
is  semisimple. Further,
we apply  $\hh_c$-representation theory 
to show that, for integral  `$c$', the
algebra $\hh_c$ is simple and Morita equivalent
to $\dd(\h)\#W$, the cross product of $W$ with the algebra $\dd(\h)$
of  polynomial 
differential operators on $\h$.

Our motivation to study the case of integral   `$c$' comes also
from an interesting
connection with the theory of Calogero-Moser integrable systems
and the theory of differential operators on singular algebraic varieties.
In more detail, let $\C[\h]^W\subset \C[\h],$ be the subalgebra
of $W$-invariant polynomials.
For each  non-negative integral  value of `$\,c\,$',
Chalykh, Feigin, and Veselov [CV], [FV],
have introduced an algebra  $Q_c$ of so-called
$W$-{\it quasi-invariant} polynomials,
 such that $\C[\h]^W\subset Q_c\subset \C[\h].$
We use representation theory of Cherednik algebras
to study differential operators on
quasi-invariants, that is,
the algebra $\dd(Q_c)$ 
 of  differential operators on
the singular variety $\Spec^{\,}Q_c$.
The Calogero-Moser differential operator may be
viewed as
an analogue of the second order Laplacian
on the  variety $\Spec^{\,}Q_c$.

It turns out that the algebra $\dd(Q_c)$
is almost as nice as the algebra of  differential operators on
a smooth variety. Specifically, 
we prove that the algebra $\dd(Q_c)$ is
 Morita equivalent to the algebra of  polynomial 
differential operators on the vector space $\h$.
The subalgebra  $\dd(Q_c)^W$
of $W$-invariant differential operators
will be shown to be
isomorphic to the spherical subalgebra $\ehe$.
Further,  $\dd(Q_c)$ is a simple algebra
equipped with a natural involution,
an analogue of the Fourier transform for differential
operators on $\Spec^{\,}Q_c$.
The algebra  $\dd(Q_c)$ contains  $Q_c$ and
 its
`Fourier dual', $Q_c\ff$,
as two maximal commutative subalgebras.
We show that  $\dd(Q_c)$ is 
generated, as an algebra,  by  $Q_c$ and
 $Q_c\ff$,
and that $\dd(Q_c)$ is a rank one projective
$Q_c\otimes Q_c\ff$-module (under multiplication-action on $\dd(Q_c)$ on
opposite sides).

The rings of differential operators on 
general singular algebraic
 varieties typically have rather unpleasant behavior
(e.g. are not Noetherian, cf. [BGG1]).
The question of simplicity and
 Morita equivalence of such rings
 has been studied by several authors, see [BW], [Sm], [Mu], [SS], 
 [HS], [CS], [VdB], (and also [MvdB] for a  result similar
in spirit to ours). 
The
 varieties $ \Spec^{\,}Q_c $ studied in the present paper 
seem to be  natural generalizations of the
 one-dimensional examples constructed in [Sm] and [Mu].

The paper is organized as follows. Sections 2,3,5 are devoted
to `pure' representation theory of Cherednik algebras.
In \S2 we exploit an 
 idea  due to Opdam relating highest weight modules over the
 Cherednik algebra to finite dimensional representations of the Hecke
algebra
$H_W(e^{2\pi i c})$; thus
\S2 has some overlap with [DO] and [OR].
 In section 3 a new and quite useful
notion of a {\it Harish-Chandra $\hh_c$-bimodule} is introduced.
In \S4 we prove a Cherednik algebra counterpart
of an important
theorem (due to Levasseur-Stafford [LS]), saying that the 
spherical subalgebra in $\hh_c$ is generated, for all regular
values of `$c$', by its two invariant commutative
subalgebras,
$\C[\h]^W$ and $\C[\h^*]^W$. This result has numerous applications.
 In \S5 we compute
a trace on  the Cherednik algebra of type $\mathbf{A}$, and use it to derive some
applications to {\it finite-dimensional} representations of $\hh_c$.
Quasi-invariants are introduced in \S6,
and the structure of the algebra $\dd(Q_c)$ is studied in detail
in \S\S7,9. The structure of $\dd(Q_c)$
 is in turn  exploited in section 8
to get further results in  representation theory
 of Cherednik algebras, e.g., to give 
 an explicit construction
of simple Harish-Chandra bimodules.
 These results of \S8 bear some resemblance 
with the technique of {\it translation functors},
a well-known and very powerful tool in representation
theory of semisimple Lie algebras.
 \medskip

\noindent {\bf Acknowledgments.} {\footnotesize We are very
grateful to E. Opdam and R. Rouquier
for generously
sharing their ideas with us, and for making the results of [OR]
available to us before its publication.
These results play a crucial role in our arguments
and, to a great extent, has triggered the present work.
We would like to thank
T. Stafford for many interesting comments and other
useful information that was quite essential for us.
 The  first
author was partially supported by the NSF grant 
 DMS 00-71792 and A. P. Sloan Research Fellowship; the
work of the second author was  partly conducted for 
the Clay Mathematics Institute and partially
supported by the NSF grant DMS-9988796}.

\section{Standard modules over the rational Cherednik
  algebra}
\setcounter{equation}{0}

Most of the results of this section (in particular,
Theorem \ref{opdam1} and Lemma \ref{opdam2})
are due to Opdam-Rouquier [OR], and are reproduced here
for the reader's convenience only.

Fix a finite Coxeter group $W$ in a
complex vector space $\h$. Thus, $\h$ is the complexification
of a  real Euclidean  vector space and $W$ is
generated by reflections with respect to a certain finite set
$\{H_\alpha\}$ of hyperplanes in that  Euclidean  space.
We write  $ (\,\cdot\,,\,\cdot\,)\,$ for the complex bilinear
form on $\h$ extending the  Euclidean inner product.
For each hyperplane, $H_\alpha,$ we choose
 nonzero
linear functions $\pm\alpha\in\h^*$ which vanish on $H_\alpha$.
The set  $R\subset \h^*$ of all such linear functions is called
the set of roots.
Write $ s_\alpha \in \mbox{\rm GL}(\h) 
$ for the   reflection corresponding to $ \alpha\in R,$ and
$ \alpha^{\vee} \in \h $ for the corresponding coroot, 
a vector such that $ s_\alpha(\alpha^\vee) = -\alpha^\vee$.
The lengths of roots and coroots are normalized so that
$ \langle \alpha,\alpha^\vee \rangle = 2\,$.    
We make a choice of the set $R_+$ of positive roots
so that $R=R_+\sqcup (-R_+)$.
The  $ W $-action on $ \h $  induces canonical
actions on the
symmetric algebras $ \mbox{Sym}(\h) =\C[\h^*]$ and $
\mbox{Sym}(\h^{*}) =\C[\h] \,$.

Let  $\hr$ denote the complement to the union
of all the  reflection
hyperplanes, i.e., the complement to
 the zero set of the discriminant polynomial
$\delta=\prod_{\alpha\in R_+}\,\alpha\in\C[\h].$
Given a $\C[\h]^W$-module $M$, we will
write $M\mrw$ for 
$\C[\hreg]^W\otimes_{_{\C[\h]^W}}M,$
the localization to $\hreg/\!W$. Note that if $M$ is 
a $\C[\h]$-module  viewed as a $\C[\h]^W$-module by restriction
of scalars, then  $M\mrw$ is a $\C[\hreg]$-module which
coincides with
$M\mr:=\C[\hreg]\otimes_{_{\C[\h]}}M.$ In particular,
 $\hh_c\mr =\C[\hreg]\otimes_{_{\C[\h]}}\hh_c$, the localization
of the left regular $\hh_c$-module,
acquires a natural algebra structure,
hence an $\hh_c$-bimodule  structure, such that the
imbedding $\hh_c\into \C[\hreg]\otimes_{_{\C[\h]}}\hh_c$
becomes an algebra map. Alternatively, the algebra
$\hh_c\mr$ is obtained by Ore localization of $\hh_c$ with respect to
the multiplicative set~$\{\delta^k\}_{k=1,2,\ldots}$.

The group $W$ acts freely on
 $\hr$. 
Let $B_W$ be the braid group of $W$, that is
 the fundamental group of the 
variety $\hr/W$. Fix a point $*\in\hr$ inside a Weyl chamber in $\h$,
and for
 each simple  reflection $s_\alpha\in W$, let $T_\alpha$ be the
class in $B_W=\pi_1(\hr,\,*)$ corresponding to a straight  path from
the point $*$ to the point $s_\alpha(*)$ with an inserted
 little semi-circle
(oriented counter-clockwise)
around the hyperplane $\alpha=0$.
Given a $W$-invariant function $q: R \to \C^{\times}\,,\,
\alpha\mapsto q_\alpha\,$,
let $H_W(q)$ be the Hecke algebra. This algebra is  obtained by taking 
the quotient of $\C[B_W]$ by the relations
$(T_\alpha-1)(T_\alpha+q_\alpha)=0$, one
for each simple reflection $s_\alpha$. It is known that
$\dim H_W(q)=|W|,$ for any function $q: R \to \C^{\times}$.

\begin{definition}\label{Reg}
A $W$-invariant function $q: R \to \C^{\times}$ is said to be
{\sl regular} if the Hecke algebra  $H_W(q)$ is semisimple.
Write $\Reg$ for the set of regular
$W$-invariant functions $q$.
\end{definition}

Since $H_W(1)=\C W$, we see that $q=1$ is a regular function.
Moreover, the set $\Reg$ is a dense Zariski open subset in
the set of all $W$-invariant functions $q: R \to \C^{\times}\,$
(the latter set is naturally identified with $(\C^*)^l,$
where $l$ is the number of $W$-orbits in $R$).
The complement of $\Reg$ has real codimension $\geq 2$,
therefore  $\Reg$ is a connected set. Hence, using
rigidity of semisimple algebras one proves
that $H_W(q)\simeq\C{W},$ for any $q\in\Reg.$

We define  standard modules over the
Cherednik
algebra $\hh_c$, see (\ref{defrel}), as follows. 
Fix $ \lambda \in \h^* \, $, and let $W_\lambda$ be 
the stabilizer of $\lambda$ in $W\,$. 
Let $\C[\h^*]\#W_\lambda$ be the cross-product of $W_\lambda$ 
 with the polynomial
algebra.
Sending $P\in \C[\h^*]$ to $P(\lambda)\cdot 1$ yields an algebra homomorphism:
$\C[\h^*]\#W_\lambda\onto\C{W_\lambda}$. Given 
$ \tau \in \mbox{\rm Irrep}(W_\lambda)\,$,
an irreducible representation of $W_\lambda\,$,
we write  $\lambda\#\tau$ 
for the  representation of $\C[\h^*]\#W_\lambda$ obtained via
the pull-back by the
homomorphism above.
Then we set $M(\lambda,\tau)
:=\text{Ind}_{\C[\h^*]\#W_\lambda}^{\hh_c} (\lambda\#\tau)$.
The module $M(\lambda,\tau)$ is called a {\it standard module}. 
In particular, 
if $\lambda=0$ we have $W_\lambda=W$, in which case
we write $ M(0, \tau) = M(\tau) \,$ for the corresponding standard module.

Let $\CC$ denote the vector space of
all $W$-invariant functions $c: R\to\C$,
and $\ccr\subset\CC$  the set of  functions
$c\in\C[R]^W$ such that $\exp(2\pi i c)\in \Reg$,
i.e., such that the Hecke algebra
$H_W(e^{2\pi i c})$ is semisimple.
The goal of this section is to prove the following result
 due to Opdam-Rouquier in the key special case $\lambda=0$.

\begin{theorem}[\cite{OR},\cite{GGOR}]\label{opdam1} If $c \in
\ccr$ then, for any $\lambda\in\h^*$ and $\tau\in \mbox{\rm
Irrep}(W_\lambda)\,$,
the standard module $M(\lambda,\tau)$ is a simple $\hh_c$-module.
\end{theorem}

 Since $\exp(2\pi i c)=1\in \Reg$ for any integral valued
function $c$, Theorem \ref{opdam1} yields

\begin{corollary} All standard modules, $M(\lambda,\tau),$
are simple, for any 
$c\in\Z[R]^W$.\qed
\end{corollary}

By analogy with representation
 theory of semisimple Lie algebras we introduce
the following

\begin{definition}\label{cat_O} Let $\oo_{_{\hh_c}}$
be the  category of finitely-generated
$\hh_c$-modules $M$, such that the action
on $M$ of the subalgebra $\C[\h^*]\subset\hh_c$ is
locally finite, i.e., $\dim_{_\C}\C[\h^*]\cd m<\infty,$
for any $m\in M$.
\end{definition}

We say that an object $M\in \oo_{_{\hh_c}}$ has {\it type}
${\bar{\lambda}}\in \h^*/W=\Spec(\C[\h^*]^W)$ if , for any $P\in \C[\h^*]^W$,
the action on $M$ of the element $P-P(\bar{\lambda})\in
\hh_c$ is locally nilpotent. Let $\oo_{_{\hh_c}}(\bar{\lambda})$
be the full subcategory of modules having type
${\bar{\lambda}}$. Then, by a routine argument, cf. e.g. [Di],
one obtains a direct sum decomposition:
$\oo_{_{\hh_c}}=\bigoplus_{{\bar{\lambda}}\in \h^*/W}\;
\oo_{_{\hh_c}}(\bar{\lambda})$.

\begin{lemma}\label{O_prop} \vi Any object $M\in \oo_{_{\hh_c}}$
is finitely generated over the subalgebra $\C[\h]\subset\hh_c$,
in particular, $\oo_{_{\hh_c}}$ is an abelian category.

\vii For any $\lambda\in\h^*$ and
$ \tau \in \mbox{\rm Irrep}(W_\lambda)\,$,
we have $M(\lambda,\tau)\in \oo_{_{\hh_c}}$.

\viii For any $M\in \oo_{_{\hh_c}}$, there exists
a nonzero homomorphism $M(\lambda,\tau)\to M$,
for certain $\lambda\in\h^*$ and
$ \tau \in \mbox{\rm Irrep}(W_\lambda).$

${\sf{(iv)}}\;$ Every object of the category
 $\oo_{_{\hh_c}}(0)$ has finite length.
\end{lemma}

\begin{remark} Using Lemma \ref{O_prop}(iv), it is proven in [Gi] that, more generally,
every
 object of the category
 $\oo_{_{\hh_c}}$ also  has finite length.
\end{remark}
\smallskip

To prove Lemma \ref{O_prop}, we need some notation.
Let $\{x_i\}$ and $\{y_i\}$ be a pair of dual bases of $\h^*$
and $\h$, respectively. We view
$\h^*$ and $\h$ as subspaces in $\hh_c$,
and let $\sh=\frac{1}{2}\sum_i\,(x_iy_i+y_ix_i)\in\hh_c$ denote the
canonical element, which is independent of the choice
of the bases. The element $\sh$ satisfies the following
commutation relations: 
\begin{equation}\label{sh}
\sh\cdot x = x\cdot (\sh+1)\enspace,\enspace
\forall x\in \h^*,\quad\text{and}\quad\sh\cdot y=y\cdot (\sh-1)
\enspace,\enspace\forall y\in\h\,.
\end{equation}
To prove the first of these formulas, 
use the notation $[a,b]: = a\cdot b - b\cdot a$.
Given $x\in \h^*,$ we calculate
\begin{align}\label{sh1}
&[\sh, x]=\frac{1}{2}\big[\sum_i\,(x_iy_i+y_ix_i)\,,\,x\big]=
\frac{1}{2}\sum_i\,\bigl(x_i\cd [y_i,x] + [y_i,x]\cd x_i\bigr)
\nonumber\\
&=\frac{1}{2}\sum_i\,\Bigl(
\bigl(x_i\cd \langle y_i,x\rangle+
\langle y_i,x\rangle\cd x_i\bigr)\,-\,\sum_{\alpha\in R_+}\;
\frac{c_\alpha}{2}
\bigl(x_i\cd \langle y_i,\alpha\rangle\langle\alpha^\vee, x \rangle 
s_\alpha
+\langle y_i,\alpha\rangle\langle\alpha^\vee, x \rangle
s_\alpha\cd x_i\bigr)\Bigr)\nonumber\\
&=\sum_i\,\bigl(x_i\cd \langle y_i,x\rangle\bigr)-
\sum_{\alpha\in R_+}\;\frac{c_\alpha}{2}\cdot 
\langle\alpha^\vee, x \rangle\cdot  \Bigl(\sum_i\;
\bigl(\langle y_i,\alpha\rangle x_i\cd s_\alpha
+\langle y_i,\alpha\rangle
s_\alpha\cd x_i\bigr)\Bigr)\nonumber\\
&=x-
\sum\nolimits_{\alpha\in R_+}\;\frac{c_\alpha}{2}\cdot
\langle\alpha^\vee, x \rangle\cdot \bigl(\alpha\cd s_\alpha+
s_\alpha\cd\alpha\bigr)\,,
\end{align}
where in the last equality we have used the identities
$\sum_i\,\langle y_i,x\rangle\cdot x_i=x,$ and
$\sum_i\langle y_i,\alpha\rangle\cdot   x_i$
$=\alpha$.
Since $(s_\alpha)^{-1}\cdot\alpha\cdot s_\alpha
=s_\alpha(\alpha)=-\alpha,$ we find
$$\alpha\cdot s_\alpha+s_\alpha\cdot\alpha=
s_\alpha\cdot(s_\alpha)^{-1}\cdot\alpha\cdot s_\alpha+s_\alpha\cdot\alpha=
s_\alpha\cdot(-\alpha)+s_\alpha\cdot\alpha=0\,.
$$
Thus, each term in the sum on the last line of (\ref{sh1}) vanishes,
and
we deduce: $[\sh, x]=x$. This proves the first identity in (\ref{sh});
the second one is proved similarly.\qed
\medskip

Let $V$ be a bimodule over an associative
 algebra $A$. For any $a\in A$,
we  have an adjoint $a$-action on
$V$ given by $\ad a: v\mapsto av-va$.
We say that the adjoint action of $A$ on $V$ is
{\it locally nilpotent} if, for any $v\in V,$ there
exists an integer $n\ge 0$ such that
$\ad a_0\ccirc\ad a_1\ccirc\ad a_2\ccirc\ldots\ccirc\ad a_n(v)=0,$
for any $a_0,\ldots,a_n\in A$.
Let $n(v)$ be the smallest among such integers $n\ge 0$,
to be referred to as the order of nilpotency of $v$.

\begin{lemma}\label{tens2} Let $A$ be a finitely generated
commutative algebra and $V$ an $A$-bimodule, such that
the adjoint action of $A$ on $V$ is locally nilpotent. Then for any
$v\in V$, the space $AvA$ is finitely generated both as a
left and as a right $A$-module.
\end{lemma}

\begin{proof} We proceed  by induction in the order of nilpotency of $v$.
If ${n(v)}=0$, the result is clear, since $AvA=Av.$ So,
we have to prove the statement
for ${n(v)}=m,$ knowing it for ${n(v)}=m-1.$

Fix $v\in V$ of nilpotency  order $m$.
Let $a_1,...,a_d$ be generators of $A$,
and put: $u_i=[a_i,v]\,,\, i=1,...,d$.
We have: $A\cd v\cd A=\sum A\cd v\cd a_{i_1}\cd\ldots\cd a_{i_N}.$
By commuting $v$ with $a_{i_1},$ we get
$A\cd v\cd a_{i_1}\cd a_{i_2}\cd\ldots$
$\cdot
a_{i_N}
\subset A\cd v\cd a_{i_2}\cd\ldots\cd
a_{i_d}+A\cd u_{i_1}\cd A.$ Now, continuing like this
(i.e. interchanging $v$ with $a_{i_2}$ etc.), we get
$A\cd v\cd A\subset A\cd v+\sum_{i=1}^d A\cd u_i\cd 
A.$ By the induction assumption,
the module on the right is finitely generated.
Thus, by the Hilbert-Noether lemma, so is the module $AvA.$
\end{proof}

\noindent
{\sl Proof of Lemma \ref{O_prop}.\,\,}
(i) Let $ M_0\subset M$  be the finite dimensional vector space
spanned by a finite set of generators. Then $ \hh_c\cdot  M_0=M.$ 
On the other hand, the space, $ M_0':=\C[\h^*]\cdot 
M_0$  is finite dimensional, since the
action of $ \C[\h^*]$  is locally finite.
So $ M=\hh_c \cdot M_0=\C[\h]\cdot \C W\cdot \C[\h^*]\cdot M_0=\C[\h]\cdot
\C W\cdot 
M_0'.$  But $ \C W\cdot M_0'$  is finite dimensional, so
we are done.

(ii) 
The $\hh_c$-module $M(\lambda,\tau)$ is generated by
the $W$-stable finite dimensional subspace $E:=\C W\cdot\tau\subset
M(\lambda,\tau).$ Hence, given  $ v \in  M(\lambda,\tau),$
there
exists  a finite dimensional 
subspace $B\subset \hh_c$ such that
$ \C[\h^*]\cdot v\subset\C[\h^*]\cdot B\cdot E$.
Moreover, since $\C[\h^*]$ is finite over $\C[\h^*]^W$,
we may find $B$ large enough, so that
$\C[\h^*]\cdot E\subseteq\C[\h^*]^W\cdot B\cdot E$.
But $\C[\h^*]^W\cdot B\subset\C[\h^*]^W\cdot B\cdot \C[\h^*]^W.$
By Lemma \ref{tens2} applied to $V=\hh_c$,
there exists a  finite dimensional subspace $B'\subset \hh_c$,
 such that $\C[\h^*]^W\cdot B\cdot \C[\h^*]^W\subset
B'\cdot 
\C[\h^*]^W.$ Hence,
$ \C[\h^*]^W\cdot B\cdot E\subset\C[\h^*]^W\cdot B\cdot \C[\h^*]^W\cdot E
\subset B'\cdot 
\C[\h^*]^W\cdot E$. The latter space is finite dimensional since
$ \C[\h^*]^W \cdot E$  is finite dimensional, and (ii) follows.

(iii) The action of the subalgebra $\C[\h^*]\subset \hh_c$
on $M\in \oo_{_{\hh_c}}$
being locally finite, it follows that  
$M$ 
contains a vector annihilated
by the  maximal ideal $J$  of the algebra $\C[\h^*]$
corresponding to a point $\lambda\in\h^*$.
The elements of $M$ which are annihilated by $J$
clearly form an $W_\lambda$-stable vector space. 
Hence, this  vector space contains a simple $W_\lambda$-module
$\sigma$. Any vector in this $\sigma$
gives rise to a nonzero element
of $\Hom_{\C[\h^*]\#W_\lambda}(\lambda\#\sigma, M)=
\Hom_{\hh_c}\bigl(M(\lambda,\sigma),M\bigr)$.

To prove  (iv), we observe that
the $\sh$-action on any standard module $M(\tau)$ is diagonal,
with finite dimensional  eigenspaces.
It follows that, for any $\tau\in\irrep(W)$, the multiplicity
of each simple object of $\oo_{_{\hh_c}}(0)$
in $M(\tau)$ is bounded from above by the dimension
of the corresponding  $\sh$-eigenspace.
Further, the category $\oo_{_{\hh_c}}(0)$ has only finitely
many simple objects,
and therefore, each $M(\tau)$
has a finite Jordan-H\"older series.

We say that
an object $N \in \oo_{_{\hh_c}}(0)$ is a highest weight object if
it is a quotient of a standard module. It follows from the
paragraph above that  any highest weight object
has finite length.

To prove that an arbirary object $M\in \oo_{_{\hh_c}}(0)$
has finite length it suffices to show  that $M$ has a finite filtration
by subobjects
$0=F_0\subset F_1\subset\ldots\subset F_n=M,$
such that successive quotients $F_i/F_{i-1}$ are highest weight objects.
To construct such a filtration,
we use the fact that $M$ contains
a nonzero highest weight submodule $N_1$, by part (iii) of the Lemma. 
Put
$M_1:=M/N_1$, which  is also an object of $\oo_{_{\hh_c}}(0)$.
If $M=N_1$ we are done; if not, then $M_1$ also
contains a  nonzero highest weight object $N_2.$
We set $M_2=M_1/N_2,$ and iterate the procedure.
Thus, for each $i\geq 1$, we get an object $M_i$ which
is a quotient of $M$.
Let $F_i:=\ker(M\onto M_i)$.
Clearly, $F_i\in \oo_{_{\hh_c}}(0)$, and we have an increasing
chain of subobjects in $M$:
$0=F_0\subset F_1\subset\ldots\,$.
  Since $M$ is a Noetherian $\C[\h]$-module, this sequence
must stabilize. This means that at some step
$l$, the object $N_{l+1}$ cannot be chosen, which means that
$M_l=0$. Thus, the
$\{F_i\}$ provide a finite filtration with quotients being
highest weight objects, and the result follows.
\qed

\begin{remark} The argument above proves also that
if for some $c\in\CC$ and $\bar\lambda\in \h^*/W,$
all the standard modules $M(\lambda,\tau)$ are simple,
then every object of $\oo_{_{\hh_c}}(\bar\lambda)$
has finite length. As we will see below, these conditions
always hold, for instance, provided $c\in\ccr$.
\end{remark}
\medskip

\noindent
{\bf Corollary.\;\;}
{\it The action of $\sh$ on
any object $M\in\oo_{_{\hh_c}}(0)$ is locally finite,
that is 
$\dim(\C[\sh]\cdot m) < \infty,$ for any $m\in M$.}
\qed
\medskip

We need to recall a few
important results about Dunkl operators.
According to Cherednik, see [EG] Proposition 4.5,
the algebra $\hh_c$ has a faithful "Dunkl representation",
an injective algebra homomorphism
$\Theta: \hh_c \into \dd(\hreg)\#W$.
This morphism extends by $\C[\hreg]$-linearity
to a map:
$ \hh_c\mr=\C[\hreg]\otimes_{_{\C[\h]}} \hh_c\too \dd(\hreg)\#W$,
which gives an  algebra isomorphism
$\Theta:\hh_c\mr\iso \dd(\hreg)\#W$
(surjectivity is clear since the set of Dunkl operators 
generates $\dd(\hreg)\#W$ over the subalgebra $\C[\hreg]\#W$).
Now, for any $\hh_c$-module $M$, 
the
localization, $M\mr,$ has a natural
$\hh_c\mr$-module
structure, therefore, a $\dd(\hreg)\#W$-module
structure, via the Dunkl isomorphism above. Thus,
since $\hreg$ is affine, one can view $M\mr$
as a $W$-equivariant $\dd$-module on $\hreg$.

Observe further that the $\dd$-module on $\hreg$
arising from  any object
$M\in \oo_{_{\hh_c}}$ is
finitely generated over the subalgebra
$\C[\hr]$,
by Lemma \ref{O_prop}(i). Hence, $M\mr$, viewed
as a $\dd$-module on $\hr$, must be a vector
bundle with flat connection. In particular,
 the standard module $M(\tau)$  is a free $\C[\h]$-module 
of rank $\text{dim}(\tau),$ as follows
 from Poincar\'e-Birkhoff-Witt isomorphism (\ref{pbw}).
Thus $M(\tau)\mr$, viewed as a  $\dd$-module  on $\hreg,$
is
the trivial vector bundle $\C[\hreg]\otimes\tau$ equipped with a flat 
connection.

The following well-known result plays a crucial role in this
paper. 

\begin{proposition}\label{dr1}
\vi {\sf (Dunkl \cite{D}, Cherednik \cite{Ch})\;}
 The connection on  $\C[\hreg]\otimes\tau$ arising
from the standard module $M(\tau)$ is the Knizhnik-Zamolodchikov
connection with values in~$\tau$.

\vii {\sf (Opdam [O])\;} The monodromy 
representation of the fundamental group
$B_W=\pi_1(\hreg\!/W)$ corresponding to this
connection factors through the Hecke algebra $H_W(e^{2\pi ic})$.
\end{proposition}

\begin{proof} Although this is a known result, we will give a
short  proof for reader's convenience. 
Let us write down a system of linear differential equations 
which defines the horisontal sections of the corresponding
 connection. For this purpose, consider the generating subspace
$\tau\subset M(\tau)\mr$. 
We have $yv=0$ for $y\in \h$ and
$v\in \tau$. Therefore (since the isomorphism $\Theta: \hh_c\mr\to
\dd(\hr)\#W$ is defined via Dunkl operators), we have 
$$
\partial_yv=\biggl[y-\sum\nolimits_{\alpha\in
  R_+}\;c_\alpha\frac{\langle \alpha,y\rangle}
{\alpha}(s_\alpha-1)\biggr]v=
\sum\nolimits_{\alpha\in
  R_+}\;c_\alpha\frac{\langle \alpha,y\rangle}
{\alpha}(1-s_\alpha)v
$$
In the trivialization $M(\tau)\mr\simeq
\C[\hreg]\otimes\tau$, the above formulas  for the Dunkl
operators equip  the trivial bundle $\C[\hreg]\otimes\tau$
with  the flat connection
$$\nabla= d +\sum\nolimits_{\alpha\in
  R_+}\;c_\alpha\frac{d\alpha}
{\alpha}\otimes(\id_\tau-s_\alpha)\;:\quad \C[\hreg]\otimes\tau
\too \Omega^1(\hreg)\otimes\tau\,.
$$
This vector bundle with connection $\nabla$ is $W$-equivariant,
hence, descends to $\hreg/W$.
The corresponding flat sections are (multivalued) functions $
f: \hr/W\to \tau
$ satisfying the holonomic system  of differential equations:
$
\partial_yf=\sum\nolimits_{\alpha\in
  R_+}\;c_\alpha\frac{\langle \alpha,y\rangle}
{\alpha}(1-s_\alpha)f$.
An elementary theory of ordinary
differential equations says that, for generic $c$, the image 
of  any element $T_\alpha$ in the monodromy representation of
this system
 satisfies 
the equation $(T_\alpha-1)(T_\alpha+e^{2\pi ic_\alpha})=0$. By continuity, the same
holds for all $c$,  see [O] for more details.
\end{proof}

In general, given a $W$-equivariant
vector bundle $\mm$ on $\hreg$ with flat connection,
the germs of horizontal holomorphic sections of $\mm$
form a locally constant sheaf on $\hreg\!/W$.
Let $\Mon(\mm)$ be the corresponding
monodromy representation of the fundamental group
$\pi_1(\hreg\!/W, *)$
in the fiber over $*$, where $*$ is some fixed point in $\hreg\!/W$.
The assignment: $\mm\mapsto \Mon(\mm)$ gives a functor
from the category of $W$-equivariant vector bundles on
$\hreg$  with flat connection
to the category of finite dimensional
representations of the group $\pi_1(\hreg\!/W, *)=B_W$.

The following result is due to Opdam-Rouquier.

\begin{lemma}[\cite{OR},\cite{GGOR}]\label{opdam2}
Let  $N$ be an object of $\oo_{_{\hh_c}}$ 
which is torsion-free over the subalgebra $\C[\h]\subset \hh_c$. Then
the canonical
 map: $\Hom_{\hh_c}(M,N)\to
\Hom_{_{B_W}}\bigl(\Mon(M\mr), \Mon(N\mr)\bigr)$
is {\sl injective}, for any
$M\in\oo_{_{\hh_c}}$.
\end{lemma}

\begin{proof}\!\!{\sf {(borrowed from \cite{OR},\cite{GGOR}).}}\;
By assumption, the module  $N$
contains no $\C[\h]$-submodule supported on a proper
subset of $\h$. It follows that if
$f: M\to N$ is an $\hh_c$-module morphism such
that the induced morphism: $ M\mr \to N\mr$ 
has zero image, then $f$ itself has the zero image, that is $f=0$.
Thus, for any
$M\in\oo_{_{\hh_c}}$, the canonical
map $i_1: \text{Hom}_{\,\hh_c}(M,N)\to
\text{Hom}_{_{\hh_c\mr}}\bigl(M\mr,N\mr\bigr)$  is injective.

Now, as has been explained earlier, we may regard
$\mm:=M\mr$ and ${\mathcal{N}}:=N\mr$ as $W$-equivariant
vector bundles on $\hr$ with flat connections.
Hence, the map of the Lemma can be factored as 
a composition
\begin{align*}
&\text{Hom}_{\,\hh_c}(M,N)\stackrel{i_1}{\into} 
\text{Hom}_{_{\hh_c\mr}}\bigl(M\mr,N\mr\bigr)\\
&=\text{Hom}_{\dd(\hreg)\#W}(\mm,{\mathcal{N}})
\iso\enspace\left\lbrace\begin{array}{cc}\text{global $W$-invariant
horizontal}\\
\text{ sections of}\enspace (\mm^*\otimes{\mathcal{N}})
\end{array}\right\rbrace\\
&\stackrel{i_2}{\into}
\bigl(\Mon(\mm^*\otimes{\mathcal{N}})\bigr)^{B_W}=
\text{Hom}_{B_W}\bigl(\Mon(M\mr), \Mon(N\mr)\bigr).
\end{align*}

We already know that $i_1$ is injective.
Further, any $\dd(\hreg)\#W$-module morphism between vector
bundles with flat connections is given by
a horizontal section of the Hom-bundle. The latter is nothing
but a global horizontal section of
$\mm^*\otimes{\mathcal{N}},$
hence $i_2$ is injective. Finally, assigning to a horizontal
section its value in the fiber 
at a given point $*\in \hreg\!/W$ gives an injection
(which is a bijection, if the connection 
has regular singularities), and the Lemma follows.~\end{proof}

The following corollary of the above results
and its proof
 were communicated to us by E. Opdam.

\begin{corollary}[\cite{OR},\cite{GGOR}] 
\label{opdam} 
If $c\in \ccr$, then
all standard modules $M(\tau)$ are irreducible. 
\end{corollary}

\begin{proof} 
Fix $\tau\in \irrep(W)$, and given $c\in\CC$ write
$\Mon_\tau(c):=\Mon\bigl(M(\tau)\mr\bigr)$
 for the corresponding monodromy representation
of the Hecke algebra $H_W(e^{2\pi ic})$.
The assignment $c\mapsto \Mon_\tau(c)$ gives a continuous function
on  the set $\ccr.$
Furthermore,  this set is connected and
contains the point $c=0$, for which we have:
$\Mon_\tau(0)=\tau$. Since, $H_W(q)\simeq \C{W}$ for all $q\in\Reg$, a simple deformation argument shows that
$\Mon_\tau(c)$ is a simple $H_W(e^{2\pi ic})$-module, for any 
$c\in \ccr$. Since $M(\tau)$ is free over
$\C[\h]$, hence torsion-free,
Lemma \ref{opdam2} yields:
$\dim\Hom_{\hh_c}\bigl(M(\tau),M(\tau)\bigr)=1.$

Next, fix two non-isomorphic $W$-modules 
 $\tau, \sigma\in \irrep(W).$ A similar argument
shows that, for each $c\in \ccr,$ the
$H_W(e^{2\pi ic})$-modules $\Mon_\tau(c)$ and 
$\Mon_\sigma(c)$ are simple and non-isomorphic to each other.
Thus, $\text{Hom}_{B_W}(\Mon_\sigma(c), \Mon_\tau(c))=0$.
Hence, Lemma \ref{opdam2} yields:
$\Hom_{\hh_c}\bigl(M(\sigma),M(\tau)\bigr)$ $=0$.
Thus, we conclude: $\dim\Hom_{\hh_c}\bigl(M(\sigma),M(\tau)\bigr)=
\dim\Hom_W(\tau, \sigma)\,,\,\forall \sigma,\tau\in\irrep(W).$

Now, assume $M(\tau)$ is not a simple $\hh_c$-module,
and let $M\subset M(\tau)$ be a proper nonzero
submodule.
Then, Lemma \ref{O_prop}(iii) says that
there exists   $\lambda\in\h^*$ and
$ \sigma \in \mbox{\rm Irrep}(W_\lambda)$
such that
$\Hom_{\hh_c}\bigl(M(\lambda,\sigma), M\bigr)
\neq 0$. Furthermore, the direct sum decomposition
$\oo_{_{\hh_c}}=\bigoplus_{{\bar{\lambda}}\in \h^*/W}\;
\oo_{_{\hh_c}}(\bar{\lambda})$ implies that  $\lambda$
must be zero,
so that $M(\lambda,\sigma)=M(\sigma)$ and
$ \sigma \in \mbox{\rm Irrep}(W).$
Hence, the dimension of
$\Hom_{\hh_c}\bigl(M(\sigma),M(\tau)\bigr)$
is at least $1$ if $\sigma\neq\tau$, and is at least
$2$ if $\sigma=\tau$.
But this contradicts the dimension equality:
$\dim
\Hom_{\hh_c}\bigl(M(\sigma),M(\tau)\bigr)=\dim\Hom_W(\tau, \sigma)$
proved earlier.
The contradiction shows that  $M(\tau)$ is simple.
\end{proof}

The module $M(\lambda,\tau)$ can be
identified, as a vector space,
 with $\C[\h]\otimes \text{Ind}_{W_\lambda}^W\tau,$
via the PBW-theorem.
We fix such an identification.
Thus, the standard  increasing filtration
on  $\C[\h]$ by degree of polynomials
gives rise to an increasing
filtration  $F_\bullet M(\lambda,\tau)$ on $M(\lambda,\tau),$ 
such that
 $F_0M(\lambda,\tau)=\text{Ind}_{W_\lambda}^W\tau.$ It
is clear that operators corresponding to elements $y\in \h$
preserve this filtration. 

\begin{lemma}\label{commeig} If $c\in \ccr$, then
any $\h$-weight vector  in
  $M(\lambda,\tau)$ belongs to $F_0M(\lambda,\tau)$. 
\end{lemma}

\begin{proof}  The proof is based on a continuity
argument and on the fact that the
result clearly holds true, provided $M(\lambda,\tau)$
is a simple $\hh_c$-module, e.g.,
for $\lambda=0$. 

First, we form  the induced module
 $X(\lambda)=\text{Ind}_{\C[\h^*]}^{\hh_c} \lambda$.
We can (and will) canonically identify,
via the PBW-theorem, the underlying  vector space
of the module $X(\lambda)$
 with $\C[\h]\otimes \C W$. As above,
 the standard  increasing filtration
on  $\C[\h]$ by degree of polynomials
gives rise to an increasing
filtration  $F_\bullet X(\lambda)$,
such that $F_0 X(\lambda)=\C W$.

Let ${\sf Weight}_{\,}X(\lambda)$ be the sum of all the
$\h$-weight subspaces of $X(\lambda).$
It is clear
that, for any $\lambda\in\h^*$,
the space ${\sf Weight}_{\,}X(\lambda)$ contains $F_0 X(\lambda).$
 Further, one has an
$\hh_c$-module direct sum decomposition: $X(\lambda)=\bigoplus_{\tau\in
\irrep(W_\lambda)}\,\tau^*\otimes
M(\lambda,\tau),\,$ where $\tau^*$ stands for the
representation dual to $\tau$.
Therefore, to prove the Lemma
it  suffices to show that every 
$\h$-weight vector in $X(\lambda)$ belongs to 
$F_0 X(\lambda)=\C W$.

Fix $\lambda\in\h^*$, and for any $t\in\C,$
view ${\sf Weight}_{\,}X(t\cdot\lambda)$ 
as a subspace in $\C[\h]\otimes \C W$
(note that 
 the latter is {\it independent} of $t\in\C$).
For any integer $d\geq 0$, 
consider the set
$$U_d=\big\{t\in\C\enspace\big|\quad
\bigl({\sf Weight}_{\,}X(t\cdot\lambda)\bigr)\,\cap\,
F_dX(t\cdot\lambda)
\;=\;F_0 X(t\cdot\lambda)\big\}\,.
$$
It is clear that $\C=U_0\supset U_1\supset\ldots$.
Observe that by
Corollary \ref{opdam},
$M(\tau)$ is a simple $\hh_c$-module,
for any $\tau\in\irrep(W)$. It follows that for $t=0$
one has ${\sf Weight}_{\,}X(t\cdot\lambda)=F_0 X(t\cdot\lambda).$
Hence, each set $U_d$ contains the point $t=0$ and
is Zariski open, by an elementary continuity argument.
 Hence, for each $d\geq 0$, the set
$\C\ssminus U_d$ is a proper Zariski closed,
hence finite, subset of $\C$.
 It follows that the set $\cup_{d\geq 0}\,(\C\ssminus U_d)$
is a union of a countable family of finite sets, hence a countable
subset in $\C$. 

Let $T$ be the set of all $t\in \C$ such that 
${\sf Weight}_{\,}X(t\cdot\lambda)$ $=F_0 X({t\cdot{\lambda}})$. 
By definition we have $T=\cap_{d\geq 0}\,U_d.$ Therefore,
$\C\ssminus T=
\cup_{d\geq 0}\,(\C\ssminus U_d)$
is a  countable set.
Further,  this  countable set 
 must be invariant under dilations,
since the assignment $x\mapsto t\cdot x\,,\,y\mapsto t^{-1}\cdot 
y$, $\,t\in \C^*$, $x\in \h^*$, $y\in \h,$ gives a
   $\C^*$-action on $\hh_c$ by algebra automorphisms.
Finally, we know that $0\in T$, i.e.,
  $0\not\in\C\ssminus  T$.
These properties force the set  $\C\ssminus T$
to be empty, hence $T=\C$, and we are done.  
\end{proof} 

\noindent
{\sl Proof of Theorem \ref{opdam1}.\,}
Let $N$ be a nonzero 
submodule of $M(\lambda,\tau)$. Thus, $N\in \oo_{_{\hh_c}}$.
Hence, $N$ contains a nonzero $\h$-weight vector $v$.
 Lemma \ref{commeig} yields
$v\in F_0M(\lambda,\tau)$. But $F_0M(\lambda,\tau)
=\text{Ind}_{W_\lambda}^W\tau$ is an
irreducible $\C[\h^*]\#W$-module, so $v\in
F_0M(\lambda,\tau)\subset N$ implies
$N=M(\lambda,\tau)$.\qed
\smallskip

\begin{remark} It has been shown in [DO] that 
each standard module $M(\tau)\,,\,\tau\in\irrep(W),$
has a unique simple quotient $L(\tau)$. Moreover,
any object of the category
$\oo_{_{\hh_c}}(0)$ has finite length, and 
the finite set $\{L(\tau)\}_{\tau\in\irrep(W)}$ is a
complete collection of the
isomorphism classes of simple objects of $\oo_{_{\hh_c}}(0)$.
Further, it is proved in [Gu] and [OR] that
the category $\oo_{_{\hh_c}}(0)$ has enough projectives,
i.e., each module $L(\tau)$ has an indecomposable 
projective cover $P(\tau)\in \oo_{_{\hh_c}}(0)$.
Moreover, any projective  $P\in\oo_{_{\hh_c}}(0)$
has a finite increasing filtration by $\hh_c$-submodules:
$0=F_0\subset F_1\subset\ldots\subset F_l=P,$
such that $F_i/F_{i-1}\simeq M(\sigma_i),$
for some standard modules $M(\sigma_i)\,,\, \sigma_i\in\irrep(W)\,,\,
i=1,\ldots,l.$
Furthermore, for any $\tau,\sigma\in\irrep(W),$
Guay [Gu], and Opdam-Rouquier [OR], proved
the following analogue of Brauer-Bernstein-Gelfand-Gelfand
type reciprocity formula
$\,[P(\tau)\,:\,M(\sigma)]=[M(\sigma)\,:\,L(\tau)]$,  see [BGG].
This formula implies, in particular, that the category
$\oo_{_{\hh_c}}(0)$ is semisimple if and only if
all the standard modules  $M(\tau)\,,\,\tau\in\irrep(W),$
are simple, i.e. if and only if
$M(\tau)=L(\tau)\,,\,\forall \tau\in\irrep(W)$.
It would be interesting to obtain similar results for the
category $\oo_{_{\hh_c}}(\bar\lambda),$ with $\bar\lambda\neq 0$.
\end{remark}

\section{Harish-Chandra $\hh_c$-bimodules}
\setcounter{equation}{0}

One of the goals of this section is to prove 

\begin{theorem}\label{th_simple1} If $c\in \ccr,$
then $\hh_c$ is a simple algebra. 
\end{theorem}

The proof of the Theorem will be based on the concept
of a Harish-Chandra bimodule which we now introduce.

Let $V$ be an $\hh_c$-bimodule, i.e., a 
left $\hh_c\otimes \hh_c^{op}$-module.
For any $x\in\hh_c$ we have an adjoint $x$-action on
$V$ given by $\ad x: v\mapsto xv-vx$.

\begin{definition}\label{HC_def} A finitely generated $\hh_c\otimes
(\hh_c)^{op}$-module $V$ is called
a Harish-Chandra  bimodule  if,
for any $x\in \C[\h]^W$ or $x\in \C[\h^*]^W$, the $\ad x$-action
on $V$ is locally nilpotent.
\end{definition}

Here are the first elementary results about Harish-Chandra  bimodules.
\begin{lemma}\label{straight}
\vi Harish-Chandra  bimodules form a full abelian subcategory in the category of
$\hh_c$-bimodules; this subcategory is stable under extensions.

\vii Any Harish-Chandra  bimodule is finitely-generated over the subalgebra
$\C[\h]^W \otimes \C[\h^*]^W\subset \hh_c\otimes
(\hh_c)^{op},$ with $\C[\h^*]^W$ 
 acting on the right and
$\C[\h]^W$ on the left. Similarly, it 
 is also finitely-generated over the subalgebra
$\C[\h^*]^W \otimes\C[\h]^W\subset \hh_c\otimes
(\hh_c)^{op}$;

\viii {Any Harish-Chandra  bimodule is  finitely-generated
as a left, resp. right,~$\hh_c$~-~module;}

\iv For any maximal ideal $J\subset \C[\h^*]$ and any
Harish-Chandra  bimodule $V$, we have
$V/V\cd J\in \oo_{\hh_c},$
as a left module.

${\sf{(v)}}\;$  The algebra $\hh_c$ is a Harish-Chandra  $\hh_c$-bimodule.
\end{lemma}

\begin{proof} The algebra $\hh_c$ has an increasing
filtration $F_\bullet\hh_c$
such that $\grd(\hh_c)=\C[\h\times\h^*]\#W$, see [EG]. 
Hence  $\hh_c$ is a Noetherian algebra, for this 
is clearly true for the  algebra 
$\grd(\hh_c)=\C[\h\times\h^*]\#W$. Part (i) of the Lemma follows.

To prove (ii), observe that the filtration on $\hh_c$ gives rise to the
tensor product filtration:
$(\hh_c\otimes\hh_c^{op})_p
:= \sum_{p=i+j}\,F_i(\hh_c) \otimes F_j(\hh_c^{op})$,
on the algebra 
$\hh_c\otimes\hh_c^{op}$.
Choose a finite-dimensional subspace $V_0\subset V$
generating $V$ as an $\hh_c\otimes\hh_c^{op}$-module
and, for each $p\geq 0$, put
$V_p:=(\hh_c\otimes\hh_c^{op})_p\cdot V_0.$
This is
a good filtration on $V$, and
one can regard $\grd(V)$, the associated graded space,
 as a $W\times W$-equivariant
finitely-generated module over
 $\C[\h_1\times\h_2\times\h_1^*\times\h_2^*],
$ where the subscripts `$1,2$' indicate the corresponding copy of
 $\h$.
Then,
for any homogeneous element $P\in \C[\h]^W$,
the $\ad P$-action on $V$ corresponds to
the action of the element
$P^{\mathtt{left}}\otimes 1 - 1\otimes P^{\mathtt{right}}
\in \C[\h_1]^W\otimes \C[\h_2]^W$
on $\grd(V)$. The (local) nilpotency of the
 $\ad P$-action on a Harish-Chandra module thus
implies that, if $P$ is homogeneous of
degree $>0$, then the action of the element
$P^{\mathtt{left}}\otimes 1 - 1\otimes P^{\mathtt{right}}$
on $\grd V$ is nilpotent.
Hence, the support of $\grd V$,
viewed as a  $\C[\h_1\times\h_2\times\h_1^*\times\h_2^*]$-module,
is contained in the zero locus of the polynomial
$P^{\mathtt{left}}\otimes 1 - 1\otimes P^{\mathtt{right}}
\in \C[\h_1\times\h_2\times\h_1^*\times\h_2^*]$.
This way, for any Harish-Chandra module $V$,
we deduce  the following upper bound on the 
set $\supp(\grd V)\subset\h_1\times\h_2\times\h_1^*\times\h_2^*$:
\begin{equation}\label{graph}
\supp(\grd V) \;\subset\;
\Bigl(\cup_{w\in W}\;\text{Graph}(w)\Bigr)
\times\Bigl(\cup_{y\in W}\;\text{Graph}(y)\Bigr)\,,
\end{equation}
where $\text{Graph}(w)$ denotes the graph in $\h_1\times\h_2$,
resp. $\h_1^*\times\h_2^*$, of the $w$-action
map: $x\mapsto w(x)$.

In particular, the restriction of the
composite map:
$\h_1\times\h_2\times\h_1^*\times\h_2^*
\to\h_1\times\h_1^* \to \h_1/W\times\h_1^*/W$
to $\supp(\grd V)$ is a finite map.
Therefore, $\grd(V)$ is finitely generated over
the subalgebra $\C[\h]^W \otimes \C[\h^*]^W$,
hence the same holds for $V$ itself.
Part (ii) follows.

To prove (iii) we observe that (\ref{graph})
implies also that the restriction  to $\supp(\grd V)$ of the
corresponding projection: $\h_1\times\h_2\times\h_1^*\times\h_2^*
\to\h_1\times\h_2^*$ is a
 finite map.
We deduce similarly that $\grd V$ is finitely generated 
with respect to the left action of the algebra
$\grd(\hh_c)$. This forces
$V$ to be finitely generated left $\hh_c$-module,
and (iii) follows.

Clearly (ii) implies (iv). To prove (v), endow $\hh_c$
with an increasing filtration by assigning
elements of $\h$ filtration degree 1, and elements of
$\h^*$ and $\C W$  filtration degree 0
(this is not the standard filtration $F_\bullet\hh_c$ used earlier).
 The defining relations
(\ref{defrel}) of the Cherednik algebra show that
commuting with an element of $\C[\h^*]^W\subset \hh_c$
decreases  filtration degree at least by 1. Hence,
the adjoint action of $\C[\h^*]^W$ is locally nilpotent.
The other condition is proved similarly.
\end{proof}

\begin{proposition}\label{loc_triv} 
If $c\in\ccr,$ then any Harish-Chandra
$\hh_c$-bimodule is a finite rank projective module over 
the algebra
$\C[\h]^W\otimes \C[\h^*]^W\subset \hh_c\otimes\hh_c^{op}.$
\end{proposition}
\begin{proof} Recall that a finitely generated module $F$ over $A$,
a commutative finitely generated algebra  without zero-divisors,
 is projective
if and only if, for  all maximal ideals
${\mathfrak{m}}\subset A$, the geometric fibers $
F/{\mathfrak{m}}\cd F$ all
have the same (finite) dimension over $\C$. Thus, we arrive at the following
\smallskip

\noindent
{\bf Observation:}
{\it Let $X,Y$ be irreducible
 algebraic varieties, and ${\mathcal{F}}$
 a coherent sheaf on $X\times
Y.$  Suppose that for any $x\in X,$ the sheaf 
${\mathcal{F}}_x ={\mathcal{F}}|_{x\times Y}$
is locally free of finite rank,
 and for any $y\in Y$ the 
sheaf ${\mathcal{F}}_y={\mathcal{F}}|_{X\times y}$ is
 locally free of finite rank.
  Then ${\mathcal{F}}$ is itself locally free of finite rank.}
\smallskip

Now let $V$ be a  Harish-Chandra $\hh_c\otimes\hh_c^{op}$-module.
We view $V$ as a module over
$\C[\h]^W\otimes \C[\h^*]^W,$ with $\C[\h^*]$ acting on the right and
$\C[\h]$ on the left. By the previous Lemma, this module 
is finitely generated. Therefore,
by the Observation above, to prove the Proposition
 we only have to
show that, for any   maximal ideal
$J \subset \C[\h^*]^W$, the left $\C[\h]^W$-module 
$V/V\cd J$ is  locally free 
  of finite rank
(and a similar result for a right  $\C[\h^*]^W$-module 
$V/I\cd V\,,\, I\subset \C[\h]^W$).
Thus,  let $J \subset \C[\h^*]^W$ be the maximal ideal,
corresponding to a point $\bar \lambda\in\h^*/W$.
Then $V/V\cd J\in \oo_{_{\hh_c}}(\bar \lambda)$.
This module has a finite
composition series
whose successive quotients are
standard modules, since $c\in\ccr$. But all standard modules,
$M(\lambda,\tau),$
 are $\C[\h]^W$-free,
and the statement follows.
\end{proof}

Let $V\mr$ denote the localization of a  Harish-Chandra
$\hh_c$-bimodule
$V$ with respect
to the {\it left} action of the subalgebra $\C[\h] \subset \hh_c$.
We note that  $V\mr$ coincides also with  the localization of
$V$  with respect
to the {\it right} action of the subalgebra $\C[\h]$. To see this,
observe that in either case, one 
may  replace localizing with respect to $\C[\h]$ by localizing with respect
to
the smaller algebra
$\C[\h]^W$. But it follows easily from Definition
\ref{HC_def} that left and right  localizations of $V$  with respect
to $\C[\h]^W$ coincide. 

Proposition \ref{loc_triv} yields readily

\begin{corollary}\label{key} Let $c\in\ccr,$ and let $V$ 
be a  Harish-Chandra  $\hh_c$-bimodule   such that
$V\mr=0$.
 Then $V=0$.\qquad\qed
\end{corollary}

\noindent
{\sl Proof of Theorem \ref{th_simple1}.\,} Let $I\subset \hh_c$ be
a nonzero two-sided ideal. Then $V:= \hh_c/I$
is, by  Lemma \ref{straight}(v), a Harish-Chandra  bimodule.
Moreover, since $\hh_c$ is a free left $\C[\h]$-module,
and any submodule of a free $\C[\h]$-module is torsion-free,
it follows that $I\mr\neq 0$, where
$I\mr$ denotes the left (equivalently, right) localization of $I$
viewed as a Harish-Chandra 
bimodule.
We claim that $V\mr=
(\hh_c\mr)\big/(I\mr)=0$.
Indeed, as we have mentioned in \S2, there is an
algebra isomorphism $\hh_c\mr\simeq\dd(\hr)\#W$.
The algebra $\dd(\hr)\#W$ is known to be simple [Mo].
Hence,  $I\mr\neq 0$ implies $\hh_c\mr=I\mr$. Thus, 
$(\hh_c/I)\mr=0$.
Now, Corollary \ref{key} yields: 
$0=V= \hh_c/I$, and $\hh_c=I$.\qed
\smallskip

   From Proposition \ref{loc_triv} one also obtains (see also [Gi]) the following important

\begin{corollary}\label{finite}
If $c\in\ccr$, then every Harish-Chandra module has finite length.
\end{corollary}
\begin{proof} Let $V=V^0 \underset{^{\neq}}{\supset} V^1
\underset{^{\neq}}{\supset}\ldots $
be a strictly decreasing chain of Harish-Chandra sub-bimodules 
of a Harish-Chandra bimodule $V$. View it as a chain of
locally-free $\C[\h]^W\otimes \C[\h^*]^W$-modules.
Since each quotient $V^i/V^{i+1}$ is a non-zero
locally-free $\C[\h]^W\otimes \C[\h^*]^W$-module
(by Proposition \ref{loc_triv}),
the ranks of the locally-free $\C[\h]^W\otimes \C[\h^*]^W$-modules
$V^i$ form a strictly decreasing sequence:
$\rk_{_{\C[\h]^W\otimes \C[\h^*]^W}}V^0 \,>\, 
\rk_{_{\C[\h]^W\otimes \C[\h^*]^W}}V^1\,>\,\ldots.$
It follows that $V^N=0$, for $N$ big enough.
\end{proof}

Let  $\XX\in \C[\h]$ denote the squared norm function:
$x\mapsto (x,x)$, cf. \S2, and let $\YY\in \C[\h^*]$  denote a
similar function on $\h^*$ relative to the inner product transported
from $\h$ via the isomorphism $\h^*\simeq \h$.
 A straightforward computation based
on formula (\ref{sh}) shows that 
the elements  $\langle \XX,\sh,\YY\rangle$ form
 an
$\sll2$-triple in
 the algebra $\hh_c$.

\begin{proposition}\label{sltwo}
Any Harish-Chandra bimodule $V$ breaks,
under the adjoint action of $\langle \XX,\sh,\YY\rangle$,
 into a direct sum of finite-dimensional
$\sll2$-modules. In particular, the
$\ad\sh$-action on $V$ is diagonalizable.
\end{proposition}

\begin{remark}
We see that, in the situation of the Proposition
 the adjoint $\sll2$-action exponentiates
to an algebraic $SL_2$-action on $V$.
In particular, the automorphism
corresponding to the action of
the element $\left(\begin{matrix}0 &1\\-1& 0
\end{matrix}\right)\in SL_2$ may be thought of as a sort of
{\it Fourier transform}. In the special case $V=\hh_c,$
this $SL_2$-action on  $\hh_c$ by
algebra automorphisms has been studied in [EG].
\end{remark}

The $SL_2$-action on  $\hh_c$ combined with
$SL_2$-action on a Harish-Chandra bimodule
gives way to  the following more symmetric
point of view on  Harish-Chandra bimodules.

\begin{corollary}\label{alternHC}
 A finitely generated $\hh_c\otimes
(\hh_c)^{op}$-module $V$ is 
a Harish-Chandra  bimodule  if and only if,
for any $g\in SL_2(\C)$, the adjoint action
on $V$ of the subalgebra $g(\C[\h]^W)\subset \hh_c$
 is locally nilpotent.\qed
\end{corollary}

Write $U(-)$ for the enveloping algebra 
of a Lie algebra.
We need the following lemma in representation theory
of a complex semisimple Lie
algebra  ${\mathfrak{g}}$ with
triangular decomposition:
${\mathfrak{g}}={\mathfrak{n}}_+
+\h+{\mathfrak{n}}_-,\,$ 
which is a special case of a more general result due to Kac [K].
\begin{lemma}\label{lemma_sltwo1}
Let $V$ be a $U({\mathfrak{g}})$-module 
such that, $\dim U({\mathfrak{n}}_+)v<\infty$ and
$\dim U({\mathfrak{n}}_-)v<\infty$,
for any $v \in V.$
Then $V$ is a locally finite $U({\mathfrak{g}})$-module.\qed
\end{lemma}

\noindent
{\sl Proof of Proposition \ref{sltwo}.\,\,}
We apply Lemma \ref{lemma_sltwo1} to the Lie
algebra ${\mathfrak{g}}=\sll2$ acting
on a Harish-Chandra $\hh_c$-bimodule $V$ via
the adjoint action of
 the $\sll2$-triple  $\langle \XX,\sh,\YY\rangle$.
The conditions of the Lemma clearly hold in this case.
We conclude that the adjoint action on $V$ of our
$\sll2$-triple is locally finite. The result follows.
\qed\smallskip

\section{The spherical subalgebra $\e\hh_c\e$}
\setcounter{equation}{0}
We begin with the following simple ring-theoretic result. 

\begin{lemma}\label{mor_simple}
If $A$ is a simple unital algebra and $\e$ is an idempotent
  of $A$. Then the algebra
$\e A\e$ is simple and Morita equivalent to $A$.
\end{lemma}

\begin{proof} Let $I$ be a nonzero two-sided ideal in $\e A\e$. 
Then $AIA=A$ since $A$ is simple, and hence $\e AIA\e=\e A\e$. 
But $I=\e I\e$, so $\e A\e I\e A\e =\e A\e$, which implies $I=\e
  A\e$. Hence $\e A\e$ is simple. The second claim follows from
  the general fact that if $A\e A=A$ then $A$ is Morita
  equivalent to $\e A\e$ (see e.g. \cite{MR}).
\end{proof}

Now, let $W$ be our Coxeter group, and
 $\e=\frac{1}{|W|}\sum_{w\in W}\,w$  the symmetrizer 
 in
$\C{W}.$ 
In [EG] we have introduced the subalgebra $\e\hh_c\e\subset \hh_c,$
called the spherical 
subalgebra. 
    From Lemma \ref{mor_simple} and Theorem \ref{th_simple1} we get
\begin{corollary} \label{cor_simple1} 
If $c\in \ccr,$ then $\e\hh_c\e$ is a
simple algebra, Morita equivalent to $\hh_c$.\qed 
\end{corollary}

One defines the category $\oo_{\e\hh_c\e}$ to be the abelian category
formed by 
the finitely-generated  $\oo_{\e\hh_c\e}$-modules $M$, such that
the action on $M$ of the subalgebra $\C[\h^*]^W\subset \e\hh_c\e$
is locally finite.

    From Corollary \ref{cor_simple1} we deduce
\begin{proposition}\label{O=O} For any $c\in\ccr,$
the functor $M\mapsto \e\cdot M$ gives an equivalence:
$\oo_{\hh_c} \stackrel{_\sim}{\to}\oo_{\e\hh_c\e}.$
\qed
\end{proposition}
 
Put $M_\e(\lambda,\tau):=\e M(\lambda,\tau)$.
This $\e\hh_c\e$-module corresponds to the standard
$\hh_c$-module via Morita equivalence.
The modules $M_\e(\lambda,\tau)$
 will be referred to as standard $\e\hh_c\e$-modules.
These modules can be explicitly described as follows.
Given $\lambda\in\h^*$, one has an algebra
homomorphism: $\C[\h^*]\#W_\lambda \onto \C{W_\lambda}$ sending
$P\in \C[\h^*]$ to $P(\lambda)$
and identical on $W_\lambda$. We pull-back the left
regular representation of $W_\lambda$ to $\C[\h^*]\#W_\lambda$ 
via this homomorphism, and form
an induced $\hh_c$-module $N(\lambda)
:=\text{Ind}_{\C[\h^*]\#W_\lambda}^{\hh_c}\C{W_\lambda}$.
Observe further that the multiplication-action
of the group $W_\lambda$ on the group algebra $\C{W_\lambda}$
on the right commutes with the left
$\C[\h^*]\#W$-module structure, hence, gives rise
to a right $W_\lambda$-action  on $N(\lambda)
:=\text{Ind}_{\C[\h^*]\#W_\lambda}^{\hh_c}\C{W_\lambda}$.
We leave to the reader to prove
\begin{lemma}
The  $\e\hh_c\e$-module
$M_\e(\lambda,\tau)$ is isomorphic
to the $\tau$-isotypic component (with respect to the
right   $W_\lambda$-action) in the left  $\e\hh_c\e$-module
$\e\cd N(\lambda)$.\qed
\end{lemma}

There is also a  notion of Harish-Chandra  bimodule over
the spherical algebra $\e\hh_c\e$ instead of $\hh_c$,
defined by the same 
 Definition \ref{HC_def}.
An obvious analogue of Lemma \ref{straight} holds
for Harish-Chandra  bimodules over  $\e\hh_c\e$.
Further note that, for any  $c\in\ccr,$
all standard $\ehe$-modules $M_\e(\lambda,\tau)$
 are simple, due to  Morita equivalence
of Proposition \ref{O=O}.
Repeating  the proof of Proposition \ref{loc_triv}
 we deduce

\begin{corollary} \label{cor_eHe}
If $c\in\ccr,$ then any Harish-Chandra
$\e\hh_c\e$-bimodule $V$ is a finite rank projective module over 
$\C[\h]^W\otimes \C[\h^*]^W\subset \ehe\otimes(\ehe)^{op}.$
In particular,
$V\mrw=0\;$ $\Longrightarrow\;
V=0$.\qed
\end{corollary}

\begin{remark} It has been explained in
[EG,\S8] that, in the case of a root system of type ${\mathbf{A_1}}$,
the algebra $\ehe$ is isomorphic to
a quotient of the algebra $U(\sll2)$.
Specifically,
if $C=ef+fe+h^2/2\in U(\sll2)$ is the standard quadratic Casimir,
then $\ehe\simeq U(\sll2)/I_c,$
where $I_c$ is the  ideal in
 $U(\sll2)$ generated by the
central element $C- \frac{1}{2}(c+\frac{3}{2})\cdot
(c-\frac{1}{2}).$
Proposition \ref{sltwo} and Lemma \ref{lemma_sltwo1} imply 
that a Harish-Chandra  $\ehe$-bimodule
is in this case  nothing but a Harish-Chandra bimodule
 over $U(\sll2)$ (with equal\footnote{non-equal 
case will be treated in \S8.} left and right central characters),
 in the sense of representation theory of semisimple
Lie algebras, see  [Di]. Further,
the isomorphism $\ehe\simeq U(\sll2)/I_c$ suggests 
to view 
the standard module $M(\lambda,\sigma)\,,\,\lambda \neq 0,$
as a generalized
{\it Whittaker module}.\qed
\end{remark}
\medskip

We now prove
the following important result that reduces, in the special case
$c=0$ and $\e\hh_c\e=\dd(\h)^W,$ to Theorem 5 in
the Levasseur-Stafford paper [LS] (for Coxeter groups). 

\begin{theorem}[Levasseur-Stafford theorem for $\e\hh_c\e$]\label{th_simple2}
If the algebra $\e\hh_c\e$ is simple for some $c\in\CC,$ then it is generated 
as an algebra by $\C[\h]^W\!\!\cd\e$ and $\C[\h^*]^W\!\!\cd\e$. 
\end{theorem}

To prove the theorem we begin with
the following {\it weak localized} version
\footnote{the `strong' version of the Lemma,
i.e., the non-localized one where $\hreg$ is replaced
by $\h$, is known, due to Wallach [Wa],
 to be {\it false}, in general.}
of a Poisson 
analogue
of Theorem 5 of \cite{LS}.

\begin{lemma}\label{poisgen} 
The algebra $\C[\hreg\times \h^*]^W$ of regular functions 
on the symplectic manifold $(\hreg\times \h^*)/W$ is generated 
by $\C[\hreg]^W$ and $\C[\h^*]^W$ as a Poisson algebra. 
\end{lemma}

\begin{proof} (of Lemma \ref{poisgen}) Denote by $A$ the Poisson subalgebra of 
$\C[\hreg\times \h^*]^W,$ Poisson generated by $\C[\hreg]^W$ and $\C[\h^*]^W$. 
By Hilbert-Noether lemma,
 $A$ is a finitely generated commutative algebra without zero divisors.
Let $Y=(\hreg\times \h^*)/W$, and $Y'=\text{Spec}(A)$. 
Then $Y'$ is an irreducible affine 
algebraic variety, and the natural map 
$f:Y\to Y'$ is a finite map. In particular, $f$ is 
surjective. 

Let us show that $f$ is injective, i.e. that 
the algebra $A$ separates points of $Y$. For this, consider the 
natural map $h: Y\to \hreg\!/W\times \h^*/W$. 
Fix $(q,p)\in \hreg\times \h^*$, let
$W_p\subset W$ be the stabilizer of $p$, and 
$(\bar q,\bar p)$ be the projection of $(q,p)$ to 
$\hreg\!/W\times \h^*/W$. The set $h^{-1}(\bar q,\bar p)$ 
consists of equivalence classes of points $(q,gp)$, 
where $g\in W/W_p$. Thus, to prove the injectivity of $f$, it suffices to find 
an element $a\in A$ such that the numbers $a(q,gp)$ 
are distinct for all $g\in W/W_p$.

Let $z\in \h^*$ be a vector such that the numbers 
$\langle z,gp\rangle$ are distinct for all 
$g\in W/W_p$ (this is the case for a generic vector, so such $z$ exists). 
Let $b\in \C[\hreg]^W$ be a function such that 
$(db)(q)=z$ (such exists since the stabilizer of $q$ is trivial). 
Set $a=\lbrace b,\YY/2\rbrace,$
where  $\YY$ is the squared norm function on $\h^*$.
Then $a(q,gp)=\langle z,gp\rangle$, 
which are distinct, as desired. 

Thus, $f$ is bijective. To conclude the proof, we need to show that 
moreover $f$ is an isomorphism. Since $Y=(\hreg\times\h^*)/W$ is smooth,
 it is enough to find,
for any $(q_0,p_0)\in Y$, a set of elements of $A$ which form a system 
of local coordinates near $(q_0,p_0)$. 

Let $a_1,...,a_r\in \C[\hreg]^W$ be functions such that 
the vectors $da_i(q)$ form a basis of $\h^*$. 
Let $b_i=\lbrace a_i,\YY/2\rbrace $.
 It is easy to see 
by computing the Jacobian that the functions 
$a_i$ and $b_i$ form a system of coordinates on $Y$ near 
$(q_0,p_0)$. The lemma is proved. 
\end{proof}

\begin{lemma}[\cite{LS}, Lemma 9]\label{l_simple2} 
Let $R \subset S$ be two (not necessarily commutative)
Noetherian domains, 
such that $S$ is simple, and
finite both as a left and a right $R$-module. Then, if $R$ and $S$ have the 
same skew field of fractions, 
then $R=S$.\qed
\end{lemma}

\noindent
{\sl Proof of Theorem \ref{th_simple2}.}\,\, 
The proof is analogous to the proof 
of Theorem 5 in \cite{LS}.
We apply Lemma \ref{l_simple2}  to the situation when $S=\e\hh_c\e$, and $R$ is
the subalgebra generated by 
$\C[\h]^W\!\!\cd\e$ and $\C[\h^*]^W\!\!\cd\e.$ 
It remains to  check the conditions of the Lemma. 

We know that
$\grd(S)=\C[x,y]^W$ is a finite rank $\C[\h^*]^W\otimes \C[\h]^W$ module. 
Hence  $\grd(S)$, and its submodule
$\grd(R)$, are both Noetherian. 
Thus, so are $R$ and $S$. It is also clear that $R$ and $S$ are
domains (have no zero divisors).
By  assumption
 $S$ is a simple algebra. Further, as we have seen above,
$\grd(S)$ is a finite $\grd(R)$ module. Hence, $S$ is a finite 
left and right $R$-module. 

Finally, let us check that the quotient fields of $R$ and $S$ coincide.
 It suffices 
to check that $R\mrw=S\mrw$. Consider the filtration on 
$S\mrw$ defined by ``degree in $y$'' (i.e. order of differential
operators in the Dunkl-Cherednik realization).
This filtration is nonnegative, so to check that 
$R\mrw=S\mrw$, it suffices to check that $\grd(R\mrw)=\grd(S\mrw)$. 
The algebras $\grd(R\mrw)$ and $\grd(S\mrw)$ are Poisson algebras. 
Moreover, it is easy to see that $\grd(S\mrw)=\C[\hreg\times \h^*]^W$, 
while $\grd(R\mrw)$ contains $\C[\hreg]^W$ and $\C[\h^*]^W$. 
Hence, it follows from Lemma \ref{poisgen} that $R\mrw=S\mrw$. 
Theorem \ref{th_simple2} is proved.
\qed\smallskip

We observe next that the   elements $\langle \XX\e,\sh\e,\YY\e\rangle$
form an $\sll2$-triple in the algebra $\e\hh_c\e$,
and an analogue of Proposition \ref{sltwo} holds for
Harish-Chandra $\ehe$-bimodules.
In particular, the adjoint
action of $\sh$ on a Harish-Chandra $\hh_c$-bimodule,
resp., the adjoint
$\sh\e$-action on a Harish-Chandra 
$\ehe$-bimodule, $V$ gives a grading
$V=\oplus_{k\in\Z}\,V(k),$ where
 $V(k)=\{v\in V\mid \ad\sh(v)=k\cdot v\}.$ Therefore,
the structure theory of finite dimentional representations
of $\sll2$ implies that, for any $k\geq 0$ the operator
$(\ad\YY)^k$ induces an isomorphism
$\,(\ad\YY)^k: \Ker(\ad\XX) \cap V(-k) \iso \Ker(\ad\YY)\cap V(k).$
In particular, for $V=\hh_c$, resp. $V=\ehe$, 
the map $(\ad\YY)^k$ takes the space $\C^k[\h]\subset \hh_c$
into the space formed by elements of filtration degree $\leq k$ and
$\ad\sh$-weight $(-k)$, which is exactly the space $\C^k[\h^*]\subset \hh_c$.
Thus,
we deduce
\begin{corollary}\label{Lefschetz}
For any $k\geq 0$ the operator
$(\ad\YY)^k$ induces isomorphisms
$$\C^k[\h] \iso \C^k[\h^*]\quad
\text{and}\quad \C^k[\h]^W\!\!\cd\e \iso \C^k[\h^*]^W\!\!\cd\e\quad,
\quad\C[\h]\,,\,\C[\h^*]\subset \hh_c\,.\quad\qed
$$
\end{corollary}

The natural $\C^*$-action on $\hreg$ by dilations
gives a $\Z$-grading on $\dd(\hreg)$ by the {\it weight}
of the $\C^*$-action.
Also, there is a standard filtration on $\dd(\hreg)$ by 
the {\it order} of differential operators, and
we let $\dd(\hreg)_\minus$ denote the subset in $\dd(\hreg)$
spanned by all homogenious differential operators
$u\in \dd(\hreg)$, such that $\text{\it weight}(u) +
\text{\it order}(u)\leq 0.$ Clearly,
$\dd(\hreg)_\minus$ is an associative subalgebra in $\dd(\hreg).$

Next, introduce
the  rational Calogero-Moser operator
$\LL_c:
 = \Delta - \sum_{\alpha \in R_+} \frac{2c_{\alpha}}
{\alpha}\, \partial_{\alpha},\,$
where $\Delta$ denotes the 2-d order Laplacian on $\h$.
It is clear that $\LL_c\in \dd(\hreg)_\minus$.
Following [O],
we let $\cc_c$ be the {\it centralizer} of $\LL_c$
in the algebra $\dd(\hreg)_\minus^W=\dd(\hreg)^W\cap \dd(\hreg)_\minus,$
and let $\bb_c$ denote the subalgebra in $\dd(\hreg)^W$
generated by  $\cc_c$ and by $\C[\h]^W$, the
 subalgebra of $W$-invariant polynomials.

 Recall that in [EG] we have used
 the Dunkl representation of the algebra $\hh_c$
to construct an injective algebra homomorphism
$\Theta\spher: \ehe \into \dd(\hreg)^W$ that
sends the element $\YY\e$ to $\LL_c$.
The following result is a strengthening of [EG], Theorem 4.8.
\begin{proposition}\label{spher} For any $c\in\ccr$, the image of the
map: $\ehe \into
\dd(\hreg)^W$ equals $\bb_c$. Thus, one has an algebra isomorphism
$\,\Theta\spher: \ehe\iso \bb_c$.
\end{proposition}

\begin{proof} It is obvious that the map
$\Theta\spher$ takes the subalgebra $\C[\h]^W\cd\e\subset\ehe$ into
$\C[\h]^W\subset\dd(\hreg)^W$,
and  takes the subalgebra $\C[\h^*]^W\!\!\cd\e\subset\ehe$ into
$\cc_c\subset\dd(\hreg)^W$. Hence, $\Theta\spher(\ehe)\supset\bb_c$.
But, for $c\in\ccr$, the
algebra $\ehe$ is generated by  $\C[\h]^W\!\!\cd\e$ and
$\C[\h^*]^W\!\!\cd\e$, by the Levasseur-Stafford theorem.
Thus, the algebra $\Theta\spher(\ehe)$ is 
 generated by  $\C[\h]^W$ and $\cc_c$.
\end{proof}

Let $\varepsilon: W \to \Z/2\Z$ be 
a multiplicative  character 
of the group $W$
(if $W$ is the Weyl group of a simply laced
root system, then
there is only one such character $\varepsilon(w)=\det_\h(w)$). 
 We consider the
set:
$
R_{\varepsilon} := \{ \alpha \in R \mid \varepsilon(s_{\alpha}) = -1\}$.
Denote by $\mathbf{1}_\varepsilon$ the characteristic 
function of the subset $R_{\varepsilon}\subset R, $
that is $\mathbf{1}_\varepsilon(\alpha) := 1$ if $\alpha 
\in R_{\varepsilon},$ and $0$ otherwise. 
Clearly, $\mathbf{1}_\varepsilon\in\CC$.
Also, we put $\e_{\varepsilon}:= \frac{1}{|W|}\sum_{w\in W}\,
\varepsilon(w)\cd w, $ a central idempotent in $\C{W}$.

\begin{proposition}\label{H_minus}
There is an algebra isomorphism:
$ \e_{\varepsilon}\hh_{c}\e_{\varepsilon}
 \simeq \e\hh_{c -\mathbf{1}_\varepsilon}\e,$
for any $c\in\ccr.$
\end{proposition}

\noindent
{\sl Proof.\;}
Since $c\in\ccr$,
Corollary
\ref{Lefschetz} implies that the algebra  $\e\hh_{c}\e$
is generated by  $\C[\h]^W\!\cdot\e$ and by the
element $\YY\e$. Hence, 
the algebra  $\Theta\spher(\ehe)$ is generated by  $\C[\h]^W$
and by the Calogero-Moser  operator
$\Theta\spher(\YY\e)=\LL_c$.
But $\Theta\spher(\ehe)=\bb_c,$
by Proposition \ref{spher}.
Hence, $\bb_c$ is generated by  $\C[\h]^W$
and $\LL_c$ as an algebra.

Now, given a character $\varepsilon: W\to \{\pm 1\}$, we define a 
polynomial 
$$
\delta_{\varepsilon} = \prod\nolimits_{\alpha \in R_{+}}\,
\alpha^{\mathbf{1}_\varepsilon(\alpha)} =
\prod\nolimits_{\alpha \in R_\varepsilon \cap R_{+}}\, \alpha \;\;\in \C[\h]\,.
$$
To emphasize dependence on the parameter `$c$',
write $\YY={\mathbf{y}^\mathbf{2}_c}\in \hh_c$,
and $\Theta_c: \hh_c\to \dd(\hreg)\#W,$ for the Dunkl representation.
Then, an easy  calculation, see [He], yields the following equation: 
$$
\delta_{\varepsilon}^{-1}\ccirc\Theta_c({\mathbf{y}^\mathbf{2}_c}\e_{\varepsilon})\ccirc
 \delta_{\varepsilon}\,=\, 
\Theta_{c - \mathbf{1}_\varepsilon}
({\mathbf{y}^\mathbf{2}_{c -\mathbf{1}_\varepsilon}}\e)\qquad
\text{holds in }\enspace \dd(\hreg)\#W\,.
$$
Since $e^{2\pi i c}= e^{2\pi i (c - \mathbf{1}_\varepsilon)}\in\Reg$,
the algebra $ \bb_{c - \mathbf{1}_\varepsilon}\simeq
\e\hh_{c - \mathbf{1}_\varepsilon}\e$ is simple. Hence,
it is 
generated by $\C[\h]^W$ and $\LL_{c -
\mathbf{1}_\varepsilon},$
and we obtain algebra isomorphisms:
$$
\e\hh_{c - \mathbf{1}_\varepsilon}\e\,\simeq\,
\bb_{c - \mathbf{1}_\varepsilon}\,=\,
\delta_{\varepsilon}^{-1}\ccirc\Theta_c(\e_{\varepsilon}\hh_c\e_{\varepsilon})\ccirc
 \delta_{\varepsilon}
\,\simeq\,\Theta_c(\e_{\varepsilon}\hh_c\e_{\varepsilon})
\,\simeq\,\e_{\varepsilon}\hh_c\e_{\varepsilon}\,.\qquad\square
$$

\section{A trace on the Cherednik algebra}
\setcounter{equation}{0}

In this section we compute a trace on    Cherednik algebras 
 of type ${\mathbf{A_{n-1}}}$, that is on the  algebra
$\hh_c$ corresponding
to
$W=S_n$, the Symmetric group.
Write $[\hh_c,\hh_c]$ for the $\C$-vector subspace
in $\hh_c$ spanned by all the commutators $[a,b]\,,\,a,b\in \hh_c$.
 Let $\Tr(g)$ denote
 the image of $g\in \C[S_n]$ under the composite map
$\Tr:\,\C[S_n]\into \hh_c \onto \hh_c/[\hh_c,\hh_c].$ 
Given a permutation $g\in S_n$, write $\cyc(g)$ for the number of
cycles in $g$.

\begin{proposition}\label{trace1} If $c\neq 0$ then,
for any permutation $g\in S_n,$ in $\hh_c/[\hh_c,\hh_c]$ we have: 
$$\Tr(g)=(1/nc)^{n-\cyc(g)}\cdot\Tr(1)\;\;.$$
\end{proposition}
\begin{proof}
 The proof proceeds by induction in $m=n-\cyc(g). $
For $m=0,$ the result is obvious. Suppose the result is known for $m=p,$ and 
let us establish it for $p+1.$ Let $g$ have $n-p-1$ cycles. Choose 
two indices $1\leq i,j\leq n,$ such that $g(i)=j.$ Consider the permutation 
$\sigma=g\cd s_{ij}.$ Then $\sigma(j)=j,$ and $\sigma$ has $n-p$
 cycles, so for $\sigma$ the statement is known. 

Let $x_l,y_l$ be the elements in the first and second copy of $\h\cong
\h^*\cong \C^{n-1}\sset \C^n,$ respectively, 
corresponding to the vector $(-\frac{1}{n},...,1-\frac{1}{n},...,-\frac{1}{n}),$ where 
$1-\frac{1}{n}$ stands on the $l$-th place. 
The commutation relations in $\hh_c$, see \eqref{defrel}, imply readily that,
for any $l$ different from $m,$ one has
$[y_m,x_l]=c\cd s_{ml}-\frac{1}{n}.$
Thus, we find
$
[\sigma\cd y_i,x_j]=\sigma\cd [y_i,x_j]=c\cd g-\frac{1}{n} \sigma\,. $
Hence,
 $\Tr(g)=\frac{1}{nc}\cd \Tr(\sigma),$ and we are done. 
\end{proof}

Recall that the (isomorphism classes of)
irreducible representations of the Symmetric group
are naturally labelled by Young diagrams.
Write  $Y(\tau)$ for  the Young diagram (with $n$ boxes) corresponding
to
a simple $S_n$-representation
 $\tau$, and  let $s_{\tau}$  denote
the Schur function associated to $Y(\tau).$
Given a box $u\in Y(\tau)$ let $\cont(u)$ denote
the signed distance of $u$ from the diagonal of $Y(\tau).$ 
We introduce a polynomial $F_\tau\in \Z[z]$ by the following
two equivalent formulas
\begin{equation}\label{F-pol}
F_\tau(z):=\;\prod\nolimits_{u\in Y(\tau)}\,
(1+\cont(u)z)\;=\; h(\tau)\cdot s_{\tau}(1,z, ..., z^{n-1}),
\end{equation}
where $h(\tau)$ is the product of hook-lengths of elements of $Y(\tau)$
and $s_{\tau}(p_1, p_2, ..., p_n)$ is the Schur polynomial of $\tau $
written in terms of power sums of symmetric variables.
The equality  of the two expressions in \eqref{F-pol} follows from
\cite[Example 4, Ch. I]{Ma}.

\begin{remark} The polynomial $F_\tau(z)$ has an interesting  interpretation
in terms of the Calogero-Moser space ${\mathcal{M}}_n$, cf. [Wi]. It is
known that fixed points of a
natural $\C^*$-action on  ${\mathcal{M}}_n$ are labelled by
Young diagrams, cf.
e.g. [EG, Prop. 4.16], [Ku]. Let $(X_\tau,Z_\tau)$ denote a pair of
$n\times n$-matrices in the conjugacy class corresponding
to the  fixed point in   ${\mathcal{M}}_n$ labelled by the Young diagram $Y(\tau)$
(for a direct construction of this point in terms of the
$S_n$-module $\tau$ see [Go], and [EG, Conjecture 17.14]).
Now, using formula \cite[(6.14)]{Wi} for Schur functions (and formula (2.3) in
Appendix to [Wi]), one finds
$$      F_{\tau}(z) = \det(-zX_\tau + (1-Z_\tau)^{-1})\,.
$$
\end{remark}
Let $\e_\tau\in\C{S_n} \subset \hh_c$ denote the central idempotent in
the group algebra corresponding to a simple $S_n$-module $\tau$.

\begin{theorem}\label{trace2} For any $c\in \C^*$ we have:
$\Tr(\e_\tau)=\frac{(\dim\tau)^2}{n!}\cdot F_{\tau}(\frac{1}{nc})\cd\Tr(1)$.
\end{theorem}

\begin{proof}
Let $\chi_\tau(g)$ denote the trace of $g\in S_n$ in the representation
$\tau$. Then one has an identity
\begin{equation}\label{orth}
\e_\tau=\frac{\dim\tau}{n!}\cdot\sum\nolimits_{g\in S_n}\;\chi_\tau(g^{-1})\cd g. 
\end{equation}
To prove this, observe that
 both sides are central elements of $\C[S_n], $
and computing their traces in an arbitrary module  $\sigma\in \irrep(S_n)$, 
we get 
$$
\frac{\dim\tau}{n!}\cdot
\sum\nolimits_{g\in S_n}\; \chi_\tau(g^{-1})\cd \chi_\sigma(g)= 1
\enspace\text{if}\enspace\tau=\sigma\enspace\text{and}\enspace0\enspace\text{otherwise}.
$$
The latter equation is nothing but the character orthogonality,
and \eqref{orth} follows.
Therefore, using Proposition \ref{trace1} we find
\begin{equation}\label{tr_eq2}
\Tr(\e_\tau)=\frac{\dim\tau}{n!}\cd
\sum\nolimits_{g\in S_n}\;\chi_\tau(g^{-1})\cdot(1/nc)^{n-\cyc(g)}
\cd\Tr(1).
\end{equation}
Now, using \cite[p.515, Exercise 7.50]{S2},
 we deduce:
$\dis
\sum_{g\in S_n}\chi_\tau(g^{-1})\cd z^{n-\cyc(g)}=\dim\tau\cd F_\tau(z)\,. 
$
\end{proof}

\begin{remark} If the parameter 
 $c\in\C^*$ is  generic, then it has been shown in [EG] that
 the vector space
$\hh_c/[\hh_c,\hh_c]$ is 1-dimensional.
Moreover, we claim that for  generic `$c$' we have  $\Tr(1)\neq 0$.
Indeed, assume the contrary, that is: $1\in [\hh_c,\hh_c]$.
This implies that there exists a positive integer $l=l(c)$
such that $1\in [F_l\hh_c,F_l\hh_c]$,
where $F_\bullet\hh_c$ is the standard increasing filtration on $\hh_c$.
Note that for any given $l$, the set of  `$c$' such that
  $1\in [F_l\hh_c,F_l\hh_c]$  is semialgebraic.
Hence, one can find such  an $l$ that  $1\in [F_l\hh_c,F_l\hh_c],$
for all  $c\in \C^*$ except possibly a finite set.
But we will see below that there are  infinitely many values of  `$c$' of the form
$c=1/n+ \text{positive integer}$ such that the algebra $\hh_c$
has a non-zero finite dimensional representation $V$.
In that case the element $1\in \hh_c$ acts as $\id_{_V}$,
hence has a nonzero trace.
Therefore, one cannot have  $1\in [\hh_c,\hh_c]$, and our claim is proved.

Thus, choosing,  for  generic `$c$', the element $\Tr(1)$ as a basis in
the 1-dimensional vector space $\hh_c/[\hh_c,\hh_c]$, one may
view the map $\Tr$ as a `trace': $\hh_c\to \C$.
\end{remark}\medskip

We are going to apply the results above to study finite dimensional 
representations of $\hh_c.$
Proposition \ref{trace1} clearly yields
\begin{corollary}\label{fin_dim} For any $g\in S_n$ and any  finite
dimensional 
$\hh_c$-module $V$, we have: $\;\Tr|_V(g)=(1/nc)^{n-\cyc(g)}\cd\dim V.$
\hfill\qed
\end{corollary}
\begin{remark} For $c=\infty,$ the Corollary implies that $V|_{S_n}$ is a multiple of the regular 
representation of $S_n,$ which has been already shown in [EG]. 
\end{remark}

Let $\mathbb{S}$ denote the set of all $c\in \C$ such that the algebra
 $\hh_c$ (of type ${\mathbf{A_{n-1}}}$) has a nonzero  finite dimensional 
representation. It is obvious from the above that $\mathbb{S}$ is a subset 
of $\Bbb Q^*$. 

\begin{proposition}\label{trace3} There is a subset
$S\subset \{\frac{1}{n},\frac{2}{n},...,\frac{n-1}{n}\}$,
such that  $\mathbb{S}= -(S+\Z_{\ge 0})\,\bigcup\,(S+\Z_{\ge 0}).$ 
\end{proposition}

\begin{proof} Throughout the proof we will freely exploit
an isomorphism $\hh_c\simeq \hh_{-c}$ 
given by $s\mapsto -s$ for each simple reflection $s$, and
$x\mapsto x$, $y\mapsto y$ for all $x\in \h^*,y\in \h$. In particular, we may assume $c>0$,
whenever necessary.

\noindent
{\sc Claim 1.}  {\it If $c>0$ is in $\mathbb{S}$ then either $c=p/n$
for some $p=1,...,n-1$, or else $c>1-\frac{1}{n}$. }

To prove this,  fix $c>0$, and a non-zero finite dimensional $\hh_c$-module $V$.
Using Theorem \ref{trace2} we obtain the following formula for the
multiplicity of occurence of
a simple $S_n$-module $\tau$ in $V$ 
\begin{equation}\label{multiplicity}
[\tau:V]=\frac{1}{n!}\sum \chi_\tau(g^{-1})\Tr|_V(g)=
\frac{\dim V\cd\dim\tau}{n!}\cdot\prod\nolimits_{u\in
Y(\tau)}\,(1+\frac{\cont(u)}{nc})\;.
\end{equation}
This means that the product $\prod_{u\in Y(\tau)}(1+\cont(u)/nc)$
has to be nonnegative for all $\tau.$ We deduce that either 
$c=p/n,$ where $1\le p\le n-1$ is a nonzero integer, or else
$c>1-\frac{1}{n}$. 

\noindent
{\sc Claim 2.}
{\it Let  $c\in {\mathbb{S}}$, and $c>1-\frac{1}{n}$. 
Then $c-1,c+1\in \mathbb{S}$. }

To prove this, we first
deduce from \eqref{multiplicity} that $[\triv:V]\neq 0$.
Thus $\e V$ is a nonzero
$\ehe$-module. 

Next,
let $\varepsilon: S_n \to \{\pm1\}$ be the sign character,
and $\e_{\varepsilon}\in \C[S_n]$ the corresponding central idempotent.
Using that  $c\ge 0$,
we deduce similarly from \eqref{multiplicity} that $[\varepsilon:V]\neq 0$.
It follows that  $\e_{\varepsilon}V$  is a nonzero
$\e_{\varepsilon}\hh_c\e_{\varepsilon}$-module.

Now, Proposition \ref{H_minus} implies that the  algebras $\ehe$ and
$\e_{\varepsilon}\hh_c\e_{\varepsilon}$
are  isomorphic 
to $\e_{\varepsilon}H_{c+1}\e_{\varepsilon}$ and 
$\e\hh_{c-1}\e$, respectively.
 Thus, the  latter algebras have   non-zero finite dimensional representations.
Hence, so do $\hh_{c+1}, \hh_{c-1}$ (by taking
induced module from the corresponding spherical
subalgebra). Claim 2 is proved. 

By a similar argument one proves

\noindent$\bullet$\quad {\it
 Let  $0<c\le 1-\frac{1}{n}$. If $\hh_c$  has a  non-zero finite dimensional representation
then so does} $\hh_{c+1}$. 

Now we complete the proof of the Proposition. Fix an element $c\in
\mathbb{S}$
such that
$c>0$.
If $c>1-1/n$ then $c\ge 1+1/n$, since in this case by Claim 2, $c-1$ should
be in $\mathbb{S}$, while we know from Claim 1 that $\mathbb{S}$
has no intersection with the open interval $(-1/n,1/n)$. So let
$c'$ be the fractional part of $c$. Then $c'$ is in $\mathbb{S}$,
by Claim 2. Let $S=\mathbb{S}\cap
\lbrace{1/n,2/n,...,1-1/n\rbrace}$. 
Then $c'\in S$, so $\mathbb{S}\subset -(S+\Bbb Z_{\ge 0})\cup
(S+\Bbb Z_{\ge 0})$. The fact that this inclusion is an
equality follows from the statement above that for $c>0$, if $c\in
\mathbb{S}$ then $c+1\in \mathbb{S}$. 
\end{proof}

\begin{remark} Note that $1/n \in \mathbb{S}$
 since for $c=1/n$ there exists an action of
$\hh_c $
on 
the trivial representation of $S_n,$ with $x_i$ and $y_i$ acting by zero. 
\end{remark}

Proposition  \ref{trace1} yields a nice closed expression for
 the character
of $S_n$-action on any  finite dimensional $\hh_c$-module, up to a constant factor.
To explain this, 
write $\h_{_{\Z}} := \Z^{n-1}\subset \Z^n$ for the set of $n$-tuples of integers with
total sum zero, a $\Z$-lattice in $\h$.
The group  $S_n$ acts naturally on  $\h_{_{\Z}}$ by permutation.
This induces, for any integer $r>0$, an  $S_n$-action on
the finite set $\h_{_{\Z}}/r\cd \h_{_{\Z}}\subset \Z^n/r\cd \Z^n$,
making the $\C$-vector space $\C[\h_{_{\Z}}/r\cd \h_{_{\Z}}]$
an $S_n$-module.

Let $K(S_n)$ be  the Grothendieck group of finite  dimensional representation
of 
$S_n$.
Given a  finite dimensional  $\hh_c$-module $V$,
let $[V|_{S_n}]\in K(S_n)$ denote the class of its restriction
to $S_n$.
Note  further that $V\neq 0$ implies, by Proposition \ref{trace3},
 that $cn=r$ is an integer, not divisible by $n$.

\begin{theorem}\label{traceCorollary} Assume $c=r/n>0,$ and let  $V$
be a  finite dimensional  $\hh_c$-module.

\vi
There exists an integer $\ell>0$
such that in  $K(S_n)$ we have: $\,
[V|_{S_n}] = \frac{\ell}{n}\cdot\text{\sf class of } \C[(\Z/r\Z)^n]$.

\vii Moreover, if $r$ is coprime to $n$
then $[V|_{S_n}]=\ell\cdot\text{\sf class of }\C[\h_{_{\Z}}/r\cd \h_{_{\Z}}]$. 
\end{theorem}

\begin{proof} To prove (i) it suffices to show that the characters
of $S_n$ representations on both sides are proportional.
But it is straightforward to see that the character of
$S_n$-action on $\C[(\Z/r\Z)^n]$ is given,
by $\Tr(g)= r^{\cyc(g)}= (cn)^{\cyc(g)}$,
which up to a scalar is exactly  the formula of Proposition
\ref{trace1}.
The integrality of the coefficient $\ell$ in part (i)
follows from the fact that $\ell$ is equal to the
trace of the cyclic permutation in $(\ell/n)\cdot\C[(\Z/r\Z)^n]$.
Part (ii) readily follows from (i) by comparing
characters. 
\end{proof}

\begin{conjecture} Let $c=r/n>0$. The algebra  $\hh_c$
has a nonzero  finite dimensional  representation $V$ if and only if
  $r$ is coprime
to $n$, moreover, in such a case
 in Theorem \ref{traceCorollary}(ii) we have $\ell=1$.
\end{conjecture}

This conjecture has been proved in \cite{BEG1}.

\section{The $\e\hh_c\e$-module structure 
on quasi-invariants}
\setcounter{equation}{0}
In this section we review the definition and prove some basic 
properties of algebras of quasi-invariants of finite reflection 
groups. 
Most of these results have already appeared in the literature
(see \cite{CV}, \cite{CV2}, \cite{VSC}, \cite{FV}, \cite{EG2}). 

Let $ \Hol(U)$ be the space of  holomorphic functions
on an open set $U\subset \h$, and $\Hol(\h)^W \subset  \Hol(\h)$, 
the space of $W$-invariant holomorphic functions on $\h$.
Observe that $f\in \Hol(\h)$ is  $W$-invariant if and only if 
$\, s_{\alpha}(f) = f \,$ for all $\, \alpha \in R_{+},\,$
since the set $ \{\,s_{\alpha}\,|\, \alpha \in R_{+}\,\} 
\subset W \,$ generates $ W $.
We will extend the algebra $ \Hol(\h)^W  $ by relaxing the notion
of reflection invariance in the following way.
Let $\, m \in \Z_{+} \,$ be a non-negative integer,
$\, s_\alpha \in W\,$ a reflection,
and $x\in\h$ a point on the corresponding hyperplane,
i.e., $s_\alpha(x)=x.$ We say that
a holomorphic function $f$ on an $s_\alpha$-stable
open neighborhood $U_x\ni x$
is $m$-{\it quasi-invariant} 
with respect to $\, s_\alpha \,$ if
\begin{equation}
\label{1}
s_{\alpha}(f) \equiv f \ \mbox{mod}\, (\alpha)^{2m + 1}\ ,
\end{equation}
where $\, (\alpha) $ stands for the (principal) ideal
in $ \Hol(U_x) $ generated by the linear form $ \alpha\,$.
To compare, a function $\, f \in  \Hol(U_x) $ being invariant 
under $ s_\alpha $ means that all {\it odd} powers of $ \alpha $
in the Taylor expansion of $ f $ in the
direction $ \alpha $ vanish, 
while  $\, f \in  \Hol(U_x) $ being $m$-{\it quasi-}invariant 
means that such powers vanish only up to 
order $\, 2m-1 \,$ (inclusive).

Let $ \Z(R)^{W}_{+} $ be the set of all $ W$-invariant functions
on the root system $ R $ with values in non-negative integers.
\begin{definition}[cf. \cite{CV}, \cite{CV2}]
\label{qin}
Given $ c \in \Z(R)^{W}_{+}, $ a holomorphic function $ f \in \Hol(\h)$
is said to be $c$-{\it quasi-invariant} with respect to $ W $
if,  for each $\, \alpha \in R,$ one has 
$s_{\alpha}(f) \equiv f \ \mbox{mod}\, (\alpha)^{2c_\alpha + 1}.$
\end{definition}
Given $ c \in  \Z(R)^{W}_{+},$ let
 $ Q_c$ denote the set of all $c$-quasi-invariant 
{\it polynomials}
 with respect to $ W $.
The following lemma collects some basic algebraic properties of 
quasi-invariants.
\begin{lemma}
\label{L1} \vi The vector space $Q_c$ is $W$-stable, and
we have: $ \C[\h]^{W} \subseteq Q_c \subseteq \C[\h]\,$.

\vii $ Q_c $ is a  finitely generated graded 
subalgebra of $ \C[\h] $, such that
the integral closure of $ Q_c $ equals $ \C[\h]\,$.
\end{lemma}
\begin{proof} (i) is clear.
Further, if $ p, q \in \C[\h] $ satisfy (\ref{1})
with $m=c_\alpha$, 
for some $ \alpha \in R_{+} ,$ so do their linear combination and 
product: for example, $\, s_{\alpha}(pq) = 
s_{\alpha}(p) s_{\alpha}(q) \in pq  + (\alpha)^{2c_{\alpha}+1}\,$.
Now, if (\ref{1}) holds for all $ \alpha \in R_{+} \,$
and $ \beta \in R_{+} $ then 
$\, s_{\alpha}[s_{_\beta}(p)] = s_{\beta}[ s_{s_{_\beta}(\alpha)}(p)]
\, \in \, s_{\beta}(p) + (\alpha)^{2c_{s_{_\beta}(\alpha)}+1} =
s_{\beta}(p) + (\alpha)^{2c_{\alpha}+1}
\,$, where the last equality holds due to the $W$-invariance of $ c\,$.
Thus $ Q_c $ is a subalgebra of $ \C[\h]\,$ invariant under the action 
$ W\,$. 
Since $ W $ is a finite group, $\, \C[\h] $ is finite as an
$ \C[\h]^{W}$-module, see e.g. \cite{S}, hence,
as a $Q_c$-module. Then finiteness of $ Q_c $ 
as a $\C$-algebra follows from the Hilbert-Noether Lemma 
(see \cite{AM}, Proposition~7.8).
Further, $\, \C[\h] $ being finite as a module over $ Q_c $
is equivalent to $\, \C[\h] $ being integral as a ring over
$ Q_c \,$. On the other hand, the subalgebra 
$ Q_c $ contains an ideal of $ \C[\h]\,$, namely a 
sufficiently high power of 
$\,\delta=\prod_{\alpha \in R}\, \alpha\,$. Hence, $ Q_c $
and $ \C[\h] $ have the same field of fractions. 
It follows that $ \C[\h] $ is the integral closure of 
$ Q_c\,$.
\end{proof}

We recall the following elementary but basic observation [C],[FV]:
 \begin{lemma}\label{L4}
The action of the Calogero-Moser operator
$\LL_c$ preserves the subspace $Q_c\subset
\C[\hreg]$, i.e., $\LL_c(Q_c)\subseteq Q_c$.\qed
\end{lemma}

Recall that according to Opdam [O],
for any $c\in\CC$, there is an algebra isomorphism
$\C[\h^*]^W\iso \cc_c,$ where $\cc_c$ stands for the centralizer 
in $\dd(\hreg)_\minus^W$ of the Calogero-Moser operator $\LL_c$.
Identify $\C[\h^*]^W$ with $\C[\h]^W$ via an invariant form on $\h$.
In [CV], 
Chalykh-Veselov  extended Opdam's isomorphism:
$\C[\h]^W = \C[\h^*]^W\iso \cc_c$
to an injective algebra homomorphism: $Q_c 
\into \dd(\hreg)_\minus
\,,\,
P\mapsto \ll_{_P}.$ There is an
 explicit formula  for the restriction of the map
$P\mapsto \ll_{_P}$  to the subspace $Q_c\cap \C^k[\h^*]$ of
homogeneous polynomials in $Q_c$ of degree $k\geq 1$.
 Specifically, according to [Be],
cf. also \S7 below, one has
\begin{equation}
\label{P3}
\ll_{_P}=\frac{1}{2^{k}\,k!}\,
(\ad\LL_c)^{k}(P),
\end{equation}
where $P$ is a homogeneous polynomial in $Q_c$ of degree $k\ge 1$.
In the special case when $P\in \C^k[\h^*]^W$ is $W$-invariant,
this formula is an immediate (up to computing  the 
constant factor, $\frac{1}{2^{k}\,k!}$, which is also easy)
consequence of Corollary \ref{Lefschetz}, which 
yields the isomorphism 
$\,\Theta\spher(\ad\YY)^k: \Theta\spher(\C^k[\h]^W\!\!\cd\e) 
\iso\Theta\spher(\C^k[\h^*]^W\!\!\cd\e)$. 

Let $Q_c\ff\subset \dd(\hreg)_\minus$
denote the image of the Chalykh-Veselov homomorphism
$P\mapsto \ll_{_P}.$ Thus $Q_c\ff$ is a commutative 
subalgebra in $\dd(\hreg)$ isomorphic to $Q_c$.
As has been observed  in [FV], formula (\ref{P3})
 and
Lemma \ref{L4} imply that the
action on $\C[\hreg]$ of the algebra $Q_c^\flat$ preserves 
the space $Q_c$, that is, $u(Q_c)\subset Q_c,$ for any $u\in
Q_c^\flat$. In particular, 
$\cc_c(Q_c)\subset Q_c.$ Thus, there is a natural
action on $Q_c$ of the algebra
${\mathcal{B}}_c\subset \dd(\hreg)^W$, 
generated by $\cc_c$ and $\C[\h]^W=Q_c^W$.
This $\bb_c$-action clearly
commutes with the $W$-action on $Q_c$.
Thus, since $\bb_c=\Theta\spher(\ehe)$ by Proposition \ref{spher},
 the space $Q_c$ acquires
an $\C W \otimes \e\hh_c\e$-module structure.
The structure of this module is completely
described by the following proposition. 

For any irreducible representation $\tau$ of $W$, let $\tau_c'$ 
be the representation of $W$ for which 
the monodromy representation of the Dunkl connection with 
values in $\tau_c'$ is $\tau$ (see \cite{O1} and references
therein).\footnote{A priori, the monodromy representation 
is a representation of the
braid group of $W$, but since $c$ is integer-valued, it factors
through $W$. It is known (\cite{O1}) that $\tau_c'=\tau$ 
for all irreducible Coxeter groups
except $E_7,E_8,H_3,H_4$, and even dihedral groups; otherwise 
$\tau_c'=\gamma^{c}(\tau)$, where $\gamma$ is a certain
involution on the set of irreducible representations of $W$, and 
for even dihedral groups $\gamma^c$ stands for
$\gamma^{c_1+c_2}$. Also, it is known that characters of $\tau$
and $\tau_c'$ coincide on reflections.}

\begin{proposition}\label{isot}
There is a $\C W \otimes \e\hh_c\e$-module isomorphism:
$Q_c=\bigoplus\limits_{\tau\in \irrep(W)}\tau^*\otimes M_\e(\tau_c').$
\end{proposition}

\begin{proof} First of all,  $Q_c$ is  finitely generated
over $\C[\h]^W$, hence over $\e\hh_c\e$. Further, the action
of any homogeneous element of $Q_c^\flat$ lowers the degree
of polynomials in $Q_c$. Hence, the action on  $Q_c$
of the augmentation ideal
of the algebra $\C[\h^*]^W \subset \ehe$ is locally
nilpotent. It follows that $Q_c\in \oo_{\ehe}(0)$.
Hence, Lemma \ref{O_prop}(iv) implies that $Q_c$
has finite length, as a $\ehe$-module.
Recall further that, for any $c\in\ccp$,
 the standard modules $M(\tau)$  are 
irreducible. By Morita equivalence of $\hh_c$ and $\e\hh_c\e$, 
for any $\tau$, the modules 
$M_\e(\tau):=\e M(\tau)$ over $\e\hh_c\e$
are also irreducible. Thus, $Q_c$ has a
finite composition series whose simple subquotients
are standard modules  over $\e\hh_c\e$.

Observe next that
the direct sum decomposition of the Proposition
is easily checked directly for $c=0$. 
Now, consider the localized algebra $(\e\hh_c\e)\mrw$
and the  module 
 $Q_c\mrw$ over this algebra. 
It is clear from the existence of the Dunkl representation that 
both the localized algebra 
and module are independent of $c$, i.e. are the same as those
for $c=0$ (namely, they are the algebra of differential operators
on ${\mathfrak h}_{\rm
reg}/W$ and the D-module on ${\mathfrak h}_{\rm
reg}/W$ corresponding to the regular representation of $W$). 
Therefore, we have $Q_c\mrw=\bigoplus_{\tau\in \irrep(W)}\,\tau^*\,\otimes\, 
M_\e(\tau_c')\mrw$ (by the definition of $\tau_c'$).

Now, define $M_\e'(\tau):=(\tau\otimes Q_c)^W$. This is an 
$\e\hh_c\e$-module, such that the localized module 
$M_\e'(\tau)\mrw$ coincides with $M_\e(\tau_c')\mrw$.  
On the other hand,  
$M_\e'(\tau)\in \oo_{\e\hh_c\e}(0)$, so it has a finite
Jordan-H\"older series $0=F_0\subset F_1\subset \ldots\subset 
F_l=M_\e'(\tau),$ with simple subquotients $F_i/F_{i-1}$
being isomorphic to
standard modules $M_\e(\sigma_i)\,,\, i=1,\ldots,l,$ 
see Proposition \ref{O=O}. The localization functor being exact, 
we see that $M_\e'(\tau)\mrw$ has a
composition series with the subquotients $(F_i/F_{i-1})\mrw=
M_\e(\sigma_i)\mrw.$ But since  
$M_\e'(\tau)\mrw=M_\e(\tau_c')\mrw$,
the latter composition series must consist of a single term,
i.e., $l=1$. Thus,
we find that $M_\e'(\tau)=M_\e(\tau_c')$, and we are done. 
\end{proof}

Given a diagonalisable endomorphism $A$ on a vector
space $E$ with finite-dimensional eigenspaces $E_a\,,\,a\in\C,$
we put $\Tr\big|_E(t^A):= \sum_{a\in\C}\, \dim_{_\C} E_a\cdot
t^a,$ viewed as a formal infinite sum.
For example,  let $\eu$ denote the Euler operator
on $\C[\h]$ such that $\eu(P)=(\deg P)\cdot P,$
for any homogeneous polynomial $P$. Write 
$\chi_\tau(t)=\Tr\big|_{(\tau\otimes \C[\h])^W}(t^\eu)$ 
for the Poincar\'e series of the graded space 
$(\tau\otimes \C[\h])^W$. These series are classically known.

Let $\kappa_c=\frac{1}{2}\sum_{\alpha\in R}c_\alpha(1-s_\alpha)$ 
be the central element of $\C W$ canonically attached
to $c\in\CC$. This element acts by a scalar, say 
$\kappa_c(\tau)$, in each irreducible representation $\tau\in\irrep(W)$.
We have
\begin{equation}\label{verma_char}
\Tr|_{M_\e(\tau)}(t^\sh)=t^{\frac{\dim\h}{2}-\frac{1}{2}
\sum\nolimits_{\alpha\in R}c_\alpha
+\kappa_c(\tau)}\cdot \chi_\tau(t).
\end{equation}
To prove this formula, one first  finds
 the eigenvalue of $\sh$-action on the "top" subspace
$\tau\subset M(\tau)$. A straightforward calculation
gives
$$ 
\sh|_\tau= \left(
\frac{\dim\h}{2}-\frac{1}{2}
\sum\nolimits_{\alpha\in R}\,c_\alpha+\kappa_c(\tau)\right)\cdot\id_\tau\,.
$$ 
Further, it is easy to see from the Dunkl representation 
that $\sh|_{\C[\h]}=\eu+\frac{\dim\h}{2}-\frac{1}{2}
\sum_{\alpha\in R}c_\alpha$. Formula (\ref{verma_char})
follows.\qed

Applying Proposition \ref{isot}, we obtain
a formula for the Poincar\'e series  
of $Q_c$, i.e. for the function $\Tr|_{Q_c}(t^{\eu})$.
The same formula 
has been found by  Felder and Veselov (\cite{FeV}),
using a different method.

\begin{corollary} 
$\quad
\Tr|_{Q_c}(t^{\eu})=\sum_{\tau\in \irrep(W)}
\dim\tau\cdot t^{\kappa_c(\tau)}\cdot \chi_\tau(t)
.$\qed
\end{corollary}

\begin{proof}
This follows from Proposition \ref{isot} and the equality that
$\kappa_c(\tau)=\kappa_c(\tau_c')$, which follows from the fact
that the characters of $\tau$ and $\tau_c'$ coincide on
reflections.  
\end{proof}

\section{Differential operators on quasi-invariants}
\setcounter{equation}{0}

We recall some basic results and notation regarding
differential operators on algebraic varieties.
Let $\, A $ be a commutative $ \C$-algebra, and let
$\, M $ and $ N $ be a pair of $A$-modules.  The 
space $ \dd(M,N) $ of linear differential operators 
from $ M $ to $ N $ over $ A $ is a filtered subspace 
$\dd(M,N)=\cup_{n\geq n}\,\dd_{n}(M,N)\subset \mbox{Hom}_{\C}(M,N) $ 
defined inductively by 
$\, \dd_{-1}(M,N) = 0 \,$, and 
\begin{equation}
\label{1.1}
\dd_{n}(M,N) := 
\{u \in \mbox{Hom}_{\C}(M,N)\mid u\ccirc a - a\ccirc u \in 
\dd_{n-1}(M,N)\ \mbox{for all}\ a \in A\}\;,\;n> 0 .
\end{equation}
The elements of $ \dd_{n}(M,N) \setminus  \dd_{n-1}(M,N) $
are called {\it differential operators of order} $ n \,$.

If $ N = M=A \,$, we write $\, \dd(A) := \dd(A,A), \,$
which is clearly an algebra over $ \C\,$. In particular, 
if $ A = \C[X] $ is the ring of regular functions on
an irreducible  affine algebraic variety $ X \,$,
then $\, \dd(A) $ is nothing but  $\, \dd(X), \,$
the  algebra of differential 
operators on $ X.$

If $ X $ is singular, the structure of $ \dd(X) $ 
can be quite complicated and is not well understood.
A natural approach to study the algebra $ \dd(X) $ on 
an arbitrary (irreducible) variety $ X $ is to compare it 
with the ring of differential operators on the normalization 
of $ X\,$ (see \cite{SS}). In more detail, let
 $ \pi: \tilde{X} \to X $ be the 
normalization map. Write $ \BK := \C(X) $ for the field of 
rational functions on $ X \,$
(so that $\, \C[X] \subseteq \C[\tilde{X}] \subset \BK \,$, and
$ \C[\tilde{X}] $ is the integral closure of $ \C[X] $ in $ \BK \,$) 
and denote by $ \dd(\BK) $ the algebra of differential
operators over $ \BK \,$. Since $ \dd(\BK) \cong \BK \otimes_{\C[X]} \dd(X) 
\,$, we can identify $ \dd(X) $ with its image in $ \dd(\BK).$
This way we get
\begin{equation}\label{3}
\dd(X) = \{\, u \in \dd(\BK)\ \mid u\bigl(\C[X]\bigr) \subseteq \C[X]\,\}
\quad,\quad 
\dd(\tilde{X}) = \{\, u \in \dd(\BK)\ \mid u\bigl(\C[\tilde{X}]\bigr)
\subseteq \C[\tilde{X}]\,\} \ .
\end{equation}
Thus, both $ \dd(X) $ and $ \dd(\tilde{X}) $ can be viewed
as subalgebras in $ \dd(\BK)\,$.

    From now on, we fix $c\in\ccp$ and 
consider
the scheme $ X_c := \Spec\,Q_c \,$.
Given $x\in\h$, let $W_x\subset W$ 
denote the stabilizer of $x$ under the natural action of $W$ on $\h$.
\begin{lemma}
\label{C1}\vi
$\; X_c $ is a rational
irreducible  affine algebraic variety with $W$-action.
The normalization  of $ X_c $
is isomorphic to $ \h\,$.

\noindent
\vii The normalization map $\pi_c: \h\onto X_c$
is bijective.

\noindent
\viii The local ring $\co_x(X_c)$ at a point $x \in X_c$ can be identified with
the subring of $\C(\h)$ consisting of rational functions on $\h$ regular at 
$\pi_c^{-1}(x)$ and quasi-invariant under $W_x$.
\end{lemma}

\begin{proof} Part (i) is a reformulation  of Lemma~\ref{L1}.
To prove (ii) we only need to check injectivity
of the  normalization map. Thus, we have
to show that for any points
$x,x'\in \h\,,\, x'\neq x,$ there exists a 
quasi-invariant polynomial $p\in Q_c$ such that $p(x)\neq
p(x').$ To this end, fix $x\in \h$ and put
\begin{equation}\label{delta}
\delta_{2c+1,x}=\!\prod_{\{\alpha\in R_+\; \mid\; (\alpha,x)\ne 0\}}\!
\alpha^{2c_\alpha+1},\;\;\begin{array}{l}
{\text{\footnotesize this is
a quasi-invariant polynomial under all}}\\ 
{\text{\footnotesize  reflections in}}\enspace W\smallsetminus
W_x,\enspace{\text{\footnotesize  and also invariant under $W_x$\,.}}
\end{array}
\end{equation}
Now, given $x'\neq x,$
choose any polynomial $g$
such that $g(x)\neq 0$, and $g(x')=0.$
Define the polynomial
$p(z):=\delta_{2c+1,x}(z)\cdot\!\prod\nolimits_{w\in W_x}g(wz).$
We have:
$p(x)\neq 0$, and  $p(x')=0$. Further,
$p$ is invariant under $W_x$ and quasi-invariant under all reflections in
$W\ssminus W_x,$ by construction. Thus, $p$ is quasi-invariant and
(ii) is proved.

To prove (iii) observe that $\co_x(X_c) \subset \{r\in\C(\h)\mid
r\;\;\text{\tt is regular at}\;\;x\}.$
Conversely, let $r$ be a rational function  regular at $x$
which is also quasi-invariant with respect to the group $W_x$.
We can write $r=p/q$, where $p, q$ are polynomials, 
and $q(x)\neq 0$. To show that  $r=p/q\in \co_x(X_c)$
we  set $f=\prod_{\{w\in W_x\mid w\ne 1\}}\, w(q),$
and write $r=\frac{p}{q}=\frac{pf\delta_{2c+1,x}}{qf\delta_{2c+1,x}}$,
where
 $\delta_{2c+1,x}$ was defined in (\ref{delta}). The denominator
in the last fraction is $W_x$-invariant, and also
quasi-invariant under all the
reflections which are not in $W_x$ by construction.
Thus, the denominator is quasi-invariant under $W$ and
nonvanishing at $x.$ The numerator is also quasi-invariant, since so is the
whole function $r$. Thus, $r\in \co_x(X_c),$ as desired.
\end{proof}

According to  (\ref{3}) we may identify
the algebra $\dd(X_c)$  with the set
$\{u\in\dd(\hreg)\mid
u(Q_c)\subset Q_c\}$. In the previous section we have
introduced the subalgebra
$\,\dd(\hreg)_\minus=\text{\sf span of}$ $\{u\in \dd(\hreg)\mid
\text{\it weight}(u)+\text{\it order}(u)\leq 0\}.\,$
Put  $\ddx_\minus:= \ddx\cap\dd(\hreg)_\minus$.
Recall further the subalgebra
$Q_c\ff\subset \dd(\hreg)_\minus,$ the image of $Q_c$ under
the
Chalykh-Veselov homomorphism $P\mapsto \ll_{_P}$.

\begin{theorem}\label{DDX}
\vi We have: $\ddx_\minus=Q_c^\flat$. This is a maximal
commutative subalgebra in $\ddx$,  equal  to
  the centralizer in 
$\dd(\hreg)_\minus$
of the Calogero-Moser
operator $\LL_c$.

\vii The left multiplication by $Q_c$ and  right 
multiplication by $Q_c\ff$ make $\ddx$ into a rank one
{\sl projective} $Q_c\otimes Q_c\ff$-module.
\end{theorem}

The rest of the section is effectively devoted to
the proof of this Theorem. The argument exploits crucially
the
notion of  the 
Baker-Akhiezer (holomorphic)
function $\psi_c: (x,k)\mapsto \psi_c(x,k) = \Phi(x,k)\cd e^{(x,k)},$
 where $x\in\h,\,k\in\h^*,$ and $\Phi$ is a polynomial on $\h\times\h^*$,
see \cite{CV,FV,EG2}.
\smallskip

Whenever necessary we will identify below  $\h^*$ with $\h$
without further comment,
using some fixed $W$-invariant form. 
Given  $x\in\h$, put  
$$Q_{c,x}:=
\{f\in \C[\h]\;\mid\; k\mapsto f(k)e^{kx}
\enspace \text{ is a}\enspace \text{quasi-invariant holomorphic
function}\}\,.
$$
\begin{proposition}\label{kth1} For any   $x\in\h$, 
the space $Q_{c,x}$ is a projective rank
one $Q_c$-module (not free, in general).
\end{proposition}
\begin{proof}
Fix $x_0\in\h$. 
We will show that the module $Q_{c,x_0}$
 is locally-free in a neighborhood of any point
$k_0\in\h.$ 
Let $ g(t,k) := \psi_{c}(t\cd x_0, k),$ where $t\in\C$ is a complex variable.
Then, for each $ k \in \h,$
the map: $ t \mapsto g(t,k)$  is an entire function on $ \C$ . Write
$ g(t,k) = \sum_{n=0}^{\infty} g_{n}(k)\cd\frac{(t-1)^n}{n!}  .$  For
all $ t \in \C,$
we know that 
$ \psi_{c}(t\cd x_0, k)$  is a quasi-invariant holomorphic function in $ k.$ 
Hence,
the Taylor coefficients $ g_{n}(k) = (d/dt)^n g(t,k)|_{t=1}$ are
quasi-invariant holomorphic functions in $k$. Thus, 
$ g_{n}(k) \cd e^{-(x_0 k)} \in Q_{c,x_0},$  for each $ n = 0,1,2,...$ 
Further, it is known that $\psi_c(0,0)\neq 0,$ see e.g. [EG2].
We deduce that $ g_{n_0}(k_0) \not= 0$  for some $ n = n_0,$  for otherwise $ g(t,
k_0) = 0$ 
would imply  $ g(0, k_0) = \psi_{c}(0, k_0) = \psi_{c}(0,0) = 0$.
Set $ q_0(k) := g_{n_0}(k)\cdot e^{-(x_0 k)}.$  Then for
any 
$ p \in Q_{c,x_0},$  we have $ p(k) = q_0(k)\cd f(k),$  where $ f(k)\in \C(\h)$  is a rational 
function  regular at $ k = k_0.$  Observe that $ f$  can be written as
$ \frac{p(k)\cdot e^{(x_0 k)}}{g_{n_0}(k)}$,   a
ratio of two quasi-invariant functions in $ k.$  Hence, 
$ f(k) \in \co_{k_0}(X_c)$  by Lemma \ref{C1}(iii). 
It follows that $\co_{k_0}(X_c)\otimes_{\co(X_c)}Q_{c,x_0}$ is a 
free rank 1 module generated by~$q_0.$~\end{proof}

\begin{definition}\label{kk}
 Let $\kk\subset \C[\h\times \h^*]$ be the space of all
polynomials $P$ such that the holomorphic
function: $(x,k)\mapsto P(x,k)e^{(x,k)}$
is 
quasi-invariant with respect to both the variable
  $x\in \h$ and  the variable
$k\in\h^*,$ separately. 
\end{definition}

Observe that
the space $\kk$ has a natural $Q_c\otimes Q_c$-module
structure induced by multiplication in $\C[\h]\otimes\C[\h^*]$.

\begin{proposition}\label{kth2}
The space $\kk$ is a  projective rank
one $Q_c\otimes Q_c$-module.
\end{proposition}

\begin{proof} The proof is analogous to the argument above.
 Specifically, we choose as in Proposition
\ref{kth1}, a function
$g_0\in \kk$ that is nonvanishing at a given point $(x_0,k_0).$
Then, for any $g\in\kk$, the function
 $\frac{g}{g_0}$ is a {\it rational} function (not involving the
exponents) on
$\h\times\h^*$. Moreover, this function is regular
at the point $(x_0,k_0)\in \h\times\h^*$ and $W_{x_0}\times 
W_{k_0}$-quasi-invariant.
Hence,
 by 
Lemma \ref{C1}(iii) applied to
the Coxeter group $W\times W$ in the vector space  $\h\times\h^*$,
we obtain $\frac{g}{g_0}\in \co_{x_0, k_0}(X_c\times X_c).$
It follows that the function $g_0$ generates the $\co_{_{x_0,
k_0}}(X_c\times X_c)$-module
$\co_{_{\!(x_0, k_0)}}(X_c\times X_c)\bigotimes_{_{\!Q_c\otimes Q_c}}\kk$.
Thus,  $\kk$ is a locally free $Q_c\otimes Q_c$-module of
rank one.
\end{proof}
\begin{corollary}\label{kth3} 
For any $a\in X_c$, we have: $\kk|_{X_c\times \{a\}} = Q_{c,a}.$
\end{corollary}

\begin{proof}
We have an obvious evaluation map ${\mathtt{ev}}: 
\kk|_{X_c\times \{a\}}\to Q_{c,a}.$
Both sides are rank one locally free sheaves
on $X_c\times \{a\}$.
 Hence, to prove bijectivity we only
need to check that the map is surjective.
This reduces to checking that for any $x$ there exists 
$f\in \kk|_{X_c\times \{a\}}$
such that ${\mathtt{ev}}(f)(x)$ is nonzero.
But ${\mathtt{ev}}(f)(x)=f(x,a)$, so one can take $f=g_0$,
the same function as
has been taken in the proof of Proposition \ref{kth2}.
\end{proof}

In terms of the $\psi_c$-function, the Chalykh-Veselov map
$Q_c \to Q_c\ff\,,\,P\mapsto\ll_{_P}, $ is defined as follows.
Given $P\in Q_c$, let
$\ll_{_P}\in\dd(\hreg)$ be the unique differential
operator such that $(\ll_{_P}^{(x)}\psi_c)(x,k)=
P(k)\cdot\psi_c(x,k)$, where 
$u^{(x)}$ stands for the operator $u$ acting in the
"$x$-variable". 

Below, we will use the  following result due to
Chalykh-Veselov (we use the notation 
 $\sigma^{(k)}(p)$ for the leading term of $p\in\C[\hr\times\h^*]$,
viewed as a polynomial in the variable `$k$' only,
and write $\sigma(u)$  for the principal symbol of a differential
operator $u$).

\begin{lemma}[\cite{CV2}, p.33]\label{cv2}
Let $p\in \C[\hr\times\h^*]$ be a function
in the variables $(x,k)\in \hr\times\h^*$
such that
the function $F: (x,k) \mapsto p(x,k)\cdot e^{(x,k)}$
is quasi-invariant with respect to the variable
`$k$'. Then

\vi There exists a
differential operator $u\in \dd(\hr)$
such that $F=u^{(x)}\psi_c$.

\vii We have:
$\displaystyle\sigma^{(k)}(p)\in\C[\hr\times\h^*]\cdot\delta_c(k)$ 
and, moreover, 
$\sigma(u)=\frac{\sigma^{(k)}(p)}{\delta_c(x)\delta_c(k)}$.
\qed
\end{lemma}

Recall next that in
[EG2] we have shown that formula 
$\,\boldsymbol{(} f,g\boldsymbol{)}_c
:= (\ll_f g)(0)\,$
gives a non-degenerate symmetric $\C$-bilinear
pairing on the vector space $Q_c$.

 Further, given $P\in \kk$ we write $\deg\fx P$ for the {\it degree} of
$P$, viewed as a polynomial in the $x$-variable only.
Also,
 we write ${\it weight}(u)$ for the weight
of  a $\C^*$-homogeneous differential operator
 $u\in \dd(X_c)$  with respect to the $\C^*$-action and ${\it order}(u)$
for the order of $u$ as a differential operator.
We let the ring $Q_c\otimes Q_c\ff$ act on
$\kk$ via the identification
$Q_c\otimes Q_c\simeq Q_c\otimes Q_c\ff$ induced by
the Chalykh-Veselov isomorphism $Q_c\iso Q_c\ff$ on the second
tensor factors.

We define a linear map $K: \dd(X_c)\to \kk$, 
given by $u\mapsto
K(u)=e^{-(x,k)}\cdot\bigl(u^{(x)}\psi_c(x,k)\bigr)$.

\begin{proposition}\label{kernel1}\vi
The map $K: \dd(X_c)\to \kk$ is a $Q_c\otimes Q_c\ff$-module
isomorphism.

\vii For any $u\in \dd(X_c)$, we have:
$\deg\fx\bigl(K(u)\bigr)-\deg\fx(\Phi)= \text{\it weight}(u) + 
\text{\it  order}(u)$, moreover, this integer is $\geq 0$.
\end{proposition}

\begin{proof} First of all, $K$ is obviously injective. 
So, we need to show that for any $P\in \kk$, there exists 
$u\in \dd(X_c)$ such that $P=K(u)$. 
Since $P(x,k)e^{(x,k)}$ is quasi-invariant with respect to $k$, 
it follows from Lemma \ref{cv2} that there exists a
differential operator $u\in\dd(\hreg)$  such that 
$(u^{(x)}\psi_c)(x,k)$
$=P(x,k)e^{(x,k)}$. 

 Further,
by \cite{EG2}, one has an infinite expansion:
 $\psi_c(x,k)=\sum \psi_i(k)\psi^i(x)$, 
 where $\psi_i,\psi^i\in Q_c$ are homogeneous polynomials
which form
dual bases with respect to the nondegenerate form 
$\boldsymbol{(}-,-\boldsymbol{)}_c$. 
   From the expansion we get:
$P(x,k)e^{(x,k)}=(u^{(x)}\psi_c)(x,k)
=\sum_i\,\psi_i(k)\cdot(u^{(x)}\psi^i(x))$.
Thus,
since $P(x,k)e^{(x,k)}$ is
also quasi-invariant with respect to $x$, 
we find that $u\psi^i\in Q_c\,,\,\forall i$.
Hence the operator 
$u$ maps $Q_c$ to $Q_c$. So, $u\in\ddx$ and $P=K(u)$, as desired.
Part (i) is proved.

Now, given $P\in Q_c$ and  $u\in\ddx$,
we find:
$$
\begin{array}{ll}
K\bigl(u\ccirc \ll_{_P}\bigr)(x,k)\cdot e^{(x,k)}&=
(u^{(x)}\ccirc \ll_{_P}^{(x)}\psi_c)(x,k)=
u^{(x)}\bigl(P(k)\cdot  \psi_c(x,k)\bigr)
\\
&=
P(k)\cdot  \bigl(u^{(x)}\psi_c(x,k)\bigr)=
P(k)\cd K(u)(x,k)\cdot  e^{(x,k)}\,.
\end{array}
$$
Compatibility with
$Q_c\otimes Q_c\ff$-module structures follows.

Observe further that if $ u$ is  homogeneous under the $\C^*$-action
(of weight $m$, say), then $K(u)$ is homogeneous of the
same weight. It follows that the highest powers in
$x$
always come
in $K(u)$ together with the
highest powers in $k,$ i.e. come from the leading term in both $k$ and $x$
simultaneously.
Now, the leading term of $K(u)$ is equal to
the principal  symbol of $u$ times the leading
term of $\Phi.$ This implies (ii).
\end{proof}

We introduce a {\it non-standard} increasing
filtration on $\ddx$ by 
\begin{equation}\label{filt_flat}
\dd_j^\flat(X_c) := 
\text{\sf span of}_{\,}\{u\in \ddx\enspace\big|\quad
\text{\it weight}(u)+\text{\it order}(u)\leq j\}\enspace,\,j\in\Z.
\end{equation}
For any
homogeneous differential operators $u$ and $u',$ we have
$\text{\it weight}(u\ccirc u')=\text{\it weight}(u)+
\text{\it weight}(u'),$ and
$\text{\it order}(u\ccirc u')
\leq \text{\it order}(u)+\text{\it order}(u')$.
Hence, for any $j,l\in \Z$, we get:
$\dd_j^\flat(X_c)\cdot \dd_l^\flat(X_c)\subseteq
\dd_{j+l}^\flat(X_c)$. 
Moreover, since $\text{\it order}(u\ccirc u'-
u'\ccirc u) <\text{\it order}(u)+\text{\it order}(u')$,
we deduce that the associated graded algebra
$\grd\ff\ddx$ is commutative.
Further,
by Proposition \ref{kernel1}(ii), we have:
$\dd_{-1}^\flat(X_c)=0,$ so that our filtration
is by {\it non-negative} integers $j\geq 0$.
 Thus,
we have proved the inclusion on the right of (\ref{filt_2}):
\begin{equation}\label{filt_2}
\ddx_\minus=\dd_0^\flat(X_c) = Q_c^\flat\,,
\quad\text{and}\quad\,\, [u,\dd_l\ff(X_c)]\subset 
\dd_{j+l-1}\ff(X_c)\enspace,
\;\forall u\in \dd_j\ff(X_c).
\end{equation}
To prove equations on the left of  (\ref{filt_2}),
 let $u\in\ddx$ be a homogeneous 
differential operator such that $\text{\it weight}(u) + 
\text{\it  order}(u)=0,$
and  let $K(u)$ be its  kernel.
Since $K(u)$ is quasi-invariant with
respect to the `$x$'-variable, Lemma \ref{cv2}(i)
implies that there exists a differential
operator $p\in\dd(\hr)$ (in the variable $k\in \hr$)
such that $K(u)=p^{(k)}\psi_c$. 
Now, the proof of Proposition \ref{kernel1}(ii)
shows that $\text{\it weight}(u) + 
\text{\it  order}(u)=0$
 if and only if
$\deg^{(x)}(p^{(k)}\psi_c)= \deg^{(x)}\Phi$.
If the latter equation holds, then
Lemma \ref{cv2}(ii) forces the differential
operator $p\in\dd(\hr)$ to have order zero,
i.e., $p=p(k)\in \C[k]=\C[\h^*]$ is a multiplication operator.
Hence, $K(u)\cdot e^{(x,k)}=\bigl(p^{(k)}\psi_c\bigr)(x,k)=
p(k)\cdot\psi_c(x,k)=p(k)\cdot\Phi(x,k)\cdot e^{(x,k)}$.
Thus, we have proved that $\text{\it weight}(u) + 
\text{\it  order}(u)=0$ holds
 if and only if 
$K(u)$ has the form $p\cdot \Phi$, where
$p$ is a polynomial in the variable `$\,k\,$' only.
But then we must have $p\in Q_c$
and $u=\ll_{p}$. This completes the proof of
(\ref{filt_2}).\qed

Next, identify $\h^*$ with $\h$ via an invariant form,
and write
$P^\flat\in \C[\h\times\h^*]$ for the polynomial obtained from $P$
by the flip of  variables: $x\leftrightarrow k$.
It is clear that $P^\flat\in\kk$. We write $u^\flat:=K^{-1}(P^\flat)
\in \dd(X_c)$
 for the differential operator corresponding to
the kernel $P^\flat$.

Given a linear map $u: Q_c\to Q_c$, let
$u^\dag:  Q_c\to Q_c$ denote the adjoint map with respect
to the   nondegenerate form 
$\boldsymbol{(}-,-\boldsymbol{)}_c$ on $Q_c$.
 
\begin{lemma}\label{udag} 
\vi $\boldsymbol{(}
u^\flat(f),g\boldsymbol{)}_c =\boldsymbol{(} f,
u(g)\boldsymbol{)}_c,$
$\,\forall f,g\in Q_c$ and $u\in\ddx$. That is, $u^\flat=u^\dag$,
and the assignment: $u\mapsto u^\dag=u^\flat$ gives a
well-defined algebra anti-involution on $\ddx$.

\vii We have: $\ll_{_P}=P^\flat,$ for any $P\in Q_c$.
\end{lemma}

\begin{proof} 
By the results of [EG2],
it is sufficient to prove the identity $\boldsymbol{(}
u^\flat(f),g\boldsymbol{)}_c =\boldsymbol{(} f,
u(g)\boldsymbol{)}_c\,$ for $f(x)=\psi_c(x,a)\,,\,\forall a\in\h^*.$
Recall that $(u^{(x)})\ff\psi_c(x,a)=\bigl(u^{(k)}\psi_c(x,k)\bigr)|_{k=a}.$
Thus $\boldsymbol{(}
u^\flat(f),g\boldsymbol{)}_c=(ug)(a)=\boldsymbol{(}
f, u(g)\boldsymbol{)}_c,$
and the first equation of part (i) follows.
This yields  $u^\flat=u^\dag,$
for any $u\in\ddx.$
Thus, we have proved that
$u\in\ddx\;\Longrightarrow\; u^\dag \in\ddx\,.$
Since $u^\dag\ccirc v^\dag= (v\ccirc u)^\dag,$ 
part (i) follows. Part
(ii) follows from the formula $(\ll_p)^\dag=p\,,\,\forall p\in Q_c$,
proved in [EG2].
\end{proof}
\smallskip

\noindent
{\sl Proof of Theorem \ref{DDX}.\,\,} 
We prove first that
$ Q_c^\flat$ is the centralizer of the Calogero-Moser operator in
$\dd(\hreg)_\minus$. To this end, write
$\sigma(u)\in \C(\h)[\h^*]$ for the principal symbol of a
differential operator $u\in \dd(\BK)$.
Observe that, for any $p\in Q_c$, we
have: $\sigma(\ll_p)=p$. Now,
let $u\in \dd(\hreg)_\minus$ be a homogeneous operator
such that $[u,\LL_c]=0.$ We will
prove that $u\in Q_c^\flat$.  It suffices to show that, 
$\sigma(u)=p\in Q_c$.
Then, the operator  $u'= u-\ll_p$  commutes with $\LL_c.$
Since $\text{\it  order}(u')<\text{\it  order}(u)$ and
$\text{\it weight}(u') = \text{\it weight}(u)$,
we deduce that 
$\text{\it weight}(u') + 
\text{\it  order}(u')< 0$. Thus, $u'=0$ by (\ref{filt_2}),
and $u=\ll_p$.

Observe further that the leading term of $u$  must have
 constant coefficients (i.e.,
for the principal symbol we have: $\sigma(u)\in\C[\h^*]$) since
$\text{\it weight}(u) + 
\text{\it  order}(u)\leq 0.$ 
Hence $\sigma(u)=p,$ for some polynomial
$p\in \C[\h].$ Also, we claim that $u$ commutes with $Q_c^\flat.$ Indeed,
if $u'\in Q_c^\flat$ then for the operator $[u,u'],$
 the sum of order and degree is
strictly negative, so the operator must be zero. Thus,
$u$ preserves each of
the weight-spaces of $Q_c^\flat$ in 
the vector space $\bigoplus_{k\in\C}\;\C[\hr]\cdot e^{(k,x)}.$
These weight spaces are known, cf. e.g. [CV],
to be 1-dimensional and spanned by
the function
$\psi_c(-,k).$ It follows that $u^{(x)}\psi_c(x,k)=p(k)\psi_c(x,k).$
Therefore, $p$ is quasi-invariant.
Thus, we have  proved  that
the centralizer of $\LL_c$ 
 in
$\dd(\hreg)_\minus$ equals
$ Q_c^\flat$.

All the remaining claims of
 part (i) of the Theorem now follow from (\ref{filt_2}).
Part (ii) follows from Proposition \ref{kernel1}(i)
and Proposition \ref{kth2}.\qed
\medskip

\begin{remark}
 The adjoint action on $\ddx$ of the subalgebra
$Q_c=\C[X_c]$ is locally nilpotent by the definition of differential
operators. 
Since $u\mapsto u^\dag$ is an anti-homomorphism,
Lemma \ref{udag}(i) implies  that the adjoint action 
on $\ddx$ of the subalgebra  $Q_c^\flat$ is locally nilpotent.
This yields formula (\ref{P3}) in full generality.

An alternative proof of nilpotency of the adjoint
$Q_c^\flat$-action  on $\ddx$ follows from  formulas
(\ref{filt_2}), which yield:
$[u,\dd_l\ff(X_c)]\subset 
\dd_{l-1}\ff(X_c),$ for any $u\in Q_c^\flat$ and $l\geq 0$.
\hfill\qed
\end{remark}
\smallskip


Let $\dd_c$ denote the subalgebra in $\dd(X_c)$
generated by $Q_c=\C[X_c]$ and $Q_c^\flat$.
In the previous section we have introduced
 the subalgebra ${\mathcal{B}}_c$
generated by $Q_c^W$ and $(Q_c^\flat)^W$.
Thus we have 
\begin{equation}\label{inclusions}
\ehe\simeq {\mathcal{B}}_c\subset
\dd_c \subset \dd(X_c)\,.
\end{equation}
The inclusions in (\ref{inclusions})
give an algebra imbedding: $\ehe \into \ddx$, making
$\ddx$ an $\e\hh_c\e$-bimodule. 
\begin{proposition}\label{HCD}
\vi $\;\ddx$ is a Harish-Chandra $\ehe$-bimodule.

\noindent
\vii The group $SL_2(\C)$ acts on this Harish-Chandra bimodule
by algebra automorphisms; the  
Fourier transform 
interchanges $Q_c$ and $Q_c\ff$.
\end{proposition}

\begin{proof}
To show that $\ddx$ is a   Harish-Chandra $\ehe$-bimodule,
 recall that we have proved  that $\ddx\simeq \kk$,
as a  $\C[\h]^W\otimes \C[\h^*]^W$-module.
Further,  $\kk$, being a submodule
in $\C[\h\times\h^*]$, is clearly 
 a finitely generated $\C[\h]^W\otimes \C[\h^*]^W$-module.
It follows that $\ddx$ is a finitely generated
$\ehe\otimes (\ehe)^{op}$-module. By a remark above
(\ref{inclusions}) we know that the adjoint
action on $\ddx$ of any element of either $\C[\h]^W\subset Q_c$
or $\C[\h^*]^W\subset Q_c^\flat$ is locally nilpotent.
This proves (i).

To prove (ii) note that the Fourier automorphism on  $\ehe$
interchanges
$\C[\h]^W\e$ with $\C[\h^*]^W\e$, see [EG].
But $Q_c$ is obviously the centralizer of
$\C[\h]^W\e$ in $\ddx.$ Thus, the image of
$Q_c$ under the Fourier automorphism equals the centralizer
of $\C[\h^*]^W\e,$ which is exactly $Q_c^\flat.$
\end{proof}

Further, observe that the natural $W$-action on $X_c$ induces
a $W$-action on  $\dd(X_c)$ by algebra automorphisms. This
 $W$-action clearly commutes with
 both the left and  right action of $\e\hh_c\e$
on $\ddx$. 
Thus, we have a direct sum decomposition
$$
\dd(X_c)=\bigoplus_{\tau\in \text{Irrep}(W)}\tau^*\otimes D(\tau),
\quad\text{where}\quad
D(\tau)=(\tau\otimes \dd(X_c))^W
\quad\text{are}\quad \e\hh_c\e\text{-bimodules}.
$$
Note that, for any $\tau\in\irrep(W)$,
we have: $D(\tau)\neq 0$,
since $\,0\neq (\tau\otimes Q_c)^W\subset 
D(\tau)$.

\begin{theorem}\label{bimstructure}
The bimodules $\,\{D(\tau)\}_{\tau\in
\irrep(W)}\,,$ are 
pairwise non-isomorphic simple Harish-Chandra $\ehe$-bimodules,
for any $c\in\ccp$.
\end{theorem}

We start the proof  with the following lemma, which is perhaps well known. 

\begin{lemma}\label{kashiwara} Let $X$ be a smooth affine algebraic 
variety over $\C$ with a free action of a finite group $\G$. 
Let $\tau$ be an irreducible $\G$-module. 
Then the space $\dd_\tau(X):=\Hom_\G(\tau,\dd(X))$ is 
a simple bimodule over $\dd(X/\G)$.
\end{lemma}

\begin{proof} Using the equivalence between left and right 
$\dd$-modules, we can regard $\dd_\tau(X)$ as a $\dd$-module 
on $X/\G\times X/\G$. It is concentrated on the diagonal, and 
is the direct image of the irreducible local system on the diagonal, 
corresponding to the representation $\tau$ of $\G$. By Kashiwara's theorem, 
the direct image functor is an equivalence between the category of 
$\dd$-modules on the subvariety  and $\dd$-modules on the ambient variety 
supported on the subvariety. Thus, $\dd_\tau(X)$ is an irreducible $\dd$-module, 
as desired. 
\end{proof} 
 \smallskip

\noindent
{\sl Proof of Theorem \ref{bimstructure}.\,\,}
We know that $\ddx$, hence any its $W$-isotypic component $D(\tau)$,
 is a   Harish-Chandra $\ehe$-bimodule.
We first prove  simplicity of the bimodule
 $D(\tau)$, which is parallel to the proof 
of  simplicity of the algebra  $\hh_c$. 
Consider the localization $D(\tau)\mrw$
of $D(\tau)$.
This is a bimodule over $(\e\hh_c\e)\mrw=\dd(\hreg)^W$, 
which is isomorphic to $(\tau\otimes \dd(\hreg))^W$.

Let $J$ be a nonzero subbimodule of $D(\tau)$. Consider the 
Harish-Chandra $\ehe$-bimodule
$V=D(\tau)/J$, and its localization 
$V\mrw$.
By Lemma \ref{kashiwara}, $D(\tau)\mrw$ is a simple
$(\e\hh_c\e)\mrw$-module. Hence, 
$V\mrw=(D(\tau)\mrw)\big/(J\mrw)=0$.
Corollary \ref{cor_eHe} imples that $V=0$.
Hence $J=D(\tau),$ and $D(\tau)$ is a simple 
 $\ehe$-bimodule.

Finally, the bimodules $D(\tau)\,,\,\tau\in\irrep(W),$
are pairwise  nonisomorphic since
 the localized bimodules $D(\tau)\mrw=
(\tau\otimes \dd(\hreg))^W$ 
are already nonisomorphic.
The latter is clear since these modules, viewed as $\dd$-modules
on
$\hr\times \hreg$ give rise to  pairwise non-isomorphic local systems 
on the diagonal in $\hreg\!/W\times \hreg\!/W$. 
\qed

The following is an
analogue of the Levasseur-Stafford theorem
\begin{theorem}
\label{th_simple4} For any $c\in\ccp$, the algebra
$\dd(X_c) $ is generated by the two commutative subalgebras,
$Q_c$ and $Q_c^\flat$.
\end{theorem}

\begin{proof}
Recall, see (\ref{inclusions}), that $\dd_c$ denotes the subalgebra of 
$\dd(X_c) $ generated by $Q_c$ and $Q_c^\flat$.
Since $\e\hh_c\e$ is generated by $\C[\h]^W\!\!\cd\e$ and $\C[\h^*]^W\!\!\cd\e$, 
the subalgebra $\dd_c$ is a $\C W\otimes \e\hh_c\e\otimes 
(\e\hh_c\e)^{op}$-submodule 
in $\dd(X_c)$. The $\C W\otimes \e\hh_c\e\otimes 
(\e\hh_c\e)^{op}$-module decomposition:
$
\dd(X_c)=\oplus_{\tau\in \irrep(W)}\tau^*\otimes D(\tau),$
and Theorem \ref{bimstructure} imply that $\dd_c=\bigoplus_{\tau\in S}\;
\tau^*\otimes D(\tau),$ where summation runs over a certain
set $S\subseteq\irrep(W)$ of irreducible $W$-modules. 
However, $\dd_c$ contains $Q_c$, hence
 every  irreducible $W$-module occurs in $\dd_c$
with nonzero multiplicity. Thus, the set $S$ has to be 
the entire set $\irrep(W)$, and we are done. 
\end{proof}

\begin{theorem}
\label{th_simple3} The algebra
$\dd(X_c)$ is simple, for any $c\in\ccp$. 
\end{theorem}

\begin{proof} We begin with the following general

\noindent
{\bf Claim.\;} {\it
Let $ X $ be an irreducible affine variety, and
let $ A \subset  \C[X] $ be a subalgebra 
such that $ \C[X] $ is a finite module over $ A $.
Then any nonzero two-sided ideal of $ \dd(X) $ 
intersects $ A $ non-trivially.}

To prove the Claim,
let $ J $ be a non-trivial two-sided ideal of $ \dd(X) $. Choose
an nonzero differential operator, say $ u $, in $ J $ of the minimal
order. Since $ [u,f] \in J $ and $ \text{\it{order}}
([u,f]) < \text{\it{order}}(u), $ for any 
$ f \in \C[X] $, we conclude that $ [u,f] = 0 $ and hence $ u \in \C[X]$. 
It follows that $ u \in J \cap \C[X] $. 
Now, being finite as a module over $ A $, the ring $\C[X]$ is integral
over $ A $. Hence, $ u $ satisfies an equation 
$ u^n + a_{1}\cdot u^{n-1} + \ldots + a_{n-1}\cdot u + a_0 = 0, $ 
with coefficients in $ A $. In particular, there is a nonzero 
$ g \in \C[X] $ such that $ g\cdot u \subset A $. Since $ u \in J \;
\Longrightarrow\; g
\cdot u \in J $, 
we have $ J \cap A \not= {0}, $
and the Claim follows.

To prove the Theorem, assume
$J\subseteq \dd(X_c)$ is a nonzero two-sided ideal. 
The Claim above applied to 
$A=\C[\h]^W\subset \C[X_c],$ yields
$\C[\h]^W\cap J\neq 0$. Thus, $J$ 
nontrivially intersects $\e\hh_c\e$. 
But $\e\hh_c\e$ is simple, so $J\cap \e\hh_c\e=\e\hh_c\e$. 
Hence $1\in J,$  and $J=\dd(X_c).$
\end{proof}

The following corollary of Theorem~\ref{th_simple3}
has been  first
proved
in [EG2] by a different method.
\begin{corollary}
\label{CM}
For any $ c \in \Z_{+}(R)^W \,$, the variety $ X_c $ is 
Cohen-Macaulay.
\end{corollary}
\begin{proof}
Follows from a theorem of 
Van den Bergh (see \cite{VdB}, Theorem~6.2.5) stating
that $ X $ is a Cohen-Macaulay variety whenever the ring
$ \dd(X) $ is simple.
\end{proof}

\begin{remark}
$ X_c\, $ is Cohen-Macaulay $\;\Longleftrightarrow\;$
$ Q_c = \C[X_c]\, $ is a {\it free} $ \C[\h]^W$-module.
The freeness of $Q_c$ over $ \C[\h]^W$
 has been conjectured by Feigin-Veselov  and
proved  in \cite{FV} for dihedral groups,
and in \cite{EG2} in general.
\end{remark}

\begin{proposition} 
\label{ChD}
The imbedding: $\e\hh_c\e\into \dd(X_c)^W$ 
is
an algebra isomorphism.
\end{proposition}
\begin{proof} 
Let $u\in\ddx$ be a homogeneous
operator and  $K_0(u)$ be the highest term of 
the corresponding kernel
$K(u)
\in\kk,$ cf. Definition \ref{kk}.
It is clear that $K_0(u)(x,k)=\delta_c(x)\delta_c(k)
\cdot\text{\it Symbol}(u).$
The  principal symbol is always 
polynomial in $k,$ so $K_0(u)(x,k)$ is always divisible
by $\delta_c(k).$ But by symmetry of $x\longleftrightarrow
k,$ it is also always divisible by
$\delta_c(x).$ Thus, we have proved that the principal
symbol of any differential operator $u\in \ddx$ is polynomial, and 
hence
$ \grd\bigl(\dd(X_c)^W\bigr) \subset 
\C[\h^*\times \h]^W \,$. 

On the other hand, $\,
\grd(\e\hh_c\e)=\C[\h^*\times \h]^W$. 
Hence, the associated graded map: $\grd(\e\hh_c\e)\to
\grd(\dd(X_c)^W)$ is surjective.
The Proposition follows.
\end{proof}

\section{Translation functors and Morita equivalence}
\setcounter{equation}{0}

\begin{theorem}\label{tran} Let $c\in \ccr$.
Then, for any 
$m\in \Z[R]^W$,  the algebras $\hh_{c}$ and $\hh_{c-m}$
(resp. the algebras $\e\hh_{c}\e$ and $\e\hh_{c-m}\e$)
 are Morita equivalent.
\end{theorem}

\begin{proof} We may assume without loss of generality that
$m\in\ccp$.
Recall that the algebra $\hh_{c}$ is simple for any $c\in\ccr$,
therefore both $\hh_c$ and $\hh_{c-m}$ are simple.
Thus, using Proposition \ref{H_minus} and Corollary
\ref{cor_simple1} 
we deduce 
$$
{
\diagram
\hh_{c}\enspace
\rrdouble^<>(.5){_{{\tt Morita\; equivalence}}}_<>(.5){{\tt Lemma\;
 \ref{mor_simple}}}
&&
\enspace\e_{\varepsilon}\hh_{c}\e_{\varepsilon}
\enspace\rrdouble^<>(.5){\tt Prop.\; \ref{H_minus}}
&&\enspace \e\hh_{c-{\mathbf{1}_\varepsilon}}\e 
\enspace\rrdouble^<>(.5){_{{\tt Morita}}}_<>(.5){^{\tt equivalence}}
&&\enspace
\hh_{c-{\mathbf{1}_\varepsilon}}\,.
\enddiagram}
$$
Recall further that the functions of the
form ${\mathbf{1}_\varepsilon},$ where
$\varepsilon: W\to \{\pm 1\}$ is
a multiplicative character, see above Proposition \ref{H_minus},
are  known to generate $\ccp$ as a semi-group.
Thus $m$ is a sum of functions
of type ${\mathbf{1}_\varepsilon}$, and an easy induction
on the number of summands completes the proof.
\end{proof}

\begin{lemma} For any $c\in \Z[R]^W$,  the algebra $\hh_{c}$ 
is  Morita equivalent to $\dd(\h)\#W.$
\end{lemma}

\begin{proof}  For  $c\in \Z[R]^W$, Theorem \ref{tran} yields  Morita equivalence
of the algebras $\hh_c$ and $\hh_0$. But $\hh_0\simeq \dd(\h)\#W$,
and the result follows.
\end{proof}

\begin{remark} If $\h=\C$ and $W=\Z/2\Z,$ then $\e\hh_c\e$ is 
a quotient of $U(\sll2)$, and Theorem \ref{tran}
is known and due to Stafford [St].\qed
\end{remark}
\smallskip

By analogy with Lie  theory, 
we  denote by $\HC_c(\hh)$
the category of Harish-Chandra $\hh_c$-bimodules.
In a similar way, we may consider the  category
$\HC_c(\eHe)$ of Harish-Chandra
bimodules over the algebra $\ehe$. 
Recall that the category of bimodules over any algebra $A$
has a monoidal structure given by the tensor product of
bimodules over $A$. 

\begin{lemma}\label{tens} The tensor product of two  Harish-Chandra
bimodules is again a Harish-Chandra
bimodule. 
\end{lemma}

\noindent
{\sl Proof.\,\,} Let $V,V'\in \HC_c(\hh)$.
The Leibniz rule for derivations shows that the adjoint action on
$V\otimes_{\hh_c} V'$ of any element of  $\C[\h]^W$ 
and of any element of
$\C[\h^*]^W$ is locally nilpotent.
Thus, it suffices to prove that $V\otimes_{\hh_c} V'$ 
is a finitely generated $\hh_c\otimes \hh_c^{op}$-module.

To this end, choose  a finite set $\{v_i\}$, resp.  $\{v'_j\}$,
of generators of $V$, resp. of $V'$, 
as  a $\C[\h]^W$-$\C[\h^*]^W$-bimodule, see Lemma \ref{straight}(ii).
Then we get
\begin{equation}\label{V_V}
V\otimes_{\hh_c}V'=
\sum\nolimits_{i,j}\;\bigl(\C[\h]^W\cdot v_i\cdot
 \C[\h^*]^W\bigr)\otimes_{\hh_c}
\bigl(\C[\h]^W\cdot v'_j\cdot \C[\h^*]^W\bigr)\,.
\end{equation}

Now
Lemma \ref{tens2} insures, due to nilpotency of the adjoint action,
that
the subspace $\C[\h]^W\cdot v\cdot \C[\h]^W\subset V$
  is finitely generated as 
a  left $\C[\h]$-module.
Similarly, $\C[\h^*]^W\cdot v'\cdot \C[\h^*]^W$ is finitely generated as 
a   right $\C[\h^*]^W$-module, for any $v'\in V'$.
We conclude that there are finite sets
$\{u_{ik}\}\subset V$ and $\{u'_{jl}\}\subset V',$
such that: $ \C[\h]^W\cdot v_i\cdot \C[\h]^W\subset\sum_k\,\C[\h]^W\cdot u_{ik}$,
and $\C[\h^*]^W\cdot v'_j\cdot \C[\h^*]^W\subset
\sum_l\,u'_{jl}\cdot \C[\h^*]^W$.

Further,
 by the Poincar\'e-Birkhoff-Witt isomorphism (\ref{pbw}),
there exists a finite dimensional subspace
$E\subset \hh_c$ such that $\C[\h^*]^W\otimes_{\hh_c}\C[\h]^W
\subset \C[\h]^W\cdot E\cdot\C[\h^*]^W\,$
(inclusion of subspaces in $\hh_c\otimes_{\hh_c}\,\hh_c$).
The proof is now completed by the following
 inclusions
\begin{align*}
\text{RHS of } (\ref{V_V}) &\subseteq
\sum\nolimits_{i,j}\;\C[\h]^W\cdot v_i\cdot 
\C[\h]^W\cdot E\cdot\C[\h^*]^W
\cdot v'_j\cdot \C[\h^*]^W\bigr)
\quad\\
&\subseteq
\sum\nolimits_{i,j,k,l}\;\C[\h]^W\cdot u_{ik}\cdot E\cdot
u'_{jl}\cdot\C[\h^*]^W\;\subseteq\;\sum\nolimits_{i,j,k,l}\;
\hh_c\cdot u_{ik}\cdot E\cdot
u'_{jl}\cdot\hh_c\,.\qed
\end{align*}\smallskip

Thus, for each $c\in\CC,$ the category  $\HC_c(\hh),$
 resp. $\HC_c(\eHe)$, has a monoidal structure
induced by the tensor product of bimodules. 
In case of $c\in\ccp$ a complete structure of these
 monoidal categories  can be described as follows.
First, recall that for any finite dimensional
(not necessarily irreducible)
$W$-representation $\tau$, the space
$D(\tau):= (\tau\otimes\ddx)^W$ is
a Harish-Chandra $\ehe$-bimodule, cf. Theorem \ref{bimstructure}.
 Let $\Rep(W)$ denote the tensor category of finite dimensional
representations of the group $W$.

\begin{theorem} \label{equiv}
For any $c\in\ccp$,
the assignment
$\tau\mapsto D(\tau):= (\tau\otimes\ddx)^W$
gives an equivalence
 $\Xi: \Rep(W) \iso \HC_c(\eHe)$
 of monoidal categories. 
The monoidal structure on the functor $\Xi$
is given,
for any $\tau,\sigma\in\irrep(W),$
by a  canonical
$\ehe$-bimodule isomorphism:
$D(\tau)\otimes_{\ehe}D(\sigma)\iso
D(\tau\otimes\sigma)$ induced by 
multiplication in the algebra $\ddx$:
$$(\tau\otimes\ddx)^W\;\bigotimes\nolimits_{\ehe}\;
(\sigma\otimes\ddx)^W\,
\iso\,\bigl((\tau\otimes\sigma)\otimes\ddx\bigr)^W\,.
$$

Similarly,
there is an equivalence of monoidal categories
$\Rep(W) \iso \HC_c(\hh)$.
\end{theorem}

 \begin{proof} We first consider the case $c=0$.
 In this case
we have $\hh_0\simeq \dd(\h)\#W$
and $\ehe\simeq \dd(\h)^W$.

Since $\hh_0=\dd(\h)\#W$,
 a Harish-Chandra  $\hh_0$-bimodule  is in particular a
$\dd(\h)$-bimodule with a  $W\times W$-action,
that is a  $W\times W$-equivariant $\dd$-module, $V$,
on $\h\times\h$.
The bound given 
in (\ref{graph}) shows that $\Ch V$, the {\it characteristic variety}
of  $V$ (as a $\dd(\h\times\h)$-module),
is contained
in $\Bigl(\cup_{w\in W}\;\text{Graph}(w)\Bigr)
\times\Bigl(\cup_{y\in W}\;\text{Graph}(y)\Bigr),$
which is a subvariety in $T^*(\h\times\h)$ of dimension $\frac{1}{2}\dim
T^*(\h\times\h)$.
Hence,  $V$ is a holonomic $\dd$-module.
Moreover, since $\Ch^{\,}V$ must be a Lagrangian subvariety, we deduce
\begin{equation}\label{Ch}
\Ch^{\,}V\subseteq \bigcup\nolimits_{w\in W}\;
T^*_{_{\text{Graph}(w)}}(\h\times\h)\,.
\end{equation}

Given $w\in W$, let $i_w: \text{Graph}(w)\into \h\times\h$
denote the imbedding.
Kashiwara's theorem
and estimate (\ref{Ch}) imply that every irreducible subquotient of 
$V$, viewed as a $\dd$-module,  must be of the form
$(i_w)_*E$, where $E$ is an irreducible local system, i.e., a vector bundle
on $\text{Graph}(w)$ with flat connection.
The condition of $\ad$-nilpotency of the action on $(i_w)_*E\subset V$ of
the algebra $\C[\h^*]^W$ forces the connection to be trivial.
 We deduce that all  irreducible subquotients
 of 
$V$ are  of the form
$(i_w)_*\C[\text{Graph}(w)]$.
Semisimplicity  of the
category  $\HC_0(\hh)$ follows from this, since
it is well-known that
there are no extensions between $\dd$-modules of 
the form above.
We conclude that the category $\HC_0(\hh)$
is semisimple.

We have shown that any
object $V\in \HC_0(\hh)$ can be uniquely written, viewed
as a $\dd(\h\times\h)$-module, in the form
$V=\bigoplus_{w\in W}\; V_w\otimes_{_\C}\, (i_w)_*\C[\text{Graph}(w)],$
where $V_w$ are finite dimensional vector spaces. Moreover, 
the left and right actions on $V$
of the subalgebra $\C W\subset\hh_0$ induce,  
for any $y\in W,$ the natural isomorphisms
$L_y: V_w\iso V_{yw}$ and $ R_y: V_w\iso V_{wy},$ which are multiplicative,
and 
commute:
$R_yL_z=L_zR_y\,,\,\forall y,z\in W.$
 This means that the vector space $V_W:=\oplus_{w\in W}\, V_w$ is a 
representation of $W\times W,$ via 
$(y_1,y_2)v=L_{y_1}R^{-1}_{y_2}v.$ It is clear that 
$V_W=\Ind_{W_{\!_{\sf diag}}}^{W\times W}V_1,$ 
where the action  on $V_1$ of the diagonal subgroup
$W_{\!_{\sf diag}}\subset W\times W$  is given
by $w: v\mapsto L_w R^{-1}_w(v)$, as
above. Moreover, we see that $V \simeq
\displaystyle
\Ind_{A}^{\hh_0\otimes \hh_0^{op}}\!\Bigl((i_{_{\sf diag}})_*\bigl(V_1\otimes
\C[\h_{_{\sf diag}}]\bigr)\Bigr),$ where $A$ is the
preimage of $\C[W_{\!_{\sf diag}}]$ under the natural
 projection $\displaystyle\hh_0\otimes
\hh_0^{op}
=\dd(\h\times\h)\#(W\times W)\onto 
\C[W\times W].$ Here $V_1\otimes
\C[\h_{_{\sf diag}}]$ is a $W$-equivariant $\dd(\h)$-module, 
and $i_{_{\sf diag}}: \h_{\!_{\sf diag}}=\h\into\h\times\h
$ is the diagonal map. 

Next, recall the following general result on 
 Morita equivalence.
Let $A$ be an algebra and $\e\in A$ an idempotent,
such that $A\e A=A$.
Then, {\it the natural Morita equivalence
between $A$-bimodules and $\e A\e$-bimodules
is an equivalence of monoidal categories.}

In the special case $A=\hh_c$
the result above
implies that the categories $\HC_0(\hh)$ and $\HC_0(\eHe)$
are equivalent as monoidal categories, due to Morita equivalence
of the algebras $\dd(\h)\# W$ and $\dd(\h)^W$
proved in [Mo].
In particular, the category $\HC_0(\eHe)$
is also  semisimple.

 Further it is easy to see,
e.g. ([LS2], \S\S2-3) or Theorem \ref{bimstructure}
(in the special case  $c=0$)
 that, for each $\tau\in \irrep(W)$,
the $W$-isotypic
component
$\bigl(\tau\otimes_{_\C}\dd(\h)\bigr)^W$
is a simple $\ehe$-bimodule.
Moreover, any simple object of  $\HC_0(\eHe)$
is isomorphic, by Morita equivalence,
 to a bimodule of this type.
Thus, simple objects of the category $\HC_0(\eHe)$
 are
parametrized by the set $\irrep(W)$. This completes the proof
of the  Theorem
for $c=0$.

 In the general case
of an arbitrary $c\in\ccp$, the result follows from that
for $c=0$ via
Morita equivalence of Proposition \ref{transl} below
(which is independent of the
intervening material).
\end{proof}

\begin{corollary}
Both $\HC_c(\eHe)$
 and  $\HC_c(\hh)$ are
semisimple tensor categories.\qed
\end{corollary}

To proceed further, it will be useful to generalize
the setup and to consider bimodules over two Cherednik
algebras with possibly {\it different} parameters. Thus,
we fix $c,c'\in \CC$, and let $\hh_c\,,\,\hh_{c'}$ be
the corresponding Cherednik algebras.
There are canonical imbeddings of both
$\C[\h]^W$ and $\C[\h^*]^W$ into any of the two
algebras $\hh_c$ or $\hh_{c'}$.
We have, cf. Definition \ref{HC_def}

\begin{definition}\label{HC_def2} A finitely generated 
 $\hh_c\otimes\hcc^{op}$-module $V$ is called
a Harish-Chandra  $\hhcc$-bimodule  if,
for any $a\in \C[\h]^W$ and any $a\in \C[\h^*]^W$, the $\ad a$-action
on $V$ is locally nilpotent.
\end{definition}

We write $\HC_{c,c'}(\hh)$, resp. $\HC_{c,c'}(\eHe)$,
for the category of Harish-Chandra  $\hhcc$-bimodules,
resp. $\ehhcce$-bimodules. Most of the results 
on Harish-Chandra bimodules proved in \S\S3-4 extend
verbatim to this `two-parameter' setting.

Let $V$ 
be an $\hhcc$-bimodule,  resp. $\ehhcce$-bimodule.
  We say that a vector $v\in V$ is of order
$\le d$ if for any $f\in \C[\h]^W$ and $g\in \C[\h^*]^W$, one has 
$(\ad f)^{d_1+1}v=(\ad g)^{d_2+1} v=0$, 
for some $d_1,d_2$ such that $d_1+d_2=d$.
The space of elements of finite order in $V$ is
a $\hhcc$-sub-bimodule, to be denoted by
$V\fin$. 
Clearly, if $V$ is a finitely generated 
bimodule, then  $V\fin$ is the maximal Harish-Chandra  subbimodule in $V$. 
For example, $\dd(X_c)\fin=\dd(X_c)$.

Further, given a {\it left} $\hh_c$-module $M$,
and a {\it left} $\hcc$-module $N$,
the vector space $\Hom_{_\C}(M,N)$ has a canonical
$\hhcc$-bimodule structure, and we put
$\,\Hom\fin(M,N) := \bigl(\Hom_{_\C}(M,N)\bigr)_{\sf{fin}}$.
Similar definitions apply to $\ehhcce$-bimodules.

Since we will be dealing with standard modules over Cherednik
algebras with varying parameter `$c$', it will be
convenient to incorporate this parameter in the notation
and to write $M(\tau,c)$ for the standard module
over $\hh_c$, and $M_\e(\tau,c):=\e\cd M(\tau,c)$ for the
corresponding  standard module
over $\ehe$. Recall that $\triv$ denotes the
trivial 1-dimensional representation of the group $W$.

For any $\tau\in\irrep(W)$ and $c,c'\in\CC$, we
define an object $\lc \pp\rc(\tau)\in \HC_{c,c'}(\eHe)$ by
\begin{equation}\label{irr_def}
 \lc \pp\rc(\tau):= \Hom\fin\bigl(M_\e(\tau^*,c')\,,\,M_\e(\triv,c)\bigr)\,.
\end{equation}

\begin{proposition}\label{hfin}\vi
 For any $c\in\ccr,$ the $\ehe$-action on 
$M_\e(\triv,c)$ induces
a natural $\ehe$-bimodule isomorphism:
$\displaystyle\ehe\iso \lc \pp_{\!c}(\triv)$.

\vii The space $\lc \pp\rc(\triv)$ is nonzero if and only if 
$c-c'$ is an integer valued function. 

\viii For any $c,c',c''\in \ccr$
such that $c-c'\,, \,c'-c''\in\Z[R]^W$,
there is a canonical bimodule isomorphism
$$\lc\pp\rc(\triv)\;\bigotimes\nolimits_{\e\hh_{c'}\e}\;
{{}_{c'}\pp_{\!c''}(\triv)}
\iso \lc\pp_{\!c''}(\triv)\,.$$
\end{proposition}

\begin{proof} Part (i) is formulated here
for the reader's convenience only;
it is a special case of the more general Theorem
\ref{deco} formulated and proved below.

To prove (ii) we use the {\it shift operator}
$\SH_{c',c'+m}$, see [O] and also \cite{EG2}.
Suppose that $c-c'$ is integer valued, and let us show
that $\lc\pp\rc(\triv)$ is nonzero. First assume that $c=c'+m$, where 
$m\in\ccp$ is nonnegative. Then the shift operator
$\SH_{c',c'+m}$ is a (nonzero) element of
$\lc\pp\rc(\triv)$. Indeed, it is of ``degree 0'' 
with respect to the algebra of Calogero Hamiltonians, and of finite degree 
with respect to the algebra of functions (as it is a differential operator).
On the other hand, if $m\le 0$, i.e. if
$(-m)\in\ccp$, then there also exists a shift operator 
$\SH_{c',c'+m}$, which however does not preserve the space of 
symmetric polynomials (it has poles). Nevertheless, 
the operator $\delta(x)^{2N}\SH_{c',c'+m}$ is a nonzero element 
of $\lc\pp\rc(\triv)$ for large enough $N$. Since any 
integer valued
function $m$ is a sum of a nonnegative and a
nonpositive function, 
we find that $\lc\pp\rc(\triv)$ is nonzero for any integer valued 
$c-c'$. 

Conversely, if $\lc\pp\rc(\triv)\ne 0$ then it must have elements of degree $0$ 
under the adjoint action of Calogero hamiltonians. Such an element 
satisfies the axiomatics of a shift operator of Opdam. But Opdam showed that 
shift operators exist only if $c-c'$ is integer: indeed,  
if a shift operator exists then the monodromies of the Calogero system for 
parameters $c$ and $c'$ are the same, i.e. $e^{2\pi ic}=e^{2\pi ic'}$. 
Part (ii) is proved. 

To prove (iii), 
denote the left hand side of the isomorphism  by $H$.
It is clear that the action of any element
$u\in H$  gives a  map
$\hat{u}: M_\e(\triv,c'')\to M_\e(\triv,c)$
such that
$\hat{u}\in\Hom\fin\bigl(M_\e(\triv,c''),\bigr.$
$\bigl.M_\e(\triv,c)\bigr).$
This way we obtain a natural bimodule homomorphism $j: H
\to \lc\pp_{\!c''}(\triv).$
Localizing at the subset $\hr/W\subset \h/W$ we get
a homomorphism $j\mrw: H\mrw\to \lc\pp_{\!c''}(\triv)\mrw.$
The latter morphism is independent of the choice of
the parameter `$c$', and is easily seen to be bijective for
$c=0$. Hence $j\mrw$ is a bijection.
Next, we use an analogue of Corollary 3.6 (for bimodules over
Cherednik
algebras with two possibly different parameters `$c$').
The proof of this result is entirely analogous to
the proof  of Corollary 3.6. It follows from the result
that the kernel of $j$ has  no torsion.
We conclude that $j$ is  injective.
Similarly, the  {\it co}kernel of $j$  has  no torsion
by the same analogue of Corollary \ref{key}. 
Thus, $j$ is surjective, hence bijective.
\end{proof}

\begin{remark}
The map $\ehe\to \Hom\fin\bigl(M_\e(\triv,c)\,,\,M_\e(\triv,c)\bigr)$
of Proposition \ref{hfin} is also an {\it algebra} isomorphism
with respect to the composition structure on the $\Hom$-space
on the right.
\end{remark}
\smallskip

We now give an alternative,
stronger and more effective version
 of 
Morita equivalence of
Harish-Chandra categories that 
arises from Morita equivalence 
of algebras $\ehe$ and $\e\hh_{c-m}\e$ provided
 by  Theorem \ref{tran}.

\begin{proposition}\label{transl}Let $c'\in \ccr$ and $c=c'+m$, where 
$m\in\ccp$. 

\vi The equivalence of categories  $\HC_c(\eHe)$
 and
$\HC_{c'}(\eHe)$ arising from Morita equivalence 
of algebras $\ehe$ and $\e\hh_{c'}\e$ 
is given by the functor: 
$V \mapsto {{}_{\!c'}\pp_c}(\triv)
\otimes_{_{\ehe}}V\otimes_{_{\ehe}}{{}_c\pp_{\!c'}(\triv)}$.

\vii Moreover, this  functor
is an equivalence of monoidal categories.
\end{proposition}
\begin{proof} 
Clearly, it suffices to consider the case
$m={\mathbf{1}}_\varepsilon$. In that case the Morita equivalence of
of algebras $\ehe$ and $\e\hh_{c'}\e$ comes from
the identification $\e_\varepsilon \hh_c \e_\varepsilon =
\e\hh_{c-{\mathbf{1}}_\varepsilon}\e$ of Proposition
\ref{H_minus}.
We will use an isomorphism of bimodules:
$\e \hh_c \e_\varepsilon\simeq
\Hom\fin\bigl(\e_\varepsilon\cd M(\triv,c)\,,$
$\,\e\cd  M(\triv,c)\bigr),$
which is proved similarly to the case where $\e_\varepsilon$
is replaced by $\e$, see 
Proposition \ref{hfin}. Further, it is easy to see that under 
the identification $\e_\varepsilon \hh_c \e_\varepsilon =
\e\hh_{c-{\mathbf{1}}_\varepsilon}\e$ of Proposition
\ref{H_minus}, the module $\e_\varepsilon\cdot  M(\triv,c)$ 
goes to $\e\cd  M(\triv, c-{\mathbf{1}}_\varepsilon).$
So, for $m={\mathbf{1}}_\varepsilon,$ the equivalence
of Harish-Chandra categories arising from Theorem \ref{tran}
is provided by the bimodule
$\Hom\fin\bigl(\e\cdot  M(\triv,c-m)
\,,\,\e\cdot  M(\triv,c)\bigr)$
which is exactly the bimodule $\lc \pp\rc(\triv)$.
Part (i) follows.

To prove (ii) one has to verify the
isomorphism: $\lc \pp\rc(\triv)\;\bigotimes\nolimits_{\e\hh_{c'}\e}\;
{{}_{c'}\pp_{\!c}(\triv)}=\lc\pp_{\!c}(\triv)$, as well as
an associativity constraint:
$$\lc\pp\rc(\triv)\;\bigotimes\nolimits_{\e\hh_{c'}\e}\;
{{}_{c'}\pp_{\!c''}(\triv)}
\iso \lc\pp_{\!c''}(\triv)\,.$$
The first isomorphism  follows from
the equality $\e\hh_c\e_\varepsilon\hh_c\e=\ehe$,
by simplicity of the algebra $\ehe$,
and the second is part (iii) of Proposition
\ref{hfin}.
\end{proof}

Let $\varepsilon: W\to \lbrace{\pm 1\rbrace}$
be a character. 
Then the algebra $\hh_c$ 
is isomorphic to $\hh_{\varepsilon c}$.
Thus, we have established the ` {\sl if} ' part of the following 

\begin{conjecture} Let $c\in\ccr$. Then
the algebras $\e\hh_{c}\e$ and $\e\hh_{c'}\e$,
resp.  $\hh_{c}$  and $\hh_{c'}$,
are Morita equivalent if and only if 
there exists a character $\varepsilon: W\to \lbrace{\pm 1\rbrace}$
such that $c-\varepsilon c'\in\Z[R]^W.$
\end{conjecture}

For $\h=\C$ and $W=\Z/2\Z$, the conjecture is a known result due 
to Hodges \cite{Ho}. Also, for $W=S_n$ and transcendental $c$, 
this conjecture has been proved in \cite{BEG2}.

Recall that if we write $\BK=\C(\h)
$ for the field of rational functions
on
$\h$, then we have $\ddx=\{u\in\dd(\BK)\;\mid\; u(Q_c)\subset Q_c\},$
see (\ref{3}).
Motivated by this,  for any $c, c'\in\ccp$, we put
$\dd(Q_{c'},Q_c):= \{u\in\dd(\BK)\;\mid\; u(Q_{c'})\subset Q_c\}.$
The space $\dd(Q_{c'},Q_c)$
has a natural right $\dd(X_{c'})$-module structure, and left 
$\ddx$-module structure. Generalizing the proof of simplicity
of the algebra $\ddx$ one obtains

\begin{proposition}\label{bimod_simpl}
The $\ddx\mbox{-}\dd(X_{c'})$-bimodule
$\dd(Q_{c'},Q_c)$ is simple, $\,\forall c, c'\in\ccp.$~\qed
\end{proposition}

We can now consider the vector space
${\mathscr{K}}_{c,c'}$ of "kernels"  with parameters $c,c'$.
Specifically, we put
$$
{\mathscr{K}}_{c,c'}
=\left\{ P\in \C[\h\times \h^*] \quad\Big|\quad
\begin{array}{l}(x,k)\mapsto P(x,k)\cdot e^{(x,k)}
\quad\text{\footnotesize is $c$-quasi-invariant 
with}\\
\text{\footnotesize respect to $x$, and $c'$-quasi-invariant 
with respect to $k$}
\end{array}\right\}.
$$
Similarly to the 
case $c=c'$, one proves 
\begin{proposition}\label{c_c'} For any $c,c'\in\ccp$,

\vi The space $
{\mathscr{K}}_{c,c'}$ is a rank one projective $Q_c\otimes
Q_{c'}$-module.

\vii There is a natural  $Q_c\otimes
Q_{c'}$-module isomorphism:
$\dd(Q_{c'},Q_c)\iso {\mathscr{K}}_{c,c'}$.\qed
\end{proposition}

Further, the group $W$ acts naturally on $\dd(Q_{c'},Q_c)$ commuting
with
$\ehhcce$-bimodule structure. This gives a
$W$-module  and $\ehhcce$-bimodule decomposition 
$$
\dd(Q_{c'},Q_c)=\bigoplus_{\tau\in \irrep(W)}\tau^*\otimes
D_{c,c'}(\tau)
\quad\text{where}\quad
D_{c,c'}(\tau):= \bigl(\dd(Q_{c'},Q_c)\otimes\tau\bigr)^W\,.
$$ 

\begin{theorem}\label{deco}
\vi For any $\sigma,\tau\in\irrep(W),$ there is a natural
$\ehhcce$-bimodule isomorphism
$$
\Hom\fin(M_\e(\sigma_{c'}',c'),M_\e(\tau_c',c))=\bigoplus_{\xi\in \irrep(W)} 
\Hom_W(\xi\otimes \sigma,\tau)\,\otimes\, D_{c,c'}(\xi)\,.
$$
\vii The $\ehhcce$-bimodules $D_{c,c'}(\tau)\,,\,\tau\in\irrep(W),
$ are irreducible, and 
pairwise non-isomorphic. 
\end{theorem}

\begin{proof} To prove (i), denote the right hand side
of the isomorphism by $V$
and the left hand side by $V'$. 
Then $V=\bigl(\dd(Q_{c'},Q_c)\otimes \Hom_\C(\sigma,\tau)\bigr)^W.$ 
On the other hand, by Proposition \ref{isot}, 
$M_\e(\tau_c')=(\tau\otimes Q_c)^W$. 
For any $v\in \Hom_W(\xi\otimes \sigma,\tau)\,\bigotimes\, 
\bigl(\dd(Q_{c'},Q_c)\otimes\xi\bigr)^W\subset V$, the natural
action-map
$v: M_\e(\sigma_{c'}',c')\to M_\e(\tau_c',c)$ gives rise to
a bimodule morphism $\theta: V\to V'$.
The morphism $\theta$  is injective
since, for  any $v\in V\,,\,v\neq 0,$
the induced action-map
$v: M_\e(\sigma_{c'}',c')\mrw\to M_\e(\tau_c',c)\mrw$
is given by a differential operator,
which is clearly nonzero.
Thus,  $V$ may (and will) be viewed as
 a subbimodule of $V'$. 

Now, consider the localized  bimodules 
$V\mrw$ and $V'\mrw$
over the localized algebra $\ehe\mrw$.
We have:
$\ehe\mrw
\simeq  \dd(\h)^W\mrw$,
and 
$V\mrw\simeq $
$\bigl(\dd(\h)\otimes \Hom_\C(\sigma,\tau)\bigr)^W\mrw$.
 It is clear
(from nilpotency of the adjoint $\C[\h]^W$-action  on $V'$)
that any element of $V'\mrw$ 
is represented by a differential operator, so $V'\mrw=V\mrw.$
Thus, $(V'/V)\mrw=
(V'\mrw)\big/(V\mrw)=0,$
and $V'/V=0$, by an  analogue of Corollary \ref{cor_eHe}.
Part (i) is proved.
Proof of part (ii) is entirely analogous to the
proof of Theorem \ref{bimstructure}.
\end{proof}

\begin{remark} Under the equivalence of Theorem \ref{equiv}, 
the algebra $\dd(X_c)$ corresponds to an algebra 
$A$ in the tensor category
$\text{Rep}(W)$ isomorphic to the regular representation of $W$. 
It is easy to show that in fact $A=\text{Fun}(W)$. 
\end{remark}

Recall the locally finite $SL(2)$-action 
on $\dd(X_c)$. When we pass to $\grd\bigl(\dd(X_c)\bigr)$, this action extends to a 
$\GL(2)$-action, since the order operator is added to the Lie algebra.

\noindent
{\bf Question.} What is the $GL(2)$-character of 
each isotypic component of $\grd\bigl(\dd(X_c)\bigr)$ ?

This question  is equivalent to the problem of computing the 
two-variable Hilbert series 
of isotypic 
components 
$\grd\bigl(\dd(X_c)\bigr)
$, with respect to gradings both in $\h$ and in $\h^*$.  
Also, by the above results, it is equivalent to the question of computing 
the ``characters'' of all Harish-Chandra modules.  
It would be interesting to use the bijection 
$\dd(X_c)\iso \kk$ in order to find the character of
$\grd\bigl(\dd(X_c)\bigr)$.
\smallskip

Since $\e\hh_c\e$ is Morita equivalent to $\hh_c$, there is an
$\ehe$-counterpart of
Theorem \ref{deco}. Specifically, one defines
$\hh_c$-bimodules $V(\tau):=\hh_c\e\otimes_{\e\hh_c\e}D(\tau)
\otimes_{\e\hh_c\e}\e\hh_c$ and proves

\begin{theorem}\label{deco1} The bimodules $V(\tau)$ are irreducible
and pairwise non-isomorphic. 
Moreover, for any  $\sigma,\tau\in\irrep(W),$
one has a $\hh_c$-bimodule isomorphism
$$
\Hom\fin(M(\sigma_c'),M(\tau_c'))=\bigoplus_{\xi\in \irrep(W)} 
\Hom_W(\xi\otimes \sigma,\tau)\otimes V(\xi)\,.\qquad\qed
$$
\end{theorem}

\section{Applications of the shift operator}
\setcounter{equation}{0}

In this section we prove some further results about
the structure of the variety $ X_c $ and the algebra 
$ \dd(X_c)\,$ using the notion of {\it a shift operator} 
introduced by Opdam. Recall that, given an {\it integral}
$\, c \,$, a shift operator $ \SH_c $ is defined as a 
differential operator in $ \dd(\h) $ that satisfies the 
relation $ \LL_{c}\circ \SH_c = \SH_c \circ \LL_0 \,$ and has 
the principal symbol $ \delta_c(x)\,\delta_c(k) \,$. 
We remark that the operator $\SH_c$ determines (and is uniquely determined by)
the Baker-Akhiezer function $\psi_{c}(x,k)$, introduced in \S7,
via the equation $\psi_c(x,k) = \SH_c[e^{(x,k)}]$.
Below, however, we prefer to work with  the shift operator rather than with 
$\psi_{c}(x,k)$ directly.

It follows from the results of \cite{O} and \cite{He} 
(see also \cite{CV}) that such an operator exists and is 
unique for any $\, c\in \ccp \,$. In fact, $ \SH_c $ can be
given explicitly by the following formula (cf. \cite{Be}, 
(2.11)):
\begin{equation}
\label{shift}
\SH_c =\frac{1}{2^{d}\,d!}\,
\ad(\LL_c, \LL_0)^{d}\,[\delta_c(x)^{2}]\ , 
\end{equation}
where $ d := \sum_{\alpha \in R_{+}} c_{\alpha} $ and
$ \ad(\LL_c, \LL_0) $ denotes the adjoint action 
on $ \dd(\hr) $ defined by multiplying the elements of
$ \dd(\hr) $ by $ \LL_c $ on the left and 
by $ \LL_0 $ on the right. 

Our applications of the shift operator $ \SH_c $ will be 
based on the following observation due to Opdam~\cite{O} 
(see also \cite{He}, (1.21)).
\begin{lemma}
\label{ng}
$\ \SH_c[1] \in \C \setminus \{0\} \,$ for any $\, c\in \ccp \,$.
\end{lemma}

First, we use $ \SH_c $ to give 
{\sl another proof of bijectivity 
of normalization $\,\pi_c: \tilde{X}_c \onto X_c \,$ 
(cf. Lemma~\ref{C1}(ii))}:

 For  an arbitrary
smooth affine variety $X$, and $ {\mathfrak{m}}_1 ,
{\mathfrak{m}}_2$, two distinct maximal ideals of $ A = \C[X],$ we have:
\begin{equation}
\label{Hom}
\mbox{\rm Hom}_{A}(A/\mathfrak{m}_1^k \,,\,
A/\mathfrak{m}_2) = 0  \quad \mbox{for all}\ k = 1, 2,...
\end{equation}

Now, fix $\, c\in \ccp \,$. By Lemma~\ref{L1}, we 
have $\, \tilde{X}_c \cong \h \,$. 
Hence we may (and will) identify $ \C[\tilde{X}_c] = \C[\h] \,$, 
and write for short $ A := \C[\h] $ and $ B := \C[X_c] \,$.
Thus, $\, \tilde{X}_c = \Spec\,A \,$, $\, X_c = \Spec\,B \,$, 
and the normalization map $\, \pi_c : \Spec\,A \to \Spec\,B \,$ is 
given by $ \mathfrak{m} \mapsto \mathfrak{m}\, \cap \, B \,$. 
Assume that $ \pi_c $ is {\it not} injective. Then there are 
(at least) two distinct maximal ideals $ \mathfrak{m}_1 $ and 
$\mathfrak{m}_2$ in $ A \,$, such that 
$\mathfrak{m}_1 \,\cap\, B = \mathfrak{m}_2 \,\cap \, B\,$.
Let $ \SH_c $ be the shift operator introduced above. It is easy 
to see from formula (\ref{shift}) (and Lemma~\ref{L4}) 
that $ \text{\it order}(\SH_c) = d $ and $ \SH_c \in \dd(A,B) \,$.
Hence $ \SH_c [\mathfrak{m}_1^k] \subset \mathfrak{m}_1 $ for all 
$ k > d \,$, and therefore
$\, \SH_c[\mathfrak{m}_1^k] \subseteq \mathfrak{m}_1 \,\cap\, B =
\mathfrak{m}_2 \,\cap\, B \subseteq \mathfrak{m}_2 \,$. It follows 
that, for each $ k > d \,$, $ \SH_c $ induces a differential operator
$ \bar{\SH}_c \in \dd(A/\mathfrak{m}_1^k \,,\, A/\mathfrak{m}_2) \,$.
But, if $\mathfrak{m}_1 \not= \mathfrak{m}_2 $, (\ref{Hom}) 
means that $ \dd_{0}(A/\mathfrak{m}_1^k \,,\, A/\mathfrak{m}_2)=0\,$
and then $ \dd(A/\mathfrak{m}_1^k \,,\, A/\mathfrak{m}_2)= 0 $
by the inductive definition of differential operators (see (\ref{1.1})).
Thus $ \bar{\SH}_c  = 0 \,$, which implies that
$ \SH_c[A] \subseteq \mathfrak{m}_2\,$. In particular, since $ 1 \in A\,$,
we get $\, \SH_c[1] \in \mathfrak{m}_2, $ with obvious
contradiction to Lemma~\ref{ng}. The result follows.
\qed\smallskip

Now, let $X$ be an affine algebraic variety with smooth normalization 
$\pi: \tilde{X}\to X$. Since $ \C[X] \subseteq \C[\tilde{X}]\,$, 
we may set
\begin{equation}
\label{5}
\dd(\tilde{X},X) := \{\,u \in \dd(\tilde{X})\ \mid 
 u\bigl(\C[\tilde{X}]\bigr) \subseteq \C[X]\,\} \ .
\end{equation}
Then $ \dd(\tilde{X},X) $ is a right ideal in $ \dd(\tilde{X}) \,$
and (with identifications (\ref{3})) a left ideal 
of $ \dd(X)\,$. According to \cite{SS}, Lemma~2.7, 
$\, \dd(\tilde{X},X) $  is isomorphic to the space 
$ \dd(\C[\tilde{X}],\C[X]) $ of all differential operators
from $ \C[\tilde{X}] $ to $ \C[X] $ over $ \C[X] \,$, so the 
notation (\ref{5}) is consistent with (\ref{1.1}).
We write $ \mbox{\rm End}_{\dd(\tilde{X})}\,
\dd(\tilde{X}, X) $ for the endomorphism ring of $ \dd(\tilde{X}, X) $
as a right $ \dd(\tilde{X})$-module and identify it with a subalgebra
in $ {\mathcal Q} $, the quotient skew-field of $ \dd(\tilde{X}) $
(cf. \S10 below).
\begin{theorem}
\label{end}
$\, \dd(X_c) = \mbox{\rm End}_{\dd(\tilde{X}_c)}\,
\dd(\tilde{X}_c, X_c) \,$.
\end{theorem}
\begin{proof} Clearly, $\dd(X_c) \subseteq 
\mbox{\rm End}_{\dd(\tilde{X}_c)}\,
\bigl(\dd(\tilde{X}_c, X_c)\bigr).\,$
Hence, we only need to show the opposite inclusion.
To this end, observe first that $\, \dd(\tilde{X}_c, X_c)
\bigl(\C[\tilde{X}_c]\bigr) = \C[X_c]\,$. Indeed, by definition
we have $\, \dd(\tilde{X}_c, X_c)\bigl(\C[\tilde{X}_c]\bigr)$
$ \subseteq \C[X_c]\,$. 
On the other hand, by Lemma~\ref{ng} and the fact that
$ \SH_c \in \dd(\tilde{X}_c, X_c) \,$, we get
$$ 
\C[X_c] = \C[X_c]\ccirc \,\SH_c[1] \subseteq 
\dd(\tilde{X}_c, X_c)[1] \subseteq \dd(\tilde{X}_c, X_c)
\bigl(\C[\tilde{X}_c]\bigr) \ .
$$ 
Now, write $\BK=\C(X_c) \cong \C(\h)$ for the field of rational 
functions on $X_c$. Then, \linebreak
$
\mbox{\rm End}_{\dd(\tilde{X}_c)}\bigl(\dd(\tilde{X}_c,
X_c)\bigr)\,\subset\,$ $\dd(\BK)\,$ by (\ref{6}).
Further, for any $ u \in \mbox{\rm End}_{\dd(\tilde{X}_c)}\,
\dd(\tilde{X}_c, X_c)$, we have
$$
u\bigl(\C[X_c]\bigr) = u\ccirc\,\dd(\tilde{X}_c,
X_c)\bigl(\C[\tilde{X}_c]\bigr)
\subseteq \dd(\tilde{X}_c, X_c)\bigl(\C[\tilde{X}_c]\bigr) =
\C[{X}_c]\,.
$$
Hence $ u \in \dd(X_c)\,$.
\end{proof}

Next, we recall a well-known ring-theoretic result usually 
referred to as the "Dual Basis Lemma". 
Let $ {A} $ be an algebra, $\, \ms{P} \,$ a right $A$-module,
and
$ B = \mbox{\rm End}_{{A}}\,\ms{P} ,$ the endomorphism ring
of $ \ms{P}, $ and $ \ms{P}^{*} = \mbox{\rm Hom}_{{A}}(\ms{P}, {A}) $
its dual module.
Write $\ms{P}\cdot \ms{P}^*$ for the
 subspace of $\End_{A}(\ms{P})$ spanned by the endomorphisms  
$p\cdot f: x \mapsto p\cdot f(x),$ for all 
$p \in \ms{P}$ and $f \in \ms{P}^* = 
\Hom_{A}(\ms{P},A). $
Obviously $\ms{P}\cdot \ms{P}^*$ is a two-sided ideal of
$B=\End_{A}(\ms{P})$.
Similarly, $\ms{P}^{*}\cdot  \ms{P}$  stands for a two-sided ideal of $A$ 
generated by the elements $\{f(q) \mid q \in \ms{P} \quad\text{and}\quad
f \in \ms{P}^*\}.$

\begin{lemma}[\cite{MR}, 3.5.2]
\label{DBL}
 The $A$-module $ \ms{P} $ is finitely
generated and projective  if and only if 
$\, \ms{P} \cdot \ms{P}^{*} = B\,$.\qed
\end{lemma}

Now we are in a position to prove one of the main results in the present
paper. 
\begin{theorem}
\label{Mor}
The ring $ \dd(X_c) $ is Morita equivalent to $ \dd(\h) \,$.
\end{theorem}

\begin{proof}
In the setup of Lemma \ref{DBL},
write $ A := \dd(\tilde{X}_c) \,$ , $\, B := \dd(X_c) \,$, and 
$ \ms{P} := \dd(\tilde{X}_c, X_c) \,$. Then, by
Theorem~\ref{end}, we have $ B = \mbox{\rm End}_{A}\,\ms{P} \,$.
Now, $
\,\ms{P} \cd\ms{P}^{*} \,$ is a nonzero 
two-sided ideal of $ B $, and therefore, 
is equal to $ B \,$ by Theorem~\ref{th_simple3}.
By Lemma~\ref{DBL}, we then conclude that
$\, \ms{P} $ is finitely generated projective.
On the other hand, $ \ms{P}^{*}\cd\ms{P} $ is a 
nonzero two-sided ideal of $ A \cong \dd(\h) $, 
and hence, by simplicity of the latter, is equal to $ A \,$. 
This means that $ \ms{P} $ is a generator
in the category of (right) $ A$-modules.
Summing up, we have $ B = 
\mbox{\rm End}_{A}\,\ms{P} \,$, where $ \ms{P} $
is a finitely generated projective generator in 
${\text{\sf{mod-}}}A \,$. By the Morita Theorem, $ A$ and
$ B $ are then equivalent rings. 
\end{proof}

The Morita equivalence of Theorem \ref{Mor} yields

\begin{corollary}\label{finH} The homological dimension
of the algebra $\dd(X_c)$ equals $\dim\h$.\qed
\end{corollary}

We end up this section with a conjecture giving an alternative 
description of the category of modules over $ \dd(X_c) $ hinted 
by \cite{DE}. 

Let $ i: X_c \into Y $ be a closed imbedding of $X_c$
into a {\it smooth} affine algebraic variety $Y$.
Write $\dd(X,Y) := \dd(\C[Y], \C[X_c])$ for the space of differential
operators
from $\C[Y]$ to $\C[X_c]$ over $\C[Y].$
Let ${\text{\sf{mod-}}}\dd(X_c)$ denote the abelian
 category of right
$\dd(X_c)$-modules, and ${\text{\sf{mod}}}_{_{X_c}}\!\text{-}\dd(Y)$
 denote the abelian
 category of right
$\dd(Y)$-modules with support in  $i(X_c) \subset Y.$


\begin{conjecture} 
\label{Kash}
The direct image functor $i_{+}$ is an equivalence of categories, with quasi-inverse
$i^+: {\text{\sf{mod}}}_{_{X_c}}\!\text{-}\dd(Y)\too
{\text{\sf{mod-}}}\dd(X_c)
\,,\,N\,\longmapsto\, \Hom_{\dd(Y)}(\dd(X,Y)\,,\, N).$
\end{conjecture}
Conjecture~\ref{Kash} can be viewed as a generalization of a 
well-known result of Kashiwara to the {\it singular} variety $X_c\,$. 
In case when $ \dim X_c = 1 $ (or more generally, for any singular 
affine curve with injective normalization) this conjecture has 
been proved in \cite{DE}. 

Conjecture \ref{Kash} has been proved in full generality in \cite{BN}.

\section{Appendix:\, A  filtration
on differential operators}\setcounter{equation}{0}
The goal of this section is to put the filtration
(\ref{filt_flat}) into a more general context
of the ring of differential operators on a (singular) algebraic
variety $X$ whose normalisation is a vector space,
cf. also [P] in the special case $\dim X=1$.

Let $X$ be any irreducible (not necessarily smooth) algebraic variety,
$\pi: \tilde{X}\to X$ its normalisation,
and $\BK=\C(X)=\C(\tilde{X}),$ the field of rational
functions.
The ring $ \dd(\BK) $ is a Noetherian domain (see, e.g., 
\cite{MR}, 15.5.5). By Goldie's Theorem, it has 
the quotient skew-field which we denote by $ \qq \,$. 
Let $\ms{P}:= \dd(\tilde{X},X),$
regarded as 
a right ideal of $ \dd(\tilde{X}), \,$ see (\ref{5}).
We put  $ \ms{E} := 
\mbox{End}_{\dd(\tilde{X})}\,\ms{P} ,$  the endomorphism 
ring of $ \ms{P}. \,$  Then
we may (and will) identify $\, \ms{E} 
\cong \{\, q \in \qq \ : \ q\,\ms{P} \subseteq \ms{P}\,\}\,$
(see \cite{MR}, 3.1.15).
Since $\, \ms{P} $ is also a left ideal of $ \dd(X) \,$, 
we have $\, \dd(X) \subseteq \ms{E} \,$. On the other hand, 
$\, \C[\tilde{X}] $ being finite over 
$\, \C[X] $ implies that $ f\cdot \C[\tilde{X}] \subseteq  \C[X] $ 
for some nonzero $ f \in  \C[\tilde{X}] \,$, 
hence $ \ms{P} \cap \C[\tilde{X}] \not= \{0\}\,$
and therefore $ \ms{P} \in \dd(\BK)\ \Rightarrow \  \ms{E} \subset \dd(\BK) \,$.
We summarize these inclusions in the following
\begin{lemma}[see \cite{SS}]
\label{6}
$\quad\ms{P} \,\subset\, \dd(X) \,\subseteq \, \ms{E}\,
\subset \,\dd(\BK) \quad , \quad 
\ms{P}\, \subset \,\dd(\tilde{X}) \,\subset \,\dd(\BK)\ .
$\qed
\end{lemma}

By definition, the algebra $ \dd(\BK) $ carries
a natural (differential) filtration. Using (\ref{6}),
we may (and will) equip the spaces $ \ms{P},\, 
\ms{E},\, \dd(X), \,\dd(\tilde{X}) $ with the induced 
filtrations\footnote{The induced filtrations 
on $ \dd(X) $ and $ \dd(\tilde{X}) $ coincide with the 
intrinsic ones under  identifications (\ref{3}).}.
We write $ \grd\,\ms{P}, \, \grd\,\ms{E}, \, \grd\,\dd(X), \,
\grd\, \dd(\tilde{X}) $ for the associated graded objects. 
The next result is probably known to the experts, 
but we cannot find a precise reference in the literature.
\begin{proposition}
\label{P1}
Let $ X $ be an irreducible affine algebraic variety over 
$ \C\,$. Assume that the normalization $ \tilde{X} $ is smooth. Then
\begin{equation}
\label{7}
\grd\, \ms{P} \,\subset\, \grd\,\dd(X) \,\subseteq \,
\grd\, \ms{E}\,\subseteq \,\grd \,\dd(\tilde{X}) \,\subset \, \grd\,\dd(\BK)\ .
\end{equation}
\end{proposition}
\begin{proof}
Given the inclusions (\ref{6}), we need only to prove
$\, \grd\, \ms{E}\, \subseteq \,\grd \,\dd(\tilde{X}) \,$.
First, recall that $ \tilde{X} $ being smooth implies 
that $ \grd \,\dd(\tilde{X}) \cong \C[T^{*}\tilde{X}] \,$ 
is a regular Noetherian commutative domain. Hence
$ \grd\,\ms{P} \subseteq \grd \,\dd(\tilde{X}) $ is a 
finitely generated ideal of $ \grd \,\dd(\tilde{X}) \,$. 
Since $ \C[\tilde{X}] $ is finite over $ \C[X] \,$,
$\, \ms{P} $ is nonzero: in fact, $\, \ms{P} \supset 
\mbox{Ann}_{\C[X]}(\C[\tilde{X}]/\C[X]) \not= \{0\}\,$.
Next, with our description of the endomorphism ring $ \ms{E} \,$, we have 
$\, \ms{E}\cdot  \ms{P} \subseteq \ms{P}\,$ in $\, \dd(\BK) \,$, and therefore
$\, (\grd\,\ms{E}) \cdot (\grd\,\ms{P}) \,\subseteq \, \grd\, \ms{P} \,$ in 
$ \grd\,\dd(\BK) \,$. Hence, by \cite{AM}, Prop.~2.4, the 
elements of $ \grd\, \ms{E} $ are integral over $ \grd \,
\dd(\tilde{X}) $ in
$ \grd\,\dd(\BK) \,$. But, being regular, the ring $ \grd \,
\dd(\tilde{X}) $ 
is integrally closed, hence $ \grd\,\ms{E} \subseteq \,\grd 
\,\dd(\tilde{X}) \,$.
\end{proof}

\begin{remark}
\label{R1}
If $ \dim\,X = 1 $ then $ \tilde{X} $ is always smooth. 
In this case our Proposition~\ref{P1} implies 
Proposition~3.11 of \cite{SS} which establishes the 
inclusion $ \grd\,\dd(X) \,\subseteq \,\grd \,\dd(\tilde{X}) $
for any (irreducible affine) curve. 
\end{remark}

Now assume that $ X $ is an irreducible variety 
with normalization $ \tilde{X} $ being isomorphic to $ V \,$,
a vector space. We identify $\C(X)$ with $\BK=\C(V)$.
On $\BK$, introduce an increasing filtration
$\{F\ff_j\BK\}_{j\in \Z}$ as follows.
Let $\PPP(V\oplus\C)\supset V$ be the projective
completion of $V$, and $D_\infty=\PPP(V\oplus\C)
\smallsetminus V$, the divisor at infinity.
Given $f\in \BK=\C(V)$, we can think of $f$ as a rational function
on $\PPP(V\oplus\C)$, and we put $f\in F\ff_j\BK$ if
$f$ has a pole of order $\leq j$ at $D_\infty$.
The restriction of this filtration on $\BK$ to $\C[V]$
coincides with the standard filtration by degree
of polynomial. We extend the filtration
$F\ff_\bullet \BK$ to the algebra $\dd(\BK)$ of differential
operators with rational coeffitients by
assigning filtration degree zero to all constant
coefficient differential operators.
 We write $\,\{\dd_{n}^{\flat}(\BK)\}_{n \in \Z} \,$ 
for the induced filtration on $ \dd(\BK) $ and call it 
the $\flat$-{\it filtration}. 
Now, if $ E$ is a linear subspace of $ \dd(\BK) \,$, we set 
$ E_{n}^{\flat} := E \cap \dd_{n}^{\flat}(\BK) $ and write
$ \grd^{\flat}\,E := \bigoplus_{ n\in \Z} 
E_{n}^{\flat}/E_{n-1}^{\flat} $ for the associated graded space.
Clearly, by symmetry with the standard filtration, we have 
$ \grd^{\flat}\,\dd(\tilde{X}) \cong \C[V \times V^*] \,$, 
and therefore 
$ \grd^{\flat}\,\dd(\BK) $ is isomorphic to a (properly graded) 
localization of $ \C[V \times V^*] $ at $  \C[V] \setminus \{0\}\,$.
In particular, both $ \grd^{\flat}\,\dd(\tilde{X}) $ and
$ \grd^{\flat}\,\dd(\BK) $ are commutative domains, the former being
regular. Reviewing the proof of Proposition~\ref{P1} then shows 
that its argument works for the $ \flat$-filtration as well. 
Thus, with the above notation and conventions, we have the following
\begin{proposition}
\label{P2}
If $ X $ is an irreducible affine variety with normalization 
$ \tilde{X} $ being isomorphic to $ V\,$, then
$\,\grd^{\flat}\, \ms{P} \,\subset\, \grd^{\flat}\,\dd(X) \,\subseteq \,
\grd^{\flat}\,\ms{E}\,\subseteq \,\grd^{\flat} \,\dd(\tilde{X}) \,\subset \, 
\grd^{\flat}\,\dd(\BK)\,$. 
\end{proposition}

By definition, the differential operators of order zero
form a maximal commutative, (locally) ad-nilpotent subalgebra
of $ \dd(X) \,$, which is naturally isomorphic to $ \C[X]\,$. 
On the other hand,
from  Proposition~\ref{P2} we obtain the following generalization of
the result (proved in \S6) on the nilpotency
of the adjoint action on $\ddx$ of the subalgebra
$Q_c\ff=\dd\ff_0(X_c)$.
\begin{proposition}
\label{Co1} If $ \tilde{X} \cong V \,$, then
$\, \dd_{0}^{\flat}(X)$ is a commutative (locally) ad-nilpotent 
subalgebra of $ \dd(X) \,$. 
\end{proposition}
\begin{proof}
Since $ \grd^{\flat}\,\dd(X) \subset  \C[V \times V^*] \,$,
the elements of $ \dd(X) $ have non-negative $ \flat$-degrees.
Thus $\, \dd_{-1}^{\flat}(X) = \{0\} \,$.
In view of commutativity of $ \grd^{\flat}\,\dd(\BK) \,$, we have
$ [\dd_{0}^{\flat}(X)\,,\,  \dd_{n}^{\flat}(X)] 
\subset \dd_{n-1}^{\flat}(X)\,$ for all $ n \geq 0\,$. 
Hence $ \dd_{0}^{\flat}(X) $ is ad-nilpotent. 
Moreover, if $ n = 0 $, we have
$ [\dd_{0}^{\flat}(X)\,,\,  \dd_{0}^{\flat}(X)] 
\subseteq \dd_{-1}^{\flat}(X) =  \{0\} \,$, therefore 
 $ \dd_{0}^{\flat}(X) \,$ is commutative.
\end{proof}

\footnotesize{

\footnotesize{
{\bf Y.B.}: Department of Mathematics, Cornell University,
Ithaca, NY 14853-4201, USA;\\
\hphantom{x}\quad\, {\tt berest@math.cornell.edu}

{\bf P.E.}: Department of Mathematics, Rm 2-165, MIT,
77 Mass. Ave, Cambridge, MA 02139;\\
\hphantom{x}\quad\, {\tt etingof@math.mit.edu}

{\bf V.G.}: Department of Mathematics, University of Chicago,
Chicago, IL
60637, USA;\\
\hphantom{x}\quad\, {\tt ginzburg@math.uchicago.edu}}


\begin{thebibliography}{APK1}

\bibitem[AM]{AM}
M. F. Atiyah and I. G. Macdonald, \textit{Introduction to Commutative 
Algebra}, Addison-Wesley Publishing Co., Reading, 1969.
%

\bibitem[BW]{BW}
Yu. Berest and G. Wilson,
\textit{Classification of rings of differential
operators on affine curves}, Internat. Math. Res. Notices  
\textbf{2} (1999), 105--109.
%

\bibitem[Be]{Be}
Yu. Berest, \textit{Huygens' principle and the bispectral problem}
in \textit{The Bispectral Problem}, CRM Proceedings and Lecture 
Notes, \textbf{14}, Amer. Math. Soc. 1998, pp. 11--30.

\bibitem[BEG1]{BEG1}
Y. Berest, P. Etingof, and V. Ginzburg, Finite-dimensional 
representations of rational Cherednik algebras, 
Int. Math. Res. Not. (2003), no. 19, pp. 1053-1088.

\bibitem[BEG2]{BEG2}
Y. Berest, P. Etingof, and V. Ginzburg, Morita 
equivalence of Cherednik algebras, 
J. Reine Angew. Math., 568 (2004), 81–98.

\bibitem[BGG]{BGG} J. Bernstein,
I. Gelfand and S. Gelfand,
\textit{A certain category of
${\mathfrak g}$-modules.} (Russian) Funct. Anal. Appl. 10:2 (1976),  1--8.

\bibitem[BGG1]{BGG1}
J. Bernstein, I. Gelfand and S. Gelfand, 
\textit{Differential operators on the cubic cone},
Russian Math. Surveys \textbf{27} (1972), 466--488.


\bibitem[B]{B}
J. E. Bj\"ork, \textit{Rings of Differential Operators},
North-Holland Mathematical Library \textbf{21}, 
North-Holland Publishing Co., Amsterdam, 1979. 
%
%
\bibitem[BMR]{BMR} 
M. Brou\'e, G. Malle and  R. Rouquier, \textit{
Complex reflection groups, braid groups, Hecke algebras},
J. Reine Angew. Math. \textbf{500} (1998), 127--190.
%
\bibitem[BN]{BN} D. Ben-Zvi, T. Nevins, \textit{Cusps and $\mathcal D$-modules.}
  J. Amer. Math. Soc.  17  (2004),   155--179.
\bibitem[C]{C}
O. A. Chalykh, \textit{Darboux transformations for multidimensional
Schr\"odinger operators}, Russian Math. Surveys \textbf{53}(2) 
(1998), 377--379.   



\bibitem[CV]{CV}
O.~A.~Chalykh and A.~P.~Veselov, \textit{
Commutative rings of partial differential operators and Lie algebras}, 
Comm. Math. Phys. \textbf{126}(3) (1990), 597--611. 
%
\bibitem[CV2]{CV2} 
O.~A.~Chalykh and A.~P.~Veselov, \textit{Integrability in the 
theory of Schr\"odinger operator and harmonic analysis},
Comm. Math. Phys. \textbf{152}(1) (1993), 29--40. 
%
\bibitem[CS]{CS}
M. Chamarie and J. T. Stafford, \textit{When rings 
of differential operators are maximal orders},
Math. Proc. Camb. Phil. Soc. \textbf{102} 
(1987), 399--410.
%
\bibitem[Ch]{Ch} I. Cherednik, \textit{Lectures on 
affine KZ equations, quantum many-body problems, 
Hecke algebras, and Macdonald theory}, Lecture Notes 
by Date et al., Kyoto, 1995. 
 
\bibitem[Ch2]{Ch2} I. Cherednik, {\it Double affine Hecke
 algebras,
Knizhnik-Zamolodchikov equations, and Macdonald operators,}
IMRN (Duke math.J.) {\bf 9} (1992), p.171-180.



\bibitem[CG]{CG} N. Chriss, V. Ginzburg, {\it
Representation theory and complex geometry.} 
Birkh\"auser Boston, 1997.


\bibitem[Dz]{Dz}  C. Dezelee, \textit{Representations 
de dimension finie de l'algebre de Cherednik rationnelle.}\hfill\break
    {\tt{arXiv:math.RT/0111210}}.

\bibitem[DE]{DE} M. Van Doorn and A. Van den
Essen, \textit{$\dd$-Modules with support on a curve},
 Publ. RIMS \textbf{23} (1987), 
937-953.
\bibitem[Di]{Di} J. Dixmier, \textit{
Enveloping algebras.}
Graduate
Studies in Mathem. \textbf{11}, A.M.S., Providence, RI, 1996. 

\bibitem[D]{D}
C. F. Dunkl, \textit{Differential-difference operators and 
monodromy representations of Hecke algebras}, 
Pacific J. Math. \textbf{159}(2) (1993), 271--298.

\bibitem[DO]{DO} C. F. Dunkl and E. Opdam,  
\textit{Dunkl operators for complex reflection groups}, 
 {\tt{arXiv:math.RT/0108185}}.
\bibitem[DJO]{DJO} 
C. F. Dunkl,  M. F. E. de Jeu and E. Opdam,  
\textit{Singular polynomials for finite reflection groups}, 
Trans. Amer. Math. Soc. \textbf{346} (1994), 237--256.
%
\bibitem[EG]{EG}  
P. Etingof and V. Ginzburg, \textit{Symplectic 
reflection algebras, Calogero-Moser space, and 
deformed Harish-Chandra homomorphism}.
Invent. Math. {\bf 147} (2002), 243-348,
[ {\tt{arXiv:math.AG/0011114}} ].
%
\bibitem[EG2]{EG2} P. Etingof and V. Ginzburg, 
\textit{On $m$-quasi-invariants of Coxeter groups},
Preprint {\tt{arXiv:math.QA/0106175}}, to appear in
Moscow Mathem. Journ.
%
\bibitem[FV]{FV} 
M. Feigin and A. Veselov, \textit{Quasi-invariants of
Coxeter groups and $m$-harmonic polynomials}, 
Int. Math. Res. Not.
{\bf 10} (2002), 521--545. [{\tt{arXiv:math-ph/0105014}}].
%
\bibitem[FeV]{FeV} 
G. Felder and A. Veselov, 
\textit{Action of Coxeter groups on $m$-harmonic polynomials
and KZ equations}, Preprint 2001, {\tt{arXiv:QA/0108012}}.

\bibitem[Gi]{Gi} V. Ginzburg, \textit{On primitive ideals},
{\tt{arXiv:math.RT/0202079}}.

\bibitem[Go]{Go} I. Gordon, {\it 
 Baby Verma modules for rational Cherednik algebras.}
{\tt{arXiv:math.RT/0202301}}.
%
\bibitem[G]{G}
A. Grothendieck, \textit{El\'ements de G\'eom\'etrie 
Alg\`ebrique IV}, Publ. Math. \textbf{32}, IHES, Paris,
1967.
%
\bibitem[Gu]{Gu} N. Guay,
\textit{Projectives in the category $\oo$ for the Cherednik algebra.}
 Preprint, University of Chicago 2001.

\bibitem[GGOR]{GGOR}
V. Ginzburg, N. Guay, E. Opdam, R. Rouquier, On the category  
$ {\mathcal O}$ for rational Cherednik algebras, 
Invent. Math. 154 (2003), no. 3, 617-651.

\bibitem[He]{He} G. Heckman, \textit{
A Remark on the Dunkl differential-difference operators.} Harmonic
analysis on reductive groups, 181--191, Progr. Math., \textbf{101},
 Birkh\"auser Boston, Boston,
MA, 1991. 

\bibitem[Ho]{Ho}
T. Hodges, \textit{Morita equivalence of primitive factors of $U(\sll2)$,}
Contemp. Math. \textbf{139} (1992), 175-179.

\bibitem[HS]{HS}
R. Hart and S. P. Smith, \textit{Differential operators on 
some singular surfaces}, Bull. London Math. Soc. 
\textbf{19} (1987), 145--148.
%

\bibitem[K]{K} V. Kac, \textit{Constructing groups associated to
infinite
dimensional
Lie algebras,} in: Infinite dimensional groups and appl., MSRI
Publ. \textbf{4}, Springer-Verlag, 1985, 167-216.

\bibitem[Ku]{Ku} A. Kuznetsov, \textit{Quiver varieties and Hilbert
schemes.} {\tt{arXiv:math.AG/0111092}}

%


\bibitem[LS]{LS} T. Levasseur and J. Stafford,
\textit{Invariant differential operators and a 
homomorphism of Harish-Chandra}, 
J. Amer. Math. Soc. \textbf{8} (1995),  365--372.




\bibitem[LS2]{LS2} T. Levasseur and J.T. Stafford,
{\it Semi-simplicity of invariant
holonomic systems on a reductive Lie algebra}.
 Amer. J. Math. {\bf 119}
(1997),  1095--1117.


\bibitem[Ma]{Ma} I.G.Macdonald, {\it Symmetric functions and Hall polynomials,}
second edition,  Oxford University Press (1995).

\bibitem[MR]{MR}
J. C. McConnell and J. C. Robson,
\textit{Noncommutative Noetherian Rings}, J.Wiley \& Sons, NY 1987.

\bibitem[Mo]{Mo} S. Montgomery,
{\it Fixed Rings of Finite Automorphism Groups of
Associative Rings}, LN in Math., {\bf 818}
Springer-Verlag, Berlin/New-York, 1980.


\bibitem[Mu]{Mu}
I. M. Musson, \textit{Some rings of differential operators which 
are Morita equivalent to the Weyl algebra}, Proc. Amer. Math. Soc. 
\textbf{98} (1986), 29--30.

\bibitem[MvdB]{MvdB} I. M. Musson, M. Van den Bergh, \textit{Invariants under tori of 
rings of differential operators and related topics.}
Mem. Amer. Math. Soc. \textbf{136} (1998), no. 650. 
%

\bibitem[O]{O} E. Opdam, 
\textit{Dunkl operators, Bessel functions and the discriminant
of a finite Coxeter group}, Compositio Math. \textbf{85}(3) 
(1993), 333--373. 

\bibitem[O1]{O1}
E. M. Opdam, A remark on the irreducible characters 
and fake degrees of finite real 
reflection groups, Invent. Math. 120 (3) (1995), 447--454.

\bibitem[OR]{OR}
E. Opdam and R. Rouquier, private communication, July 2001. 
%

\bibitem[P]{P} P. Perkins, \textit{Commutative subalgebras of the ring of 
differential operators on a curve}, Pacific Journal of Math.
\textbf{139} (1989), 279--302.


\bibitem[Sm]{Sm}
S. P. Smith, \textit{An example of a ring Morita
equivalent to the Weyl algebra}, J. Algebra
\textbf{73} (1981), 552--555.
%
\bibitem[SS]{SS}
S. P. Smith and J. T. Stafford, \textit{Differential
operators on an affine curve}, Proc. London Math. Soc.
(3) \textbf{56} (1988), 229--259.
%

\bibitem[St]{St}
J. T. Stafford, \textit{Homological properties of the 
enveloping algebra $U(\sll2)$}, Math. Proc. Cambr. Phil. Soc. 
\textbf{91} (1982), 29--37.



\bibitem[S]{S}
R. Stanley, \textit{Invariants of finite groups and 
their applications to combinatorics}, 
Bull. Amer. Math. Soc. (N.S.) \textbf{1}(3) (1979),  475--511.

\bibitem[S2]{S2} R. Stanley, \textit{Enumerative combinatorics.} {\sl {Vol. 2}}.
 Cambridge Studies in Advanced Mathematics,  \textbf{62}
 Cambridge University Press, Cambridge, 1999.

\bibitem[VdB]{VdB}
M. Van den Bergh, \textit{Differential operators on semi-invariants 
for tori and weighted projective spaces} in 
\textit{Topics in Invariant Theory}, Lecture Notes in Math. 
\textbf{1478}, Springer, Berlin, 1991, pp. 255--272. 
%
\bibitem[VSC]{VSC}
 A. P. Veselov, K. L. Styrkas and O. A. Chalykh, 
\textit{Algebraic integrability for the Schr\"odinger equation, 
and groups generated by reflections},
Theoret. and  Math. Phys. \textbf{94}(2) (1993), 182--197.

\bibitem[Wa]{Wa} N. Wallach,
\textit{Invariant differential operators on a reductive
                Lie algebra and Weyl group representations},
                J. Amer. Math. Soc. \textbf{6} (1993), 779--816.


\bibitem[Wi]{Wi} G.Wilson, \textit{Collisions of Calogero-Moser particles and an
adelic Grassmannian,} Invent. Mathem. {\bf 133} (1998), 1-41.

%
\end{thebibliography}
\end{document}